\documentclass[12pt]{article}
\usepackage{epsfig}

\def\be{\begin{equation}}
\def\ee{\end{equation}}
\def\bea{\begin{eqnarray}}
\def\eea{\end{eqnarray}}
\def\bes{\begin{eqnarray*}}
\def\ees{\end{eqnarray*}}

\def\nn{\nonumber}
\def\<{\langle}
\def\>{\rangle}
\def\lb{\label}
\def\bs{\setminus}

\def\R{{\bf R}}
\def\C{{\bf C}}
\def\Z{{\bf Z}}
\def\N{{\bf N}}
\def\U{{\bf U}}
\def\Q{{\bf Q}}
\def\T{{\bf T}}

\def\ga{{\gamma}}

\def\th{{\theta}}
\def\om{{\omega}}
\def\Om{{\Omega}}
\def\ep{{\epsilon}}
\def\lm{{\lambda}}

\def\sg{{\sigma}}

\def\Sg{{\Sigma}}
\def\vf{{\varphi}}

\def\H{{\cal H}}
\def\T{{\cal T}}

\def\P{{\cal P}}

\def\Nn{{\cal N}}

\def\rank{{\rm rank}}
\def\Sp{{\rm Sp}}
\def\mod{{\rm mod}}
\def\dm{{\rm \diamond}}

\def\ol#1{\overline{#1}}  

\def\hb{\vrule height0.18cm width0.14cm $\,$}

\def\ol#1{\overline{#1}}  

\title{Closed characteristics on compact convex  \\hypersurfaces in $\R^8$}
\author{Wei Wang\thanks{Partially supported by National Natural Science Foundation of China No. 11222105,  Foundation for the Author of National Excellent Doctoral Dissertation of PR China No. 201017.
E-mail: alexanderweiwang@gmail.com, wangwei@math.pku.edu.cn  }\\
Key Laboratory of Pure and Applied Mathematics\\
School of Mathematical Science \\ Peking University, Beijing 100871 \\
PEOPLES REPUBLIC OF CHINA \\ }

\date{ Nov.  26th, 2013}

\begin{document}

\maketitle

\begin{abstract}
{\it  In this paper, we prove there exist at least four geometrically distinct
closed characteristics on every compact convex hypersurface $\Sg$ in $\R^8$.
This gives a confirmed answer in the case $n=4$ to a long standing conjecture
in Hamiltonian analysis since the time of A. M. Liapounov in 1892 (cf. P. 235 of \cite{Eke3}).
}
\end{abstract}



{\bf Running title}: Closed characteristics on convex
hypersurfaces 

\renewcommand{\theequation}{\thesection.\arabic{equation}}
\renewcommand{\thefigure}{\thesection.\arabic{figure}}

\setcounter{equation}{0}
\section{Introduction and main results}

In this paper, let $\Sigma$ be a  $C^3$ compact convex hypersurface in
$\R^{2n}$, i.e., $\Sigma$ is the boundary of a compact and strictly
convex region $U$ in $\R^{2n}$. We denote the set of all such
hypersurfaces by $\H(2n)$. Without loss of generality, we suppose
that $U$ contains the origin. We consider closed characteristics $(\tau,
y)$ on $\Sigma$, which are solutions of the following problem \be
\left\{\matrix{\dot{y}=JN_\Sigma(y), \cr
               y(\tau)=y(0), \cr }\right. \lb{1.1}\ee
where $J=\left(\matrix{0 &-I_n\cr
                I_n  & 0\cr}\right)$ is the standard symplectic matrix in $\R^{2n}$,
$I_n$ is the identity matrix in $\R^n$, $\tau>0$ is the period of $y$, $N_\Sigma(y)$ is
the outward normal vector of $\Sigma$ at $y$ normalized by the
condition $N_\Sigma(y)\cdot y=1$. Here $a\cdot b$ denotes the
standard inner product of $a, b\in\R^{2n}$.
A closed characteristic
$(\tau, y)$ is {\it prime}, if $\tau$ is the minimal period of $y$.
Two closed characteristics $(\tau, y)$ and $(\sigma, z)$ are {\it
geometrically distinct},  if $y(\R)\not= z(\R)$. We denote by
$\T(\Sg)$ the set of all geometrically distinct closed
characteristics on $\Sg$. A closed characteristic $(\tau,y)$ is {\it
non-degenerate}, if $1$ is a Floquet multiplier of $y$ of precisely
algebraic multiplicity $2$, and is {\it elliptic}, if all the
Floquet multipliers of $y$ locate on ${\bf U}=\{z\in\C\,|\,|z|=1\}$,
i.e., the unit circle in the complex plane. It is  {\it hyperbolic},
if $1$ is a double Floquet multiplier of it and all the other
Floquet multipliers of $y$ are away from  ${\bf
U}$.

It is surprising enough that A. M. Liapounov in \cite{Lia1} of 1892
and J. Horn in \cite{Hor1} of 1903 were able to prove the following
great result: {\it Suppose $H: \R^{2n}\rightarrow\R$ is analytic,
$\sigma(JH^{\prime\prime}(0))=\{\pm\sqrt{-1}\omega_1, \ldots
,\pm\sqrt{-1}\omega_n\}$ are purly imaginary and satisfy
$\frac{\omega_i}{\omega_j}\notin\Z$ for all $i, j$. Then there
exists $\epsilon_0>0$  small enough such that }
 \be
^{\#}\T(H^{-1}(\epsilon))\ge n, \qquad \forall \; 0< \epsilon\le
\epsilon_0.\lb{1.2}\ee
This deep result was greatly improved by A.
Weinstein in \cite{Wei1} of 1973. He was able to prove that for
$H\in C^2(\R^{2n}, \R)$, if $H^{\prime\prime}(0)$ is positive
definite, then there exists $\epsilon_0>0$ small such that
(\ref{1.2}) still holds. In \cite{EL}, I. Ekeland and  J. Lasry
proved that if there exists $x_0\in\R^{2n}$ such that
\bea r\le |x-x_0|\le R,\qquad \forall x\in\Sg
\nn\eea
and $\frac{R}{r}<\sqrt{2}$, then $^\#\T(\Sg)\ge n$.

Note that we have the following example of weakly non-resonant
ellipsoid: Let $r=(r_1,\ldots, r_n)$ with $r_i>0$ for $1\le i\le n$.
Define \bea \mathcal{E}_n(r)=\left\{z=(x_1, \ldots,x_n,
y_1,\ldots,y_n)\in\R^{2n}\left |\frac{}{}\right.
\frac{1}{2}\sum_{i=1}^n\frac{x_i^2+y_i^2}{r_i^2}=1\right\}\nn\eea¡¡
where $\frac{r_i}{r_j}\notin\Q$ whenever $i\neq j$. In this case,
the corresponding Hamiltonian system is linear and all the solutions of (\ref{1.1})
can be computed explicitly. Thus it is easy to verify that
$^{\#}\T(\mathcal{E}_n(r)) = n$ and all the closed characteristics
on $\mathcal{E}_n(r)$ are elliptic and non-degenerate, i.e.,  its
linearized Poincar\'e map splits into $n-1$ two dimensional rotation
matrix $\left(\matrix{\cos\theta &\sin\theta\cr -\sin\theta & \cos\theta\cr}\right)$
with $\frac{\theta}{\pi}\notin\Q$ and one $\left(\matrix{1 &1\cr 0 & 1\cr}\right)$
in appropriate coordinates.

Based on these facts, there is a long standing conjecture on the
number of closed characteristics on compact convex hypersurfaces in
$\R^{2n}$: \be \,^{\#}\T(\Sg)\ge n, \qquad \forall \; \Sg\in\H(2n).
\lb{1.3}\ee
Since the pioneering works \cite{Rab1} of P. Rabinowitz and
\cite{Wei2} of A. Weinstein in 1978 on the existence of at least one
closed characteristic on every hypersurface in $\H(2n)$, the
existence of multiple closed characteristics on $\Sg\in\H(2n)$ has
been deeply studied by many mathematicians.
When $n\ge 2$,  in 1987-1988, I. Ekeland-L.
Lassoued, I. Ekeland-H. Hofer, and A, Szulkin (cf. \cite{EkL1},
\cite{EkH1}, \cite{Szu1}) proved
$$ \,^{\#}\T(\Sg)\ge 2, \qquad \forall\,\Sg\in\H(2n). $$
 In \cite{HWZ1} of 1998, H. Hofer-K. Wysocki-E. Zehnder
proved that $\,^{\#}\T(\Sg)=2$ or $\infty$ holds for every
$\Sg\in\H(4)$. In \cite{LoZ1} of 2002, Y. Long and C. Zhu further proved
\bea\;^{\#}\T(\Sg)\ge \left[\frac{n}{2}\right]+1, \qquad \forall\, \Sg\in \H(2n), \lb{1.4}\eea
where we denote by $[a]\equiv\max\{k\in\Z\,|\,k\le a\}$.
 In \cite{WHL} of 2007,
W. Wang, X. Hu and Y. Long proved $\,^{\#}\T(\Sg)\ge 3$ for
every $\Sg\in\H(6)$, which gave a confirmed answer to the above
conjecture in the case $n=3$.

The following main result of this paper gives a confirmed answer to the above
conjecture for the case  $n=4$.

{\bf Theorem 1.1.} {\it   There exist at least four geometrically distinct
closed characteristics on every compact convex hypersurface $\Sg$  in $\R^{8}$, i.e., we have
$^{\#}\T(\Sg)\ge 4$ for any $\Sg\in\H(8)$.}

The proof of Theorem 1.1 is given in Section 5. Mainly ingredients
in the proof include: the critical point theory for closed
characteristics established in \cite{WHL}, Morse theory,
the index iteration theory developed by Long and his coworkers,
a new method to handle the  degenerate critical point and Kronecker's uniform distribution theorem in number theory.

In contrast  to the previous works, we introduce several new ideas in this paper.
In fact, by Theorem 1.1 of \cite{LoZ1},  the lower bound $\varrho_4(\Sg)$  for 
$^{\#}\T(\Sg)$ is $3$ if there exists a 
closed characteristic $(\tau, y)$ on $\Sg\in\H(8)$ satisfying  $i(y,\,1)=5$ together with $\gamma_{y}(\tau)$ can be connected within
$\Omega^0(\gamma_{y}(\tau))$ to $N_1(1, 1)\diamond
N_1(1,-1)^{\diamond3}$ (cf. Case B in \S4 and \S3 for notations) or $i(y,\,1)=4$ together with  $\gamma_{y}(\tau)$
can be connected within $\Omega^0(\gamma_{y}(\tau))$ to $N_1(1,
1)\diamond N_1(1,-1)^{\diamond2}\diamond M^\prime$ for some
$M^\prime\in\{R(\theta), \,D(\lambda),\,N_1(-1, b),\, I_2\}\subset\Sp(2)$
(cf. Case A in \S4).
Hence we must develop new methods to overcome the difficulties caused by these cases:

(i) We find that the orders of appropriate iterations of any two fixed prime closed characteristics in the common index jump 
intervals have certain commutative property   (cf. Proposition 4.5), i.e., given any two prime closed characteristics, there must exist two common index jump 
intervals such that  the orders of appropriate iterations of these two closed characteristics in these two   intervals interchange. 

(ii) The critical modules for iterations of closed characteristics have the periodic property (cf. Proposition 2.6 below).

(iii) By (i) and (ii), we can firstly derive some stability properties for closed characteristics.
In fact, in order to interchange the orders of two  closed characteristics in the common index jump intervals, their linearized Poincar\'e map must have enough numbers of  components of
rotation matrix with irrational angles.
 
(iv) Then we obtain the desired multiplicity  result by a combination of Morse theory,
index iteration theory and Kronecker's uniform distribution theorem in number theory. 

(v) In this paper, the main idea to prove Theorem 1.1 is  studying the relations
between the closed characteristics, i.e., the closed characteristics are dependent.
  While the  methods of Y. Long et al.
concerns firstly multiplicity, then the stability;  their methods view the closed characteristics
as independently.  

These  viewpoints are new and used firstly in this paper to handle the multiplicity problem.

Here we  give the outline of the proof of Theorem 1.1.
By Theorem 1.1 of \cite{LoZ1}, we have $^\#\T(\Sg)\ge 3$ for every $\Sg\in\H(8)$.
We prove Theorem 1.1 by contradiction, i.e., assume $^\#\T(\Sg)=3$ for some $\Sg\in\H(8)$.
Applying the Fadell-Rabinowitz index theory to the Clarke-Ekeland dual action
functional $\Phi$ (cf. (\ref{2.24})), we obtain a sequence of critical values 
$$-\infty<c_1<c_2<\dots<c_k<c_{k+1}<\dots<0$$
of $\Phi$.
Critical points of $\Phi$ correspond exactly to closed characteristics  on $\Sg$.
Since $\Phi$ is not defined on a Hilbert space, in order to apply Morse theory,
we construct a functional $\Psi_a$ (cf. (\ref{2.3})) which have isomorphic
critical modules as $\Phi$ at the corresponding critical points,
while the critical modules of $\Psi_a$ can be computed out via
Gromoll-Meyer theory. Thus there exists a
critical point $u$ of $\Phi$ satisfying
$\Phi(u)=c_i$ and  $C_{S^1,\, 2(i-1)}(\Psi_a,\,S^1\cdot u)\neq 0$
for each $i\in\N$ (cf. Proposition 2.11, here we denote also by $u$ the corresponding
critical point of $\Psi_a$).
Applying the common index jump theorem  of Long and Zhu (cf. Theorem 3.9),
we obtain infinitely many tuples $(T, m_1, m_2, m_3)$
such that
\bea &&\Phi^\prime(u_{j_k}^{l_{j_k}})=0,\quad \Phi(u_{j_k}^{l_{j_k}})=c_{T+1-k},
\qquad C_{S^1,\; 2T-2k}(\Psi_a, \;S^1\cdot u_{j_k}^{l_{j_k}})\neq 0,
\lb{1.5}\eea
for $1\le k\le 4$ (cf. (\ref{4.21})), where $u_j^m$ denotes the critical point of
$\Phi$ (or $\Psi_a$) corresponding to the $m$-th iteration $(m\tau_j,\,y_j)$ of a
prime closed characteristic $(\tau_j,\,y_j)$.
Moreover, we have $l_{j_k}=2m_{j_k}$ for $1\le k\le 3$
and $j_1, j_2, j_3$ are pairwise distinct.

Fix a tuple $(T^\ast, m_1^\ast, m_2^\ast, m_3^\ast)$ and $(j_k^\ast,
l^\ast_{j_k^\ast})$ satisfying (\ref{1.5}). By the assumption
$^\#\T(\Sg)=3$, we can derive
$l^\ast_{j^\ast_4}=2m_{j_4^\ast}^\ast-1$, (assume $j_4^\ast=1$
without loss of generality), and either $i(y_1,\,1)=5$ together with
$\gamma_{y_1}(\tau_1)$ can be connected within
$\Omega^0(\gamma_{y_1}(\tau_1))$ to $N_1(1, 1)\diamond
N_1(1,-1)^{\diamond3}$ or $i(y_1,\,1)=4$ together with  $\gamma_{y_1}(\tau_1)$
can be connected within $\Omega^0(\gamma_{y_1}(\tau_1))$ to $N_1(1,
1)\diamond N_1(1,-1)^{\diamond2}\diamond M^\prime$ for some
$M^\prime\in\Sp(2)$ with $M^\prime\in\{R(\theta), \,D(\lambda),\,N_1(-1, b),\, I_2\}$ (cf. Cases A, B in \S 4 and \S 3 for notations), where $b=\pm 1, 0$, $R(\theta)$ is a rotation matrix
with rotation angle $\theta$ and $D(\lambda)$ is hyperbolic.

Suppose $(T, m_1, m_2, m_3)$ is any tuple found by the common index jump theorem
that satisfying (\ref{1.5}).

Now we describe the ideas of proofs of two typical cases, other cases can be handled
similarly. 

(I) If $i(y_1,\,1)=5$ and  $\gamma_{y_1}(\tau_1)$ can be connected
within $\Omega^0(\gamma_{y_1}(\tau_1))$ to $N_1(1, 1)\diamond
N_1(1,-1)^{\diamond3}$ holds. Then by the periodic property of
critical modules (cf. Proposition 2.6), we have (cf. Lemma 5.1 for details) \bea C_{S^1,\;
2T-2}(\Psi_a, \;S^1\cdot u_1^{2m_1})\cong C_{S^1,\;
2T^\ast-8}(\Psi_a, \;S^1\cdot u_1^{2m^\ast_1-1})\neq 0. \lb{1.6}\eea
Hence by the critical point theory (cf. Proposition 2.7), we
have \bea C_{S^1,\; 2T-2-2l}(\Psi_a, \;S^1\cdot u_1^{2m_1})=0,\quad
\forall l\neq 0.\lb{1.7}\eea In fact, 
we have $i(y_1^{2m_1})+\nu(y_1^{2m_1})-1=2T-2$, and then $u_1^{2m_1}$ is a local
maximum of $\Psi_a$ restricted to a local characteristic manifold of
$\Psi_a$ by (\ref{1.6}), thus (\ref{1.7}) holds. Hence we have
$c_T=\Phi(u_1^{2m_1})$ by (\ref{1.5}), and then  we have $\Phi(u_1^{2m_1})>\Phi(u_i^{2m_i})$
for $i=2, 3$.
This  contradict to the commutative property for closed characteristics (cf.
Proposition 4.5). Hence Theorem 1.1 holds in this case.

(II) If $i(y_1,\,1)=4$ and  $\gamma_{y_1}(\tau_1)$ can be connected
within $\Omega^0(\gamma_{y_1}(\tau_1))$ to $N_1(1, 1)\diamond
N_1(1,-1)^{\diamond2}\diamond R(\theta)$ with $\theta/\pi\in\Q$
holds. This is the most complicated case in this paper. We can
compute out $i(y_1^{2m_1})=2T-6$ and
$i(y_1^{2m_1})+\nu(y_1^{2m_1})-1=2T-2$ (cf. (\ref{5.16})). By (\ref{1.5}), there are three sub-cases: 

(II-a) If  $C_{S^1,\; 2T^\ast-2}(\Psi_a, \;S^1\cdot u_1^{2m_1^\ast})\neq 0$,
i.e.,  $u_1^{2m_1^\ast}$ is a local maximum of $\Psi_a$ restricted to
a local characteristic manifold of $\Psi_a$ at $u_1^{2m_1^\ast}$,
then we have $C_{S^1,\; 2T^\ast-2-l}(\Psi_a, \;S^1\cdot u_1^{2m_1^\ast})=0$
for $l\neq 0$. Thus by the periodic property of critical modules, we have
$C_{S^1,\; 2T-2-l}(\Psi_a, \;S^1\cdot u_1^{2m_1})=0$
for $l\neq 0$. Hence we have $c_T=\Phi(u_1^{2m_1})$ by (\ref{1.5}), 
and then  we have $\Phi(u_1^{2m_1})>\Phi(u_i^{2m_i})$
for $i=2, 3$.
This  contradict to the commutative property for closed characteristics
and proves Theorem 1.1 in this case.

(II-b) If  $C_{S^1,\; 2T^\ast-6}(\Psi_a, \;S^1\cdot u_1^{2m_1^\ast})\neq 0$,
i.e.,  $u_1^{2m_1^\ast}$ is a local minimum of $\Psi_a$ restricted to
a local characteristic manifold of $\Psi_a$ at $u_1^{2m_1^\ast}$,
then we have $C_{S^1,\; 2T^\ast-6+l}(\Psi_a, \;S^1\cdot u_1^{2m_1^\ast})=0$
for $l\neq 0$. Thus by the periodic property of critical modules, we have
$C_{S^1,\; 2T-6+l}(\Psi_a, \;S^1\cdot u_1^{2m_1})=0$
for $l\neq 0$. Hence we have $c_{T-2}=\Phi(u_1^{2m_1})$ by (\ref{1.5}), 
and then  we have $\Phi(u_1^{2m_1})<\Phi(u_i^{2m_i})$
for $i=2, 3$.
This  contradict to the commutative property for closed characteristics
and proves Theorem 1.1 in this case.

(II-c) It remains to consider the case $C_{S^1,\; 2T^\ast-4}(\Psi_a, \;S^1\cdot u_1^{2m_1^\ast})\neq 0$,
i.e.,  $u_1^{2m_1^\ast}$ is neither a local maximum nor a
 local minimum of $\Psi_a$ restricted to
a local characteristic manifold of $\Psi_a$ at $u_1^{2m_1^\ast}$,
then $C_{S^1,\; 2T^\ast-2}(\Psi_a, \;S^1\cdot u_1^{2m_1^\ast})=0$
and $C_{S^1,\; 2T^\ast-6}(\Psi_a, \;S^1\cdot u_1^{2m_1^\ast})=0$
by critical point theory.
Thus by the periodic property of critical modules, we have
$C_{S^1,\; 2T-2}(\Psi_a, \;S^1\cdot u_1^{2m_1})=0$
and $C_{S^1,\; 2T-6}(\Psi_a, \;S^1\cdot u_1^{2m_1})=0$.

Then the proof of Theorem 1.1 in this case  contains the following steps:

 (1)  Firstly by the commutative property for closed characteristics
 in the common index jump intervals, we can show: There exist two
tuples $(T, m_1, m_2, m_3)$ and $(T^\prime, m_1^\prime, m_2^\prime,
m_3^\prime)$ such that $c_T=\Phi(u_2^{2m_2})$,
$c_{T-2}=\Phi(u_3^{2m_3})$ and
$c_{T^\prime}=\Phi(u_3^{2m_3^\prime})$,
$c_{T^\prime-2}=\Phi(u_2^{2m_2^\prime})$ (cf. Claim 1 in Lemma 5.6). This
implies $c_{T-1}=\Phi(u_1^{2m_1})$ and
$c_{T^\prime-1}=\Phi(u_1^{2m^\prime_1})$. Hence the positions of
appropriate iterations of $u_1$ in these two common index jump intervals
 are fixed and the positions of appropriate iterations of $u_2$ and $u_3$
in these intervals interchanged, while the orders of appropriate iterations of any two closed characteristics in these two   intervals interchanged.

(2) Using (1) and the precise index iteration formula of Long, we can derive: The matrix $\gamma_{y_2}(\tau_2),\;
\gamma_{y_3}(\tau_3)$ can be connected within $\Omega^0(\gamma_{y_2}(\tau_2)),\;\Omega^0(\gamma_{y_3}(\tau_3))$ to
 $N_1(1, 1)\diamond R(\vartheta_1)\diamond R(\vartheta_2)\diamond M_2^\prime$
 and $N_1(1, 1)\diamond R(\varphi_1)\diamond R(\varphi_2)\diamond M_3^\prime$
 with $\frac{\vartheta_i}{\pi},\, \frac{\varphi_i}{\pi}\notin\Q$
 for $i=1, 2$ and $M_2^\prime, M_3^\prime\in \Sp(2)$.
Moreover, $M_2^\prime, M_3^\prime\in\{I_2, N_1(1, -1), -I_2, N_1(-1, 1),
R(\vartheta)\}$ (cf. Claim 2 in Lemma 5.6). This implies that both $\gamma_{y_2}(\tau_2)$
and $\gamma_{y_3}(\tau_3)$ have  special forms (i.e., they have some stability property), the fact that
$\frac{\vartheta_i}{\pi},\, \frac{\varphi_i}{\pi}\notin\Q$
is essential in our study below, i.e., this is the condition
for us to use Kronecker's uniform distribution theorem.

(3) By the  classification in (2) of $\gamma_{y_2}(\tau_2),\;
\gamma_{y_3}(\tau_3)$,  we can show:
$C_{S^1,\,2k+1}(\Psi_a,\,S^1\cdot u_j^m)=0$ for $k\in\Z$, $m\in\N$
and $j=2, 3$ (cf. Claim 3 in Lemma 5.6). This implies that the critical
modules of iterations of both $(\tau_2, y_2)$ and $(\tau_3, y_3)$ 
have no contribution to the number 
$$M_{2k+1}=\sum_{1\le j\le 3,\,m\in\N}
\rank C_{S^1,\,2k+1}(\Psi_a,\,S^1\cdot u_j^m).$$

(4) Using (3) and a careful study on the Morse series od $\Psi_a$, we have:
\bea\sum_{i\in\Z}(-1)^i\rank C_{S^1,\,i}(\Psi_a, \,S^1\cdot u_1^m)=1,\quad\forall m\in\N.
\lb{1.8}\eea
(cf. Claim 4 in Lemma 5.6). This implies  that the critical modules of
iterations of $(\tau_1, y_1)$ behave like those of a non-degenerate
critical point in the sense that the alternative sum of their ranks
is $1$.

(5) Using (3), (4) and Morse inequality, we can derive: It is impossible that
 $C_{S^1,\,2K}(\Psi_a, \,S^1\cdot u_1^m)\neq 0$
 and $C_{S^1,\,2K}(\Psi_a, \,S^1\cdot u_j^k)\neq 0$
 hold simultaneously for some $K, m, k\in\N$
 and some $j\in\{2,\, 3\}$. This implies that the critical modules 
 of iterations of $(\tau_1, y_1)$ and $(\tau_j, y_j)$ for $j\in\{2, 3\}$
 can not hit together. In fact, in the
 Morse inequality
\bea
 M_i-M_{i-1}+\cdots +(-1)^{i}M_0 \ge b_i-b_{i-1}+\cdots +(-1)^{i}b_0,
   \qquad\forall \;i\in\Z, \lb{1.9}\eea
if $M_{i_1}=b_{i_1}$ and $M_{i_2}=b_{i_2}$ hold
for some $i_1<i_2$, then we have
\bea
 M_{i_2}-M_{i_2-1}+\cdots +(-1)^{i_2-i_1}M_{i_1} =
b_{i_2}-b_{i_2-1}+\cdots +(-1)^{i_2-i_1}b_{i_1}. \lb{1.10}\eea
Using (3), (4) and (\ref{1.10}) properly, we can derive the above result
 (cf. Claim 5 in Lemma 5.6 for details).

(6) Up to now, the problem is transformed to find appropriate $K, m, k\in\N$
such that $C_{S^1,\,2K}(\Psi_a, \,S^1\cdot u_1^m)\neq 0$
 and $C_{S^1,\,2K}(\Psi_a, \,S^1\cdot u_2^k)\neq 0$ hold
 simultaneously. Using the precise index iteration formula
 (cf. Theorem 3.7), this is transformed further to a problem in number theory,
 i.e., whether an appropriate integer valued equation has integer solutions
 (cf. Cases 1-4 in Lemma 5.6 for the precise form of the equation).
By a case-by-case study on the possible form of $M_2^\prime$ and
Kronecker's uniform distribution theorem,
this equation actually has integer solutions in each case.
As mentioned in (2), the crucial point is that
$\frac{\vartheta_1}{\pi},\, \frac{\vartheta_2}{\pi}\notin\Q$,
this enables us to use Kronecker's uniform distribution theorem
to find solutions of the equation. This proves Theorem 1.1 in this case.

In Section 2, we
review briefly the equivariant Morse theory for closed
characteristics on compact convex hypersurfaces in $\R^{2n}$
developed in \cite{WHL} and the Fadell-Rabinowitz index theorey
applied to the study of closed characteristics. In Section 3, we review the
index iteration theory developed by Long and his coworkers.
In Section 4, we prove  a commutative property for closed characteristics
in the common index jump intervals.

In this paper, let $\N$, $\N_0$, $\Z$, $\Q$, $\R$,  $\C$ and $\R^+$ denote
the sets of natural integers, non-negative integers, integers,
rational numbers, real numbers, complex numbers and positive real numbers
respectively. Denote by $a\cdot b$ and $|a|$ the standard inner
product and norm in $\R^{2n}$. Denote by $\langle\cdot,\cdot\rangle$
and $\|\cdot\|$ the standard $L^2$-inner product and $L^2$-norm. For
an $S^1$-space $X$, we denote by $X_{S^1}$ the homotopy quotient of
$X$ module the $S^1$-action, i.e., $X_{S^1}=S^\infty\times_{S^1}X$.
We define the functions \be \left\{\matrix{[a]=\max\{k\in\Z\,|\,k\le
a\}, & E(a)=\min\{k\in\Z\,|\,k\ge a\} , \cr
                   \varphi(a)=E(a)-[a],  & \{a\}=a-[a]. \cr}\right. \lb{1.11}\ee
Specially, $\varphi(a)=0$ if $ a\in\Z\,$, and $\varphi(a)=1$ if
$a\notin\Z\,$. In this paper we use only $\Q$-coefficients for all
homological modules. For a $\Z_m$-space pair $(A, B)$, let
$H_{\ast}(A, B)^{\pm\Z_m}= \{\sigma\in H_{\ast}(A,
B)\,|\,L_{\ast}\sigma=\pm \sigma\}$, where $L$ is a generator of the
$\Z_m$-action.

\setcounter{equation}{0}
\section{ Critical point theory for closed characteristics}

In the rest of this paper, we fix a $\Sg\in\H(2n)$ and assume the
following condition on $\Sg$:

\noindent (F) {\bf There exist only finitely many geometrically
distinct closed characteristics \\$\quad \{(\tau_j, y_j)\}_{1\le
j\le q}$ on $\Sigma$. }

In this section, we review briefly the equivariant Morse theory for
closed characteristics on $\Sg$ developed in \cite{WHL} and
\cite{W1} which will be used in Section 4 and 5 of this paper. All the
details of proofs can be found in \cite{WHL} or \cite{W1}.

Let $\hat{\tau}=\inf\{\tau_j|\;1\le j\le q\}$.   Then by
\S2 of \cite{WHL}, for any $a>\hat{\tau}$, we can construct a
function $\varphi_a\in C^\infty(\R,\R^+)$ which has $0$ as its
unique critical point in $[0,\,+\infty)$ such that $\varphi_a$ is
strictly convex for $t\ge 0$. Moreover,
$\frac{\varphi_a^\prime(t)}{t}$ is strictly decreasing for $t> 0$
together with $\lim_{t\rightarrow
0^+}\frac{\varphi_a^\prime(t)}{t}=1$ and
$\varphi_a(0)=0=\varphi_a^\prime(0)$ (cf. Propositions 2.2-2.4 in \cite{WHL}).

Let $j: \R^{2n}\rightarrow\R$ be the gauge function of $\Sigma$, i.e.,
$j(\lambda x)=\lambda$ for $x\in\Sigma$ and $\lambda\ge0$, then
$j\in C^3(\R^{2n}\setminus\{0\}, \R)\cap C^0(\R^{2n}, \R)$
and $\Sigma=j^{-1}(1)$.  Define the Hamiltonian function $H_a(x)=a\varphi_a(j(x))$ and
consider the fixed period problem \be
\left\{\matrix{\dot{x}(t)=JH_a^\prime(x(t)), \cr
     x(1)=x(0).         \cr }\right. \lb{2.1}\ee
Then  $H_a\in C^3(\R^{2n}\setminus\{0\}, \R)\cap C^1(\R^{2n}, \R)$
is strictly convex. Solutions of (\ref{2.1}) are $x\equiv0$ and
$x=\rho y(\tau t)$ with
$\frac{\varphi_a^\prime(\rho)}{\rho}=\frac{\tau}{a}$, where $(\tau,
y)$ is a solution of (\ref{1.1}). In particular, nonzero solutions
of (\ref{2.1}) are one to one correspondent to solutions of
(\ref{1.1}) with period $\tau<a$.

Now  we use the Clarke-Ekeland dual action principle to transform (\ref{2.1})
to a variational problem and use variational methods to study the problem.
As usual, let $G_a$ be the Fenchel transform of $H_a$ defined by
$G_a(y)=\sup\{x\cdot y-H_a(x)\;|\; x\in \R^{2n}\}$. Then $G_a\in
C^2(\R^{2n}\bs\{0\},\R)\cap C^1(\R^{2n},\R)$ is strictly convex. Let
\be L_0^2(S^1, \;\R^{2n})= \left\{u\in L^2([0, 1],\;\R^{2n})
   \left|\frac{}{}\right.\int_0^1u(t)dt=0\right\}.  \lb{2.2}\ee
Define a linear operator $M: L_0^2(S^1,\R^{2n})\to
L_0^2(S^1,\R^{2n})$ by $\frac{d}{dt}Mu(t)=u(t)$,
$\int_0^1Mu(t)dt=0$. The dual action functional on $L_0^2(S^1,
\;\R^{2n})$ is defined by \be
\Psi_a(u)=\int_0^1\left(\frac{1}{2}Ju\cdot Mu+G_a(-Ju)\right)dt.
   \lb{2.3}\ee
Then the functional $\Psi_a\in C^{1, 1}(L_0^2(S^1,\; \R^{2n}),\;\R)$
is bounded from below and satisfies the Palais-Smale condition.
Suppose $x$ is a solution of (\ref{2.1}). Then $u=\dot{x}$ is a
critical point of $\Psi_a$. Conversely, suppose $u$ is a critical
point of $\Psi_a$. Then there exists a unique $\xi\in\R^{2n}$ such
that $Mu-\xi$ is a solution of (\ref{2.1}). In particular, solutions
of (\ref{2.1}) are in one to one correspondence with critical points
of $\Psi_a$. Moreover, $\Psi_a(u)<0$ for every critical point
$u\not= 0$ of $\Psi_a$.

Suppose $u$ is a nonzero critical point of $\Psi_a$. Then following
\cite{Eke3} the formal Hessian of $\Psi_a$ at $u$ is defined by
$$ Q_a(v,\; v)=\int_0^1 (Jv\cdot Mv+G_a^{\prime\prime}(-Ju)Jv\cdot Jv)dt, $$
which defines an orthogonal splitting $L_0^2(S^1,\; \R^{2n})=E_-\oplus E_0\oplus
E_+$ of $L_0^2(S^1,\; \R^{2n})$ into negative, zero and positive
subspaces. The index of $u$ is defined by $i(u)=\dim E_-$ and the
nullity of $u$ is defined by $\nu(u)=\dim E_0$. Let $u=\dot{x}$ be
the critical point of $\Psi_a$ such that $x$ corresponds to a
closed characteristic $(\tau,\,y)$ on $\Sigma$. Then the index
$i(u)$ and the nullity $\nu(u)$ defined above coincide with the
Ekeland indices defined by I. Ekeland in \cite{Eke1} and
\cite{Eke3}. In particular,  $1\le \nu(u)\le 2n-1$ always holds.

We have a natural $S^1$-action on $L_0^2(S^1,\; \R^{2n})$ defined by
$\th\cdot u(t)=u(\th+t)$ for all $\th\in S^1$ and $t\in\R$. Clearly
$\Psi_a$ is $S^1$-invariant. For any $\kappa\in\R$, we denote by \be
\Lambda_a^\kappa=\{w\in L_0^2(S^1,\;
\R^{2n})\;|\;\Psi_a(w)\le\kappa\}.
          \lb{2.4}\ee
For a critical point $u$ of $\Psi_a$, we denote by \be
\Lambda_a(u)=\Lambda_a^{\Psi_a(u)}
  =\{w\in L_0^2(S^1,\; \R^{2n}) \;|\; \Psi_a(w)\le\Psi_a(u)\}.\lb{2.5}\ee
Clearly, both sets are $S^1$-invariant. Since the $S^1$-action
preserves $\Psi_a$, if $u$ is a critical point of $\Psi_a$, then the
whole orbit $S^1\cdot u$ is formed by critical points of $\Psi_a$.
Denote by $crit(\Psi_a)$ the set of critical points of $\Psi_a$.
Note that by the condition (F), the number of critical orbits of
$\Psi_a$ is finite. Hence as usual we can make the following
definition.

{\bf Definition 2.1.} {\it Suppose $u$ is a nonzero critical point
of $\Psi_a$ and $\Nn$ is an $S^1$-invariant open neighborhood of
$S^1\cdot u$ such that $crit(\Psi_a)\cap(\Lambda_a(u)\cap
\Nn)=S^1\cdot u$. Then the $S^1$-critical modules of $S^1\cdot u$
are defined by}
$$ C_{S^1,\; k}(\Psi_a, \;S^1\cdot u)
=H_k((\Lambda_a(u)\cap\Nn)_{S^1},\; ((\Lambda_a(u)\setminus
S^1\cdot u)\cap\Nn)_{S^1}),\qquad k\in\Z. $$

We have the following proposition for critical modules.

{\bf Proposition 2.2.} (Proposition 3.2 of \cite{WHL}) {\it The
critical module $C_{S^1,\;k}(\Psi_a, \;S^1\cdot u)$ is independent
of $a$ in the sense that if $x_i$ are solutions of (\ref{2.1}) with
Hamiltonian functions $H_{a_i}(x)\equiv a_i\varphi_{a_i}(j(x))$ for
$i=1$ and $2$ respectively such that both $x_1$ and $x_2$ correspond
to the same closed characteristic $(\tau, y)$ on $\Sigma$. Then we
have}
$$ C_{S^1,\; k}(\Psi_{a_1}, \;S^1\cdot\dot {x}_1) \cong
  C_{S^1,\; k}(\Psi_{a_2}, \;S^1\cdot \dot {x}_2), \quad \forall k\in \Z. $$

Now let $u\neq 0$ be a critical point of $\Psi_a$ with multiplicity
$mul(u)=m$, i.e., $u$ corresponds to a closed characteristic
$(m\tau, y)\subset\Sigma$ with $(\tau, y)$ being prime. Hence
$u(t+\frac{1}{m})=u(t)$ holds for all $t\in \R$ and the orbit of
$u$, namely, $S^1\cdot u\cong S^1/\Z_m\cong S^1$. Let $f: N(S^1\cdot
u)\rightarrow S^1\cdot u$ be the normal bundle of $S^1\cdot u$ in
$L_0^2(S^1,\; \R^{2n})$ and let $f^{-1}(\theta\cdot u)=N(\theta\cdot
u)$ be the fibre over $\theta\cdot u$, where $\theta\in S^1$. Let
$DN(S^1\cdot u)$ be the $\varrho$-disk bundle of $N(S^1\cdot u)$ for
some $\varrho>0$ sufficiently small, i.e., $DN(S^1\cdot u)=\{\xi\in
N(S^1\cdot u)\;| \; \|\xi\|<\varrho\}$ and let $DN(\theta\cdot
u)=f^{-1}(\th\cdot u)\cap DN(S^1\cdot u)$ be the disk over
$\theta\cdot u$. Clearly, $DN(\theta\cdot u)$ is $\Z_m$-invariant
and we have $DN(S^1\cdot u)=DN(u)\times_{\Z_m}S^1$, where the
$Z_m$-action is given by
$$ (\th, v, t)\in \Z_m\times DN(u)\times S^1\mapsto
        (\th\cdot v, \;\theta^{-1}t)\in DN(u)\times S^1. $$
Hence for an $S^1$-invariant subset $\Gamma$ of $DN(S^1\cdot u)$, we
have $\Gamma/S^1=(\Gamma_u\times_{\Z_m}S^1)/S^1=\Gamma_u/\Z_m$,
where $\Gamma_u=\Gamma\cap DN(u)$. Since $\Psi_a$ is not $C^2$ on
$L_0^2(S^1,\; \R^{2n})$, we can not use Morse theory to study $\Psi_a$ dircetly.
In order to overcome this difficulty, we  use a finite dimensional
approximation introduced by Ekeland in \cite{Eke1} and apply Morse theory
to the obtained finite dimensional submanifold.
More precisely, we can construct a finite dimensional submanifold
$\Gamma(\iota)$ of $L_0^2(S^1,\; \R^{2n})$ which admits a
$\Z_\iota$-action with $m|\iota$. Moreover $\Psi_a$ and
$\Psi_a|_{\Gamma(\iota)}$ have the same critical points.
$\Psi_a|_{\Gamma(\iota)}$ is $C^2$ in a small tubular neighborhood
of the critical orbit $S^1\cdot u$ and the Morse index and nullity
of its critical points coincide with those of the corresponding
critical points of $\Psi_a$.  Let \be D_\iota N(S^1\cdot
u)=DN(S^1\cdot u)\cap\Gamma(\iota), \quad D_\iota N(\theta\cdot
u)=DN(\theta\cdot u)\cap\Gamma(\iota). \lb{2.6}\ee Then we have \be
C_{S^1,\; \ast}(\Psi_a, \;S^1\cdot u) \cong H_\ast(\Lambda_a(u)\cap
D_\iota N(u),\;
    (\Lambda_a(u)\setminus\{u\})\cap D_\iota N(u))^{\Z_m}. \lb{2.7}\ee
Now we can apply the results of Gromoll and Meyer in \cite{GrM1} to
the manifold $D_{p\iota}N(u^p)$ with $u^p$ as its unique critical
point, where $p\in\N$ is fixed. Then $mul(u^p)=pm$ is the multiplicity of
$u^p$ and the isotropy group $\Z_{pm}\subseteq S^1$ of $u^p$ acts on
$D_{p\iota}N(u^p)$ by isometries. According to Lemma 1 of
\cite{GrM1}, we have a $\Z_{pm}$-invariant decomposition of
$T_{u^p}(D_{p\iota}N(u^p))$
$$ T_{u^p}(D_{p\iota}N(u^p))
=V^+\oplus V^-\oplus V^0=\{(x_+, x_-, x_0)\}  $$ with $\dim
V^-=i(u^p)$, $\dim V^0=\nu(u^p)-1$ and a $\Z_{pm}$-invariant
neighborhood $B=B_+\times B_-\times B_0$ for $0$ in
$T_{u^p}(D_{p\iota}N(u^p))$ together with two $Z_{pm}$-invariant
diffeomorphisms
$$\Psi :B=B_+\times B_-\times B_0\rightarrow
\Psi(B_+\times B_-\times B_0)\subset D_{p\iota}N(u^p)$$ and
$$ \eta : B_0\rightarrow W(u^p)\equiv\eta(B_0)\subset D_{p\iota}N(u^p)$$
such that $\Psi(0)=\eta(0)=u^p$ and \be
\Psi_a\circ\Psi(x_+,x_-,x_0)=|x_+|^2 - |x_-|^2 +
\Psi_a\circ\eta(x_0),
    \lb{2.8}\ee
with $d(\Psi_a\circ \eta)(0)=d^2(\Psi_a\circ\eta)(0)=0$. As
\cite{GrM1}, we call $W(u^p)$ a local {\it characteristic manifold}
and $U(u^p)=B_-$ a local {\it negative disk} at $u^p$. By the proof
of Lemma 1 of \cite{GrM1}, $W(u^p)$ and $U(u^p)$ are
$\Z_{pm}$-invariant. Then we have \bea && H_\ast(\Lambda_a(u^p)\cap
D_{p\iota}N(u^p),\;
  (\Lambda_a(u^p)\setminus\{u^p\})\cap D_{p\iota}N(u^p)) \nn\\
= &&H_\ast (U(u^p),\;U(u^p)\setminus\{u^p\}) \otimes
H_\ast(W(u^p)\cap \Lambda_a(u^p),\; (W(u^p)\setminus\{u^p\})\cap
\Lambda_a(u^p)),
  \lb{2.9}\eea
where \be H_j(U(u^p),U(u^p)\setminus\{u^p\} )
    = \left\{\matrix{\Q, & {\rm if\;}j=i(u^p),  \cr
                      0, & {\rm otherwise}. \cr}\right.  \lb{2.10}\ee
Now we have the following proposition.

{\bf Proposition 2.3.} (Proposition 3.10 of \cite{WHL}) {\it Let
$u\neq 0$ be a critical point of $\Psi_a$ with $mul(u)=1$. Then for
all $p\in\N$ and $j\in\Z$, we have \be C_{S^1,\; j}(\Psi_a,
\;S^1\cdot u^p)\cong \left(\frac{}{}H_{j-i(u^p)}(W(u^p)\cap
\Lambda_a(u^p),\; (W(u^p)\setminus\{u^p\})\cap
\Lambda_a(u^p))\right)^{\beta(u^p)\Z_p},
  \lb{2.11}\ee
where $\beta(u^p)=(-1)^{i(u^p)-i(u)}$. Thus \be C_{S^1,\; j}(\Psi_a,
\;S^1\cdot u^p)=0, \quad {\rm for}\;\;
   j<i(u^p) \;\;{\rm or}\;\;j>i(u^p)+\nu(u^p)-1. \lb{2.12}\ee
In particular, if $u^p$ is non-degenerate, i.e., $\nu(u^p)=1$, then}
\be C_{S^1,\; j}(\Psi_a, \;S^1\cdot u^p)
    = \left\{\matrix{\Q, & {\rm if\;}j=i(u^p)\;{\rm and\;}\beta(u^p)=1,  \cr
                      0, & {\rm otherwise}. \cr}\right.  \lb{2.13}\ee

We make the following definition.

{\bf Definition 2.4.} {\it Let $u\neq 0$ be a critical point of
$\Psi_a$ with $mul(u)=1$. Then for all $p\in\N$ and $l\in\Z$, let
\bea k_{l, \pm 1}(u^p)&=&\dim\left(\frac{}{}H_l(W(u^p)\cap
\Lambda_a(u^p),\;
(W(u^p)\setminus\{u^p\})\cap \Lambda_a(u^p))\right)^{\pm\Z_p}, \nn\\
k_l(u^p)&=&\dim\left(\frac{}{}H_l(W(u^p)\cap \Lambda_a(u^p),
(W(u^p)\setminus\{u^p\})\cap
\Lambda_a(u^p))\right)^{\beta(u^p)\Z_p}. \nn\eea $k_l(u^p)$'s are
called critical type numbers of $u^p$. }

We have the following periodic property for critical type numbers.

{\bf Proposition 2.5.} (Lemma 3.12 of \cite{WHL}) {\it Let
$u\neq 0$ be a critical point of $\Psi_a$ with $mul(u)=1$. Suppose that
$\nu(u^m)=\nu(u^{pm})$ for some $p, m\in\N$, then  we have
$k_{l, \pm 1}(u^m)=k_{l, \pm 1}(u^{pm})$ for all  $l\in\Z$.    }

{\bf Proposition 2.6.} (Proposition 3.13 of \cite{WHL}) {\it Let
$u\neq 0$ be a critical point of $\Psi_a$ with $mul(u)=1$. Then
there exists a minimal $K(u)\in \N$ such that
$$ \nu(u^{p+K(u)})=\nu(u^p),\quad i(u^{p+K(u)})-i(u^p)\in 2\Z.$$
Moreover,  we have $k_l(u^{p+K(u)})=k_l(u^p)$ for all $p\in \N$ and $l\in\Z$.    }

In fact, denote by $\gamma_y$ the associated symplectic path of $(\tau, y)$, where 
$(\tau, y)$ is the closed characteristic corresponding to $u$.
Suppose $\lambda_i=e^{\pm\frac{r_i}{s_i}2\pi\sqrt{-1}}$
the eigenvalues of $\gamma_y(\tau)$ possessing rotation angles which are rational
multiple of $2\pi$ with $r_i$, $s_i\in\N$ and $(r_i,s_i)=1$ for
$1\le i\le k$. Let $K^\prime(u)$ be the least common multiple of
$s_1,\ldots, s_k$. Then  we have $\nu(u^{p+K^\prime(u)})=\nu(u^p)$
for all $p\in\N$. By Theorem 3.6 below and Theorem 9.3.4 of \cite{Lon4},  we have $i(u^{m+2})-i(u^m)\in2\Z$ for any $m\in\N$. Hence we have 
\bea K(u)=\left\{\matrix{2K^\prime(u)&&{\rm if}\quad i(u^2)-i(u)\in2\Z+1\;{\rm and}\; K^\prime(u)\in 2\N-1,\cr
K^\prime(u)\quad&&{\rm otherwise.} \cr}\right.
\nn\eea

For a prime closed characteristic $(\tau,y)$ on $\Sigma$, we denote by
$y^m\equiv (m\tau, y)$ the $m$-th iteration of $y$ for $m\in\N$. Let
$a>\tau$ be large enough and choose $\vf_a$ as above. Determine $\rho$ uniquely by
$\frac{\vf_a'(\rho)}{\rho}=\frac{\tau}{a}$. Let $x=\rho y(\tau t)$
and $u=\dot{x}$. Then we define the index $i(y^m)$ and nullity
$\nu(y^m)$ of $(m\tau,y)$ for $m\in\N$ by
$$ i(y^m)=i(u^m), \qquad \nu(y^m)=\nu(u^m). $$
These indices are independent of $a$ when $a$ tends to infinity. Now
the mean index of $(\tau,y)$ is defined by
$$ \hat{i}(y)=\lim_{m\rightarrow\infty}\frac{i(y^m)}{m}. $$
Note that $\hat{i}(y)>2$ always holds which was proved by Ekeland
and Hofer in \cite{EkH1} of 1987 (cf. Corollary 8.3.2 and Lemma
15.3.2 of \cite{Lon4} for a different proof).

By Proposition 2.2, we can define the critical type numbers
$k_l(y^m)$ of $y^m$ to be $k_l(u^m)$, where $u^m$ is the critical
point of $\Psi_a$ corresponding to $y^m$. We also define
$K(y)=K(u)$. Then we have the following.

{\bf Proposition 2.7.} (Proposition 2.6 of \cite{W1})  {\it We have $k_l(y^m)=0$ for $l\notin [0,
\nu(y^m)-1]$ and it can take only values $0$ or $1$ when $l=0$ or
$l=\nu(y^m)-1$. Moreover, the following properties hold:

(i) $k_0(y^m)=1$ implies $k_l(y^m)=0$ for $1\le l\le \nu(y^m)-1$.

(ii) $k_{\nu(y^m)-1}(y^m)=1$ implies $k_l(y^m)=0$ for $0\le l\le
\nu(y^m)-2$.

(iii) $k_l(y^m)\ge 1$ for some $1\le l\le \nu(y^m)-2$ implies
$k_0(y^m)=k_{\nu(y^m)-1}(y^m)=0$.

(iv) If $i(y^m)-i(y)\in 2\Z+1$ for some $m\in\N$, then $k_0(y^m)=0$.}

Let $\Psi_a$ be the functional defined by (\ref{2.3}) for some
$a\in\R$ large enough and let $\varepsilon>0$ be small enough such
that $[-\varepsilon, +\infty)\setminus\{0\}$ contains no critical
values of $\Psi_a$. Denote by $I_a$ the greatest integer in $\N_0$
such that $I_a<i(\tau, y)$ hold for all closed characteristics
$(\tau,\, y)$ on $\Sigma$ with $\tau\ge a$. Then by P. 447-448 of
\cite{WHL}, we have \be H_{S^1,\; i}(\Lambda_a^{-\varepsilon} )
\cong H_{S^1,\; i}( \Lambda_a^\infty)
  \cong H_i(CP^\infty), \quad \forall i<I_a.  \lb{2.14}\ee
For any $i\in\Z$, let \be  M_i(\Lambda_a^{-\varepsilon})
  =\sum_{1\le j\le q,\,1\le m_j<a/\tau_j} \dim C_{S^1,\;i}(\Psi_a, \;S^1\cdot u_j^{m_j}).
  \lb{2.15} \ee
Then the equivariant Morse inequalities for the space
$\Lambda_a^{-\varepsilon}$ yield \bea M_i(\Lambda_a^{-\varepsilon})
       &\ge& b_i(\Lambda_a^{-\varepsilon}),\lb{2.16}\\
M_i(\Lambda_a^{-\varepsilon}) &-& M_{i-1}(\Lambda_a^{-\varepsilon})
    + \cdots +(-1)^{i}M_0(\Lambda_a^{-\varepsilon}) \nn\\
&\ge& b_i(\Lambda_a^{-\varepsilon}) -
b_{i-1}(\Lambda_a^{-\varepsilon})
   + \cdots + (-1)^{i}b_0(\Lambda_a^{-\varepsilon}), \lb{2.17}\eea
for $i\in\Z$, where $b_i(\Lambda_a^{-\varepsilon})=\dim H_{S^1,\;
i}(\Lambda_a^{-\varepsilon})$. Now we have the following Morse
inequalities for closed characteristics.

{\bf Theorem 2.8.} (Theorem 2.8 of \cite{W1})  {\it Suppose  $\Sigma\in \H(2n)$ satisfy
$\,^{\#}\T(\Sg)<+\infty$. Denote all the geometrically
distinct closed characteristics on $\Sg$ by $\{(\tau_j,\; y_j)\}_{1\le j\le
q}$. Let \bea
M_i&=&\lim_{a\rightarrow+\infty\atop \varepsilon\rightarrow 0}M_i(\Lambda_a^{-\varepsilon}),\quad
                  \forall i\in\Z,\lb{2.18}\\
b_i &=& \lim_{a\rightarrow+\infty\atop \varepsilon\rightarrow 0}b_i(\Lambda_a^{-\varepsilon})=
\left\{\matrix{1, & {\rm if\;}i\in 2\N_0,  \cr
                      0, & {\rm otherwise}. \cr}\right.  \lb{2.19}
\eea Then we have}
\bea M_i &\ge& b_i,\qquad\forall i\in\Z\lb{2.20}\\
 M_i-M_{i-1}+\cdots +(-1)^{i}M_0 &\ge& b_i-b_{i-1}+\cdots +(-1)^{i}b_0,
   \qquad\forall \;i\in\Z. \lb{2.21}\eea

Recall that for a principal $U(1)$-bundle $E\to B$, the Fadell-Rabinowitz index
(cf. \cite{FaR1}) of $E$ is defined to be $\sup\{k\;|\, c_1(E)^{k-1}\not= 0\}$,
where $c_1(E)\in H^2(B,\Q)$ is the first rational Chern class. For a $U(1)$-space,
i.e., a topological space $X$ with a $U(1)$-action, the Fadell-Rabinowitz index is
defined to be the index of the bundle $X\times S^{\infty}\to X\times_{U(1)}S^{\infty}$,
where $S^{\infty}\to CP^{\infty}$ is the universal $U(1)$-bundle.

As in P. 199 of \cite{Eke3}, choose some $\alpha\in(1,\, 2)$ and associate with $U$
a convex function $H$ such that $H(\lambda x)=\lambda^\alpha H(x)$ for $\lambda\ge 0$.
Consider the fixed period problem
\be \left\{\matrix{\dot{x}(t)=JH^\prime(x(t)), \cr
     x(1)=x(0).         \cr }\right. \lb{2.22}\ee
Define
\be L_0^{\frac{\alpha}{\alpha-1}}(S^1,\R^{2n})
  =\left\{u\in L^{\frac{\alpha}{\alpha-1}}(S^1,\R^{2n})\,\left|\frac{}{}\right.\,\int_0^1udt=0\right\}. \lb{2.23}\ee
The corresponding Clarke-Ekeland dual action functional is defined by
\be \Phi(u)=\int_0^1\left(\frac{1}{2}Ju\cdot Mu+H^{\ast}(-Ju)\right)dt,
    \qquad \forall\;u\in L_0^{\frac{\alpha}{\alpha-1}}(S^1,\R^{2n}), \lb{2.24}\ee
where $Mu$ is defined by $\frac{d}{dt}Mu(t)=u(t)$ and $\int_0^1Mu(t)dt=0$,
$H^\ast$ is the Fenchel transform of $H$ defined above.

For any $\kappa\in\R$, we denote by
\be \Phi^{\kappa-}=\{u\in L_0^{\frac{\alpha}{\alpha-1}}(S^1,\R^{2n})\;|\;
             \Phi(u)<\kappa\}. \lb{2.25}\ee
Then as in P. 218 of \cite{Eke3}, we define
\be c_i=\inf\{\delta\in\R\;|\: \hat I(\Phi^{\delta-})\ge i\},\lb{2.26}\ee
where $\hat I$ is the Fadell-Rabinowitz index given above. Then by Proposition 3
in P. 218 of \cite{Eke3}, we have

{\bf Proposition 2.9.} {\it Every $c_i$ is a critical value of $\Phi$. If
$c_i=c_j$ for some $i<j$, then there are infinitely many geometrically
distinct closed characteristics on $\Sg$.}

As in Definition 2.1, we define the following

{\bf Definition 2.10.} {\it Suppose $u$ is a nonzero critical
point of $\Phi$, and $\Nn$ is an $S^1$-invariant
open neighborhood of $S^1\cdot u$ such that
$crit(\Phi)\cap(\Lambda(u)\cap \Nn)=S^1\cdot u$. Then
the $S^1$-critical modules of $S^1\cdot u$ is defined by
\bea C_{S^1,\; k}(\Phi, \;S^1\cdot u)
=H_k((\Lambda(u)\cap\Nn)_{S^1},\;
((\Lambda(u)\setminus S^1\cdot u)\cap\Nn)_{S^1}),\qquad k\in\Z,\lb{2.27}
\eea
where $\Lambda(u)=\{w\in L_0^{\frac{\alpha}{\alpha-1}}(S^1,\R^{2n})\;|\;
\Phi(w)\le\Phi(u)\}$.}

Comparing with Theorem 4 in P. 219 of \cite{Eke3}, we have the following

{\bf Proposition 2.11.}  (Proposition 3,5 of \cite{W1}) {\it For every $i\in\N$, there exists a point
$u\in L_0^{\frac{\alpha}{\alpha-1}}(S^1,\R^{2n})$ such that}
\bea
\Phi^\prime(u)=0,\quad \Phi(u)=c_i,
\quad C_{S^1,\; 2(i-1)}(\Phi, \;S^1\cdot u)\neq 0. \lb{2.28}\eea

The next proposition implies that $\Psi_a$ and $\Phi$ have isomorphic
critical modules at corresponding critical points,
 thus we can compute the critical modules of $\Phi$ via that of $\Psi_a$.

{\bf Proposition 2.12.} {\it Suppose $u$ is the critical point of $\Phi$ found
in Proposition 2.11. Then we have
\be C_{S^1,\;k}(\Psi_a, \;S^1\cdot u_a)\cong C_{S^1,\; k}(\Phi,\;S^1\cdot u),\quad\forall k\in\Z, \lb{2.29}\ee
where $\Psi_a$ is given by (\ref{2.3}) and $u_a\in L_0^2(S^1,\;\R^{2n})$
is its critical point corresponding to $u$ in the natural sense.}

{\bf Proof.} Fix this $u$, we  modify the function $H$ only in a small
neighborhood $\Omega$ of $0$ as in \cite{Eke1} so that the corresponding
orbit of $u$ does not enter $\Omega$ and the resulted function $\widetilde{H}$
satisfies similar properties as Definition 1 in P. 26
of \cite{Eke1} by just replacing $\frac{3}{2}$ there by $\alpha$.
Define the dual action functional
$\widetilde{\Phi}:L_0^{\frac{\alpha}{\alpha-1}}(S^1,\R^{2n})\to\R$ by
\be \widetilde{\Phi}(v)=\int_0^1\left(\frac{1}{2}Jv\cdot
   Mv+\widetilde{H}^{\ast}(-Jv)\right)dt. \lb{2.30}\ee
Clearly $\Phi$ and $\widetilde{\Phi}$ are $C^1$ close to each other,
thus by the continuity of critical modules (cf. Theorem 8.8 of \cite{MaW1} or
Theorem 1.5.6 in P. 53 of \cite{Cha1}, which can be easily generalized to the
equivariant sense) for the $u$ in the proposition, we have
\be C_{S^1,\; \ast}(\Phi, \;S^1\cdot u)\cong C_{S^1,\; \ast}(\widetilde{\Phi},
    \;S^1\cdot u).\lb{2.31}\ee
Using a finite dimensional approximation as in Lemma 3.9 of \cite{Eke1},
we have
\be C_{S^1,\; \ast}(\widetilde{\Phi}, \;S^1\cdot u)
\cong H_\ast(\widetilde{\Lambda}(u)\cap D_\iota N(u),\;
    (\widetilde{\Lambda}(u)\setminus\{u\})\cap D_\iota N(u))^{\Z_m}, \lb{2.32}\ee
where $\widetilde{\Lambda}(u)=\{w\in L_0^{\frac{\alpha}{\alpha-1}}(S^1,\R^{2n})\;|\;
\widetilde{\Phi}(w)\le\widetilde{\Phi}(u)\}$ and $D_\iota N(u)$ is a
$\Z_m$-invariant finite dimensional disk transversal to $S^1\cdot u$ at $u$
(cf. Lemma 3.9 of \cite{Eke1}), $m$ is the multiplicity of $u$.

By Lemma 3.9 of \cite{WHL}, we have
\be C_{S^1,\; \ast}(\Psi_a, \;S^1\cdot u_a)
\cong H_\ast(\Lambda_a(u_a)\cap D_\iota N(u_a),\;
    (\Lambda_a(u_a)\setminus\{u_a\})\cap D_\iota N(u_a))^{\Z_m}.\lb{2.33}\ee
By the construction of $H_a$ in \cite{WHL}, $H_a=\widetilde{H}$ in a
$L^\infty$-neighborhood of $S^1\cdot u$. We remark here that multiplying $H$ by a constant
will not affect the corresponding critical modules, i.e., the corresponding
critical orbits have isomorphic critical modules. Hence we can assume
$H_a=H$ in a $L^\infty$-neighborhood of $S^1\cdot u$ and then the above
conclusion holds. Hence $\Psi_a$ and $\widetilde{\Phi}$ coincide
in a $L^\infty$-neighborhood of $S^1\cdot u$. Note also by Lemma 3.9 of
\cite{Eke1}, the two finite dimensional approximations are actually the same.
Hence we have
\bea
&& H_\ast(\widetilde{\Lambda}(u)\cap D_\iota N(u),\;
   (\widetilde{\Lambda}(u)\setminus\{u\})\cap D_\iota N(u))^{\Z_m}\nn\\
\cong&& H_\ast(\Lambda_a(u_a)\cap D_\iota N(u_a),\;
    (\Lambda_a(u_a)\setminus\{u_a\})\cap D_\iota N(u_a))^{\Z_m}.\lb{2.34}\eea
Now the proposition follows from (\ref{2.31})-(\ref{2.34}).
\hfill\hb

\setcounter{equation}{0}
\section{ Index iteration theory for closed characteristics}

In this section, we recall briefly an index theory for symplectic paths
developed by Y. Long and his coworkers. All the details can be found in \cite{Lon4} or \cite{LoZ1}.
Then we use this theory to study the Morse indices of the critical points $u^m$ in \S2,

As usual, the symplectic group $\Sp(2n)$ is defined by
$$ \Sp(2n) = \{M\in {\rm GL}(2n,\R)\,|\,M^TJM=J\}, $$
whose topology is induced from that of $\R^{4n^2}$. For $\tau>0$ we are interested
in paths in $\Sp(2n)$:
$$ \P_{\tau}(2n) = \{\ga\in C([0,\tau],\Sp(2n))\,|\,\ga(0)=I_{2n}\}, $$
which is equipped with the topology induced from that of $\Sp(2n)$. The
following real function was introduced in \cite{Lon3}:
$$ D_{\om}(M) = (-1)^{n-1}\ol{\om}^n\det(M-\om I_{2n}), \qquad
          \forall \om\in\U,\, M\in\Sp(2n). $$
Thus for any $\om\in\U$ the following codimension $1$ hypersurface in $\Sp(2n)$ is
defined in \cite{Lon3}:
$$ \Sp(2n)_{\om}^0 = \{M\in\Sp(2n)\,|\, D_{\om}(M)=0\}.  $$
For any $M\in \Sp(2n)_{\om}^0$, we define a co-orientation of $\Sp(2n)_{\om}^0$
at $M$ by the positive direction $\frac{d}{dt}Me^{t\ep J}|_{t=0}$ of
the path $Me^{t\ep J}$ with $0\le t\le 1$ and $\ep>0$ being sufficiently
small. Let
\bea
\Sp(2n)_{\om}^{\ast} &=& \Sp(2n)\bs \Sp(2n)_{\om}^0,   \nn\\
\P_{\tau,\om}^{\ast}(2n) &=&
      \{\ga\in\P_{\tau}(2n)\,|\,\ga(\tau)\in\Sp(2n)_{\om}^{\ast}\}, \nn\\
\P_{\tau,\om}^0(2n) &=& \P_{\tau}(2n)\bs  \P_{\tau,\om}^{\ast}(2n).  \nn\eea
For any two continuous arcs $\xi$ and $\eta:[0,\tau]\to\Sp(2n)$ with
$\xi(\tau)=\eta(0)$, it is defined as usual:
$$ \eta\ast\xi(t) = \left\{\matrix{
            \xi(2t), & \quad {\rm if}\;0\le t\le \tau/2, \cr
            \eta(2t-\tau), & \quad {\rm if}\; \tau/2\le t\le \tau. \cr}\right. $$
Given any two $2m_k\times 2m_k$ matrices of square block form
$M_k=\left(\matrix{A_k&B_k\cr
                                C_k&D_k\cr}\right)$ with $k=1, 2$,
as in \cite{Lon4}, the $\;\dm$-product of $M_1$ and $M_2$ is defined by
the following $2(m_1+m_2)\times 2(m_1+m_2)$ matrix $M_1\dm M_2$:
$$ M_1\dm M_2=\left(\matrix{A_1&  0&B_1&  0\cr
                               0&A_2&  0&B_2\cr
                             C_1&  0&D_1&  0\cr
                               0&C_2&  0&D_2\cr}\right). \nn$$  
Denote by $M^{\dm k}$ the $k$-fold $\dm$-product $M\dm\cdots\dm M$. Note
that the $\dm$-product of any two symplectic matrices is symplectic. For any two
paths $\ga_j\in\P_{\tau}(2n_j)$ with $j=0$ and $1$, let
$\ga_0\dm\ga_1(t)= \ga_0(t)\dm\ga_1(t)$ for all $t\in [0,\tau]$.

A special path $\xi_n$ is defined by
\be \xi_n(t) = \left(\matrix{2-\frac{t}{\tau} & 0 \cr
                                             0 &  (2-\frac{t}{\tau})^{-1}\cr}\right)^{\dm n}
        \qquad {\rm for}\;0\le t\le \tau.  \lb{3.1}\ee

{\bf Definition 3.1.} (cf. \cite{Lon3}, \cite{Lon4}) {\it For any $\om\in\U$ and
$M\in \Sp(2n)$, define
\be  \nu_{\om}(M)=\dim_{\C}\ker_{\C}(M - \om I_{2n}).  \lb{3.2}\ee
For any $\tau>0$ and $\ga\in \P_{\tau}(2n)$, define
\be  \nu_{\om}(\ga)= \nu_{\om}(\ga(\tau)).  \lb{3.3}\ee

If $\ga\in\P_{\tau,\om}^{\ast}(2n)$, define
\be i_{\om}(\ga) = [\Sp(2n)_{\om}^0: \ga\ast\xi_n],  \lb{3.4}\ee
where the right hand side of (\ref{3.4}) is the usual homotopy intersection
number, and the orientation of $\ga\ast\xi_n$ is its positive time direction under
homotopy with fixed end points.

If $\ga\in\P_{\tau,\om}^0(2n)$, we let $\mathcal{F}(\ga)$
be the set of all open neighborhoods of $\ga$ in $\P_{\tau}(2n)$, and define
\be i_{\om}(\ga) = \sup_{U\in\mathcal{F}(\ga)}\inf\{i_{\om}(\beta)\,|\,
                       \beta\in U\cap\P_{\tau,\om}^{\ast}(2n)\}.
               \lb{3.5}\ee
Then
$$ (i_{\om}(\ga), \nu_{\om}(\ga)) \in \Z\times \{0,1,\ldots,2n\}, $$
is called the index function of $\ga$ at $\om$. }

Note that when $\om=1$, this index theory was introduced by
C. Conley-E. Zehnder in \cite{CoZ1} for the non-degenerate case with $n\ge 2$,
Y. Long-E. Zehnder in \cite{LZe1} for the non-degenerate case with $n=1$,
and Y. Long in \cite{Lon1} and C. Viterbo in \cite{Vit2} independently for
the degenerate case. The case for general $\om\in\U$ was defined by Y. Long
in \cite{Lon3} in order to study the index iteration theory (cf. \cite{Lon4}
for more details and references).

For any symplectic path $\ga\in\P_{\tau}(2n)$ and $m\in\N$,  we
define its $m$-th iteration $\ga^m:[0,m\tau]\to\Sp(2n)$ by
\be \ga^m(t) = \ga(t-j\tau)\ga(\tau)^j, \qquad
  {\rm for}\quad j\tau\leq t\leq (j+1)\tau,\;j=0,1,\ldots,m-1.
     \lb{3.6}\ee
We still denote the extended path on $[0,+\infty)$ by $\ga$.

{\bf Definition 3.2.} (cf. \cite{Lon3}, \cite{Lon4}) {\it For any $\ga\in\P_{\tau}(2n)$,
we define
\be (i(\ga,m), \nu(\ga,m)) = (i_1(\ga^m), \nu_1(\ga^m)), \qquad \forall m\in\N.
   \lb{3.7}\ee
The mean index $\hat{i}(\ga,m)$ per $m\tau$ for $m\in\N$ is defined by
\be \hat{i}(\ga,m) = \lim_{k\to +\infty}\frac{i(\ga,mk)}{k}. \lb{3.8}\ee
For any $M\in\Sp(2n)$ and $\om\in\U$, the {\it splitting numbers} $S_M^{\pm}(\om)$
of $M$ at $\om$ are defined by
\be S_M^{\pm}(\om)
     = \lim_{\ep\to 0^+}i_{\om\exp(\pm\sqrt{-1}\ep)}(\ga) - i_{\om}(\ga),
   \lb{3.9}\ee
for any path $\ga\in\P_{\tau}(2n)$ satisfying $\ga(\tau)=M$.}

For a given path $\gamma\in {\cal P}_{\tau}(2n)$ we consider to deform
it to a new path $\eta$ in ${\cal P}_{\tau}(2n)$ so that
\begin{equation}
i_1(\gamma^m)=i_1(\eta^m),\quad \nu_1(\gamma^m)=\nu_1(\eta^m), \quad
         \forall m\in {\bf N}, \label{3.10}
\end{equation}
and that $(i_1(\eta^m),\nu_1(\eta^m))$ is easy enough to compute. This
leads to finding homotopies $\delta:[0,1]\times[0,\tau]\to {\rm Sp}(2n)$
starting from $\gamma$ in ${\cal P}_{\tau}(2n)$ and keeping the end
points of the homotopy always stay in a certain suitably chosen maximal
subset of ${\rm Sp}(2n)$ so that (\ref{3.10}) always holds. In fact,  this
set was first discovered in \cite{Lon3} as the path connected component
$\Omega^0(M)$ containing $M=\gamma(\tau)$ of the set
\begin{eqnarray}
  \Omega(M)=\{N\in{\rm Sp}(2n)\,&|&\,\sigma(N)\cap{\bf U}=\sigma(M)\cap{\bf U}\;
{\rm and}\;  \nonumber\\
 &&\qquad \nu_{\lambda}(N)=\nu_{\lambda}(M),\;\forall\,
\lambda\in\sigma(M)\cap{\bf U}\}. \label{3.11}
\end{eqnarray}
Here $\Omega^0(M)$ is called the {\it homotopy component} of $M$ in
${\rm Sp}(2n)$.

In \cite{Lon3} and \cite{Lon4}, the following symplectic matrices were introduced
as {\it basic normal forms}:
\begin{eqnarray}
D(\lambda)=\left(\matrix{\lm & 0\cr
         0  & \lm^{-1}\cr}\right), &\quad& \lm=\pm 2,\lb{3.12}\\
N_1(\lm,b) = \left(\matrix{\lm & b\cr
         0  & \lm\cr}\right), &\quad& \lm=\pm 1, b=\pm1, 0, \lb{3.13}\\
R(\th)=\left(\matrix{\cos\th & -\sin\th\cr
        \sin\th  & \cos\th\cr}\right), &\quad& \th\in (0,\pi)\cup(\pi,2\pi),
                     \lb{3.14}\\
N_2(\om,b)= \left(\matrix{R(\th) & b\cr
              0 & R(\th)\cr}\right), &\quad& \th\in (0,\pi)\cup(\pi,2\pi),
                     \lb{3.15}\end{eqnarray}
where $b=\left(\matrix{b_1 & b_2\cr
               b_3 & b_4\cr}\right)$ with  $b_i\in\R$ and  $b_2\not=b_3$.

Splitting numbers possess the following properties:

{\bf Lemma 3.3.} (cf. \cite{Lon3} and Lemma 9.1.5 of \cite{Lon4}) {\it Splitting
numbers $S_M^{\pm}(\om)$ are well defined, i.e., they are independent of the choice
of the path $\ga\in\P_\tau(2n)$ satisfying $\ga(\tau)=M$ appeared in (\ref{3.9}).
For $\om\in\U$ and $M\in\Sp(2n)$, splitting numbers $S_N^{\pm}(\om)$ are constant
for all $N\in\Om^0(M)$. }

{\bf Lemma 3.4.} (cf. \cite{Lon3}, Lemma 9.1.5 and List 9.1.12 of \cite{Lon4})
{\it For $M\in\Sp(2n)$ and $\om\in\U$, there hold
\begin{eqnarray}
S_M^{\pm}(\om) &=& 0, \qquad {\it if}\;\;\om\not\in\sg(M).  \lb{3.16}\\
S_{N_1(1,a)}^+(1) &=& \left\{\matrix{1, &\quad {\rm if}\;\; a\ge 0, \cr
0, &\quad {\rm if}\;\; a< 0. \cr}\right. \lb{3.17}\eea

For any $M_i\in\Sp(2n_i)$ with $i=0$ and $1$, there holds }
\be S^{\pm}_{M_0\dm M_1}(\om) = S^{\pm}_{M_0}(\om) + S^{\pm}_{M_1}(\om),
    \qquad \forall\;\om\in\U. \lb{3.18}\ee

We have the following symplectic additivity property for index functions: 

{\bf Theorem 3.5.} (cf. Theorem 6.1 of \cite{LoZ1} or Theorem 6.2.7 of \cite{Lon4}) {\it For
any $\gamma_j\in\P_\tau(2n_j)$ with $j=0, 1$, we have 
\bea i_\omega(\ga_0\diamond\ga_1)=i_\omega(\ga_0)+i_\omega(\ga_1).\lb{3.19}\eea}

Let $\Sigma\in\H(2n)$. Using notations in \S1,
for any closed characteristic $(\tau,y)$ on $\Sg$ and $m\in\N$, we define
its $m$-th iteration $y^m:\R/(m\tau\Z)\to\R^{2n}$ by
\be y^m(t) = y(t-j\tau), \qquad {\rm for}\quad j\tau\leq t\leq (j+1)\tau,
       \quad j=0,1,2,\ldots, m-1. \lb{3.20}\ee
Note that this coincide with that in \S2.
We still denote by $y$ its extension to $[0,+\infty)$.

We define via Definition 3.2 the following
\bea  S^+(y) &=& S_{\ga_y(\tau)}^+(1),  \lb{3.21}\\
  (i(y,m), \nu(y,m)) &=& (i(\ga_y,m), \nu(\ga_y,m)),  \lb{3.22}\\
   \hat{i}(y,m) &=& \hat{i}(\ga_y,m),  \lb{3.23}\eea
for all $m\in\N$, where $\ga_y$ is the associated symplectic path of $(\tau,y)$.
Then we have the following.

{\bf Theorem 3.6.} (cf. Lemma 1.1 of \cite{LoZ1}, Theorem 15.1.1 of \cite{Lon4}) {\it Suppose
$(\tau,y)$ is a closed characteristic on $\Sg$. Then we have
\be i(y^m)\equiv i(m\tau ,y)=i(y, m)-n,\quad \nu(y^m)\equiv\nu(m\tau, y)=\nu(y, m),
       \qquad \forall m\in\N, \lb{3.24}\ee
where $i(y^m)$ and $\nu(y^m)$ are the index and nullity
defined in \S2. }

The following is the precise index iteration formulae for
symplectic paths, which is due to Y. Long
(cf. Chapter 8 of \cite{Lon4} or Theorems 6.5 and 6.7 of \cite{LoZ1}).

{\bf Theorem 3.7.} {\it Let $\ga\in\P_{\tau}(2n)$. Then there exists a path $f\in
C([0,1],\Omega^0(\gamma(\tau))$ such that $f(0)=\gamma(\tau)$ and
\bea f(1)=&&N_1(1,1)^{\diamond p_-} \diamond I_{2p_0}\diamond
N_1(1,-1)^{\diamond p_+}
\diamond N_1(-1,1)^{\diamond q_-} \diamond (-I_{2q_0})\diamond
N_1(-1,-1)^{\diamond q_+}\nn\\
&&\diamond R(\theta_1)\diamond\cdots\diamond R(\theta_r)
\diamond N_2(\omega_1, u_1)\diamond\cdots\diamond N_2(\omega_{r_*}, u_{r_*}) \nn\\
&&\diamond N_2(\lm_1, v_1)\diamond\cdots\diamond N_2(\lm_{r_0}, v_{r_0})
\diamond M_0 \lb{3.25}\eea
where $ N_2(\omega_j, u_j) $s are
non-trivial and   $ N_2(\lm_j, v_j)$s  are trivial basic normal
forms; $\sigma (M_0)\cap U=\emptyset$; $p_-$, $p_0$, $p_+$, $q_-$,
$q_0$, $q_+$, $r$, $r_*$ and $r_0$ are non-negative integers;
$\omega_j=e^{\sqrt{-1}\alpha_j}$, $
\lambda_j=e^{\sqrt{-1}\beta_j}$; $\theta_j$, $\alpha_j$, $\beta_j$
$\in (0, \pi)\cup (\pi, 2\pi)$; these integers and real numbers
are uniquely determined by $\gamma(\tau)$. Then using the
functions defined in (\ref{1.11}), we have
\bea i(\gamma, m)=&&m(i(\gamma,
1)+p_-+p_0-r)+2\sum_{j=1}^r E\left(\frac{m\theta_j}{2\pi}\right)-r
-p_--p_0\nn\\&&-\frac{1+(-1)^m}{2}(q_0+q_+)+2\left(
\sum_{j=1}^{r_*}\varphi\left(\frac{m\alpha_j}{2\pi}\right)-r_*\right).
\lb{3.26}\eea
\bea \nu(\gamma, m)=&&\nu(\gamma,
1)+\frac{1+(-1)^m}{2}(q_-+2q_0+q_+)+2(r+r_*+r_0)\nn\\
&&-2\left(\sum_{j=1}^{r}\varphi\left(\frac{m\theta_j}{2\pi}\right)+
\sum_{j=1}^{r_*}\varphi\left(\frac{m\alpha_j}{2\pi}\right)
+\sum_{j=1}^{r_0}\varphi\left(\frac{m\beta_j}{2\pi}\right)\right)\lb{3.27}\eea
\bea \hat i(\gamma, 1)=i(\gamma, 1)+p_-+p_0-r+\sum_{j=1}^r
\frac{\theta_j}{\pi}.\lb{3.28}\eea
Where $N_1(1, \pm 1)=
\left(\matrix{ 1 &\pm 1\cr 0 & 1\cr}\right)$, $N_1(-1, \pm 1)=
\left(\matrix{ -1 &\pm 1\cr 0 & -1\cr}\right)$,
$R(\theta)=\left(\matrix{\cos\th &
                  -\sin\th\cr\sin\th & \cos\th\cr}\right)$,
$ N_2(\omega, b)=\left(\matrix{R(\th) & b
                  \cr 0 & R(\th)\cr}\right)$    with some
$\th\in (0,\pi)\cup (\pi,2\pi)$ and $b=
\left(\matrix{ b_1 &b_2\cr b_3 & b_4\cr}\right)\in\R^{2\times2}$,
such that $(b_2-b_3)\sin\theta>0$, if $ N_2(\omega, b)$ is
trivial; $(b_2-b_3)\sin\theta<0$, if $ N_2(\omega, b)$ is
non-trivial. We have $i(\gamma, 1)$ is odd if $f(1)=N_1(1, 1)$, $I_2$,
$N_1(-1, 1)$, $-I_2$, $N_1(-1, -1)$ and $R(\theta)$; $i(\gamma, 1)$ is
even if $f(1)=N_1(1, -1)$ and $ N_2(\omega, b)$; $i(\gamma, 1)$ can be any
integer if $\sigma (f(1)) \cap \U=\emptyset$.}

We have the following properties in the index iteration theory.

{\bf Theorem 3.8. } (cf. Theorem 2.3 of \cite{LoZ1})
{\it Let $\gamma\in\P_\tau(2n)$ and $M = \gamma(\tau)$. Suppose that there exist
$P\in \Sp(2n)$ and $Q \in \Sp(2n -2)$ such that $M=P^{-1}(N_1(1, 1)\diamond Q)P$.
Then for any $m\in\N$, there holds
$$\nu(\gamma, m)-\frac{e(M)}{2}+1\le i(\gamma, m+1)-i(\gamma, m)-i(\gamma, 1)
\le \nu(\gamma, 1)-\nu(\gamma, m+1)+\frac{e(M)}{2},$$
where $e(M)$ is the total algebraic multiplicity of all eigenvalues of $M$ on
the unit circle $\U$ in the complex plane $\C$. }

The following is the common index jump theorem of Y. Long and C. Zhu.

{\bf Theorem 3.9.} (cf. Theorems 4.1-4.4 of \cite{LoZ1}) {\it   Let $\gamma_k\in\P_{\tau_k}(2n)$ for
$ k = 1,\ldots,q$ be a finite collection of symplectic paths. Let $M_k = \gamma_k(\tau_k)$.
 Suppose that there exist
$P_k\in \Sp(2n)$ and $Q _k\in \Sp(2n -2)$ such that $M_k=P_k^{-1}(N_1(1, 1)\diamond Q_k)P_k$ and $\hat i(\gamma_k, 1) > 0$, for all $k = 1,\ldots,q$.
Then there exist infinitely many $(T,m_1,\ldots,m_q)\in\N^{q+1}$  such that
\bea \nu(\gamma_k, 2m_k -1)&=&  \nu(\gamma_k, 1),\lb{3.29}\\
\nu(\gamma_k, 2m_k +1)&=&   \nu(\gamma_k, 1),\lb{3.30}\\
i(\gamma_k, 2m_k -1)+\nu(\gamma_k, 2m_k -1)&=&
2T-\left(i(\gamma_k, 1)+2S^+_{M_k}(1)-\nu(\gamma_k, 1)\right),\lb{3.31}\\
i(\gamma_k, 2m_k+1)&=&2T+i(\gamma_k, 1),\lb{3.32}\\
i(\gamma_k, 2m_k)&\ge&2T-\frac{e(M_k)}{2}\ge 2T-n,\lb{3.33}\\
i(\gamma_k, 2m_k)+\nu(\gamma_k, 2m_k)&\le&2T+\frac{e(M_k)}{2}-1\le 2T+n-1,\lb{3.34}
\eea
for every $k=1,\ldots,q$.
Moreover we have
\bea \min\left\{\left\{\frac{m_k\theta}{\pi}\right\},\,1-\left\{\frac{m_k\theta}{\pi}\right\}\right\}
<\delta,\lb{3.35}\eea
whenever $e^{\sqrt{-1}\theta}\in\sigma(M_k)$
and $\delta$ can be chosen as small as we want (cf. (4.43) of \cite{LoZ1}). More precisely, by (4.10), (4.40) and (4.41) in \cite{LoZ1} , we have
\bea m_k=\left(\left[\frac{T}{M\hat i(\gamma_k, 1)}\right]+\chi_k\right)M,\quad 1\le k\le q,\lb{3.36}\eea
where $\chi_k=0$ or $1$ for $1\le k\le q$ and $\frac{M\theta}{\pi}\in\Z$
whenever $e^{\sqrt{-1}\theta}\in\sigma(M_k)$ and $\frac{\theta}{\pi}\in\Q$
for some $1\le k\le q$. Furthermore, given $M_0\in\N$,
by the proof of Theorem 4.1 of \cite{LoZ1}, we may
further require  $M_0|T$ (since the closure of the
set $\{\{Tv\}: T\in\N, \;M_0|T\}$ is still a closed
additive subgroup of $\bf T^h$ for some $h\in\N$,
where we use notations  as (4.21)-(4.22)  in \cite{LoZ1}.
Then we can use the proof of Step 2 in Theorem 4.1
of \cite{LoZ1} to get $T$).

In fact, by (4.40)-(4.41) of \cite{LoZ1}, let
$\mu_i=\sum_{\theta\in(0, 2\pi)}S_{M_i}^-(e^{\sqrt{-1}\theta})$ for
$1\le i\le q$ and $\alpha_{i, j}=\frac{\theta_j}{\pi}$ where
$e^{\sqrt{-1}\theta_j}\in\sigma(M_i)$ for $1\le j\le \mu_i$ and
$1\le i\le q$. As in (4.21) of \cite{LoZ1}, let $h=q+\sum_{1\le i\le
q}\mu_i$ and \bea  v=\left(\frac{1}{M\hat i(\gamma_1,\,1)},\dots,
\frac{1}{M\hat i(\gamma_q, 1)}, \frac{\alpha_{1, 1}}{\hat
i(\gamma_1,\,1)}, \frac{\alpha_{1, 2}}{\hat i(\gamma_1,\,1)},
\dots\frac{\alpha_{1, \mu_1}}{\hat i(\gamma_1,\,1)}, \frac{\alpha_{2,
1}}{\hat i(\gamma_2,\, 1)},\dots, \frac{\alpha_{q, \mu_q}}{\hat
i(\gamma_q, 1)}\right). \lb{3.37}\eea Then by (4.22) of \cite{LoZ1},
the above theorem is equivalent to find a vertex
\bea\chi=(\chi_1,\ldots,\chi_q,\chi_{1, 1},\chi_{1,
2},\ldots,\chi_{1, \mu_1}, \chi_{2, 1},\ldots,\chi_{q,
\mu_q})\lb{3.38}\eea of the cube $[0,\,1]^h$ and infinitely many
integers $T\in\N$ such that \bea|\{Tv\}-\chi|<\epsilon\lb{3.39}\eea
 for any given $\epsilon$ small enough.}

{\bf Theorem 3.10.} (cf. Theorem 4.2 of \cite{LoZ1}) {\it Let $H$ be the closure of $\{\{mv\}|m\in\N\}$ in $\bf {T}^h=(\R/\Z)^h$ and
$V=T_0\pi^{-1}H$ be the tangent space of $\pi^{-1}H$ at the origin in $\R^h$,
where $\pi: \R^h\rightarrow  \bf{T}^h$ is the projection map.
Define
\bea A(v)=V\setminus\cup_{v_k\in\R\setminus\Q}\{x=(x_1,,\ldots,x_h)\in V | x_k=0\}.
\lb{3.40}\eea
Define $\psi(x)=0$  when $x\ge 0$ and $\psi(x)=1$ when $x<0$.
Then for any $a=(a_1,\ldots, a_h)\in A(V)$, the vector
\bea \chi=(\psi(a_1),\ldots,\psi(a_h))\lb{3.41}\eea
makes (\ref{3.39}) hold for infinitely many $T\in\N$. }

{\bf Theorem 3.11.}  (cf. Theorem 4.2 of \cite{LoZ1}) {\it We have the following properties for $A(v)$:

(i) When $v\in\R^h\setminus\Q^h$,  then $\dim V\ge 1$,
$0\notin A(v)\subset V$, $A(v)=-A(v)$ and $A(v)$ is open in $V$.

(ii) When $\dim V = 1$, then  $A(v) = V \setminus\{0\}$.

(iii) When $\dim V \ge 2$,  $A(v)$  is obtained from $V$ by deleting all the
coordinate hyperplanes  with    dimension
strictly smaller than   $\dim V$    from    $V$.}

\setcounter{equation}{0}
\section{A commutative property for closed characteristics in the common index jump intervals }

In this section, we  prove a  commutative property for closed characteristics in the common index jump intervals. This property is discovered and used firstly in this paper 
to handle the multiplicity problem. It will be essential in \S 5 below. This property is  motivated by  Theorem 5.4   of \cite{LoZ1} (which is a stability result), while we find their proof leads to our commutative property.

As Definition 1.1 of \cite{LoZ1}, we define the following:

{ \bf Definition  4.1.} For $\alpha\in(1,2)$, we define a map
$\varrho_n\colon\H(2n)\to\N\cup\{ +\infty\}$
\be \varrho_n(\Sg)
= \left\{\matrix{+\infty, & {\rm if\;\;}^\#\mathcal{V}(\Sigma,\alpha)=+\infty, \cr
\min\left\{[\frac{i(y, 1) + 2S^+(y) - \nu(y, 1)+n}{2}]\,
\left|\frac{}{}\right.\,(\tau, y)\in\mathcal{V}_\infty(\Sigma, \alpha)\right\},
 & {\rm if\;\;} ^\#\mathcal{V}(\Sigma, \alpha)<+\infty, \cr}\right.  \lb{4.1}\ee
where $\mathcal{V}(\Sigma,\alpha)$ and $\mathcal{V}_\infty(\Sigma,\alpha)$ are
variationally visible and infinite variationally visible sets respectively given
by Definition 1.4 of \cite{LoZ1} (cf. Definition 15.3.3 of \cite{Lon4}).

For a prime closed characteristic $(\tau_j, y_j)$ and  $m\in\N$,
we denote by  $u_j^m$ the unique critical point of $\Phi$ corresponding to
the closed characteristic $(m\tau_j, y_j)$ as in \S 2.

{\bf Lemma 4.2.} {\it There exists a large integer $T_0\in\N$ such that the following hold.
For every integer $i>T_0$, there exists  a prime closed characteristic $(\tau_j, y_j)$ and  $m\in\N$
such that
\bea
 \Phi^\prime(u_j^m)=0,\quad \Phi(u_j^m)=c_i, \quad
 C_{S^1,\; 2(i-1)}(\Phi, \;S^1\cdot u_j^m)\neq 0. \lb{4.2}\eea
Moreover, for any $i_1>i_2>T_0$ we have $m_{j_1}\hat i(y_{j_1})>m_{j_2}\hat i(y_{j_2})$,
where $(\tau_{j_l}, y_{j_l})$ and $m_{j_l}$ corresponds to $i_l$
via (\ref{4.2}) for $l=1, 2$.
}

{\bf Proof.} The lemma follows directly from Lemma 3.1 of \cite{LoZ1},
Proposition 2.11, Theorem 3.6 and (\ref{3.8}). \hfill\hb

By Theorem 1.1 of \cite{LoZ1} (cf. Theorem 15.4.3 of \cite{Lon4}), we have $^{\#}\T(\Sg)\ge \varrho_n(\Sg)\ge [\frac{n}{2}]+1$.
We prove Theorem 1.1 by contradiction. Hence in the following of this paper,
we  fix  a $\Sigma\in\H(8)$ and assume the following:

{\bf (C) We have $^{\#}\T(\Sg)=3$,
i.e., there are exactly three  geometrically distinct
closed characteristics  $\{(\tau_j, y_j)\}_{1\le j\le 3}$ on $\Sg$.}

Denote  by $\ga_j\equiv \gamma_{y_j}$ the associated
symplectic path of $(\tau_j,\,y_j)$ for $1\le j\le 3$. Then by Lemma 1.3 of \cite{LoZ1}
(cf. Lemma 15.2.4 of \cite{Lon4}), there exist $P_j\in \Sp(8)$ and $M_j\in \Sp(6)$ such
that
\be \ga_j(\tau_j)=P_j^{-1}(N_1(1,\,1)\dm M_j)P_j, \quad 1\le j\le 3.
  \lb{4.3}\ee
By Theorem 3.9 we obtain infinitely many
$(T, m_1, m_2, m_3)\in\N^4$
such that the following hold \bea
\nu(y_j, 2m_j -1)&=&  \nu(y_j, 1),\lb{4.4}\\
\nu(y_j, 2m_j +1)&=&   \nu(y_j, 1),\lb{4.5}\\
i(y_j,\, 2m_j) &\ge& 2T-\frac{e(\gamma_j(\tau_j))}{2}, \lb{4.6}\\
i(y_j,\, 2m_j)+\nu(y_j,\, 2m_j) &\le& 2T+\frac{e(\gamma_j(\tau_j))}{2}-1, \lb{4.7}\\
i(y_j,\, 2m_j+1) &=& 2T+i(y_j,\,1). \lb{4.8}\\
i(y_j,\, 2m_j-1)+\nu(y_j,\, 2m_j-1)
 &=& 2T-\left(i(y_j,\,1)+2S^+_{\gamma_j(\tau_j)}(1)-\nu(y_j, 1)\right).\lb{4.9}
\eea

By Corollary 1.2 of \cite{LoZ1} (cf. Corollary 15.1.4 of \cite{Lon4}), we
have $i(y_j,\,1)\ge4$ for $1\le j\le 3$. Note that
$e(\gamma_j(\tau_j))\le 8$ for $1\le j\le 3$. Hence Theorem 3.8 yields
\bea i(y_j,\, m)+\nu(y_j,\, m)
&\le&  i(y_j, m+1)-i(y_j, 1)+\frac{e(\gamma_j(\tau_j))}{2}-1\nn\\
&\le& i(y_j, m+1)-1, \quad \forall m\in\N,\;1\le j\le 3.\lb{4.10}\eea

By Theorem 3.7, the matrix $M_j$ can be connected within $\Omega^0(M_j)$ to $N_1(1,1)^{\diamond p_{j_-}} \diamond I_{2p_{j_0}}\diamond
N_1(1,-1)^{\diamond p_{j_+}}
\diamond M_j^\prime$, where $p_{j_-}, \,p_{j_0},\,p_{j_+}\in\N_0$ and $1\notin\sigma(M_j^\prime)$ for $1\le j\le 3$.
By Lemma 3.4 and (\ref{4.3}), we have
\bea
&& 2S^+_{\gamma_j(\tau_j)}(1)-\nu(y_j,\,1) \nn\\
&&\qquad = 2S^+_{N_1(1,\,1)}(1)-\nu_1(N_1(1,\,1))
   +2S^+_{M_j}(1)-\nu_1(M_j)\nn\\
&&\qquad = 1 + p_{j_-}-p_{j_+}, \qquad 1\le j\le 3.\lb{4.11}\eea
Note that by (\ref{4.1}) and (\ref{4.11}), we have $\varrho_4(\Sg)\ge 4$
if there is no closed characteristic $(\tau_j,\,y_j)$  on $\Sg$ satisfies
\bea i(y_j, 1)+2S^+_{\gamma_j(\tau_j)}(1)-\nu(y_j,\,1)
=i(y_j, 1)+1+p_{j_-}-p_{j_+}<4.\lb{4.12}\eea

 Thus in order to prove Theorem 1.1, it is sufficient to consider
 the case that there is some closed characteristic $(\tau_j,\,y_j)$ on $\Sigma$ satisfies
 $ i(y_j, 1)+2S^+_{\gamma_j(\tau_j)}(1)-\nu(y_j,\,1) <4$.
 By a permutation of $\{1, 2, 3\}$, we may assume
\be
 \left\{\matrix{i(y_j, 1)+2S^+_{\gamma_j(\tau_j)}(1)-\nu(y_j,\,1) <4,\quad{\rm if}\quad 1\le j\le K, \cr
i(y_j, 1)+2S^+_{\gamma_j(\tau_j)}(1)-\nu(y_j,\,1)\ge 4,\quad{\rm if}\quad K< j\le 3,\cr}\right.  \lb{4.13}\ee
for some $1\le K\le 3$.

Since $i(y_j, 1)\ge 4$ and $p_{j_-}+p_{j_0}+p_{j_+}\le 3$,  a closed characteristic $(\tau_j,\,y_j)$ satisfying (\ref{4.12}) must have $p_{j_+}\ge 2$
 and $p_{j_-}=0$. Hence we have the following two possible cases:

{\bf Case A. } {\it We have $p_{j_+}=2$ and $i(y_j, 1)=4$. }

 In this case, the matrix $M_j$ can be connected within $\Omega^0(M_j)$ to
 $N_1(1,\,-1)^{\diamond 2}\diamond M^\prime$   for some $M^\prime\in\Sp(2$) and
$M^\prime\in\{R(\theta), \,D(\lambda),\,N_1(-1, b),\, I_2\}$, where $b\in\{\pm 1, 0\}$.

 {\bf Case B. } {\it We have $p_{j_+}=3$ and $i(y_j, 1)=5$. }

 In fact, if $p_{j_+}=3$, the matrix $M_j$ can be connected within $\Omega^0(M_j)$ to
 $N_1(1,\,-1)^{\diamond 3}$, hence we have $i(y_j, 1)$ must be odd by (\ref{4.3}),
 the symplectic additivity property for indices (cf. Theorem 3.5) and Theorem 3.7. On the other hand, we have $i(y_j, 1)\ge 4$, hence $i(y_j,  1)=5$ holds.

Combining these two cases, we have
\bea   i(y_j, 1)+2S^+_{\gamma_j(\tau_j)}(1)-\nu(y_j,\,1) =3,\quad 1\le j\le K.\lb{4.14}\eea

By Theorem 3.6, (\ref{4.10}) and (\ref{4.14}), (\ref{4.6})-(\ref{4.9}) become
\bea
i(y_j^{2m_j}) &\ge& 2T-8, \quad 1\le j\le3,\lb{4.15}\\
i(y_j^{2m_j})+\nu(y_j^{2m_j})-1 &\le& 2T+\frac{e(\gamma_j(\tau_j))}{2}-6\le2T-2, \quad 1\le j\le 3  \lb{4.16}\\
i(y_j^{2m_j+m}) &\ge& 2T, \quad\forall\; m\ge 1, \quad 1\le j\le 3 \lb{4.17}\\
i(y_j^{2m_j-1})+\nu(y_j^{2m_j-1})-1 &=& 2T-8,\quad 1\le j\le K, \lb{4.18}\\
i(y_j^{2m_j-m})+\nu(y_j^{2m_j-m})-1 &<& 2T-8,\quad\forall\; m\ge 2,\; 1\le j\le K,\lb{4.19}\\
i(y_j^{2m_j-m})+\nu(y_j^{2m_j-m})-1 &<& 2T-8,\quad\forall\; m\ge 1, \;K< j\le 3.\lb{4.20}
\eea
Thus by Propositions 2.11 and 2.12 and Lemma 4.2,  we can find $(j_k,\,l_{j_k})_{1\le k\le 4}$ such that
\bea &&\Phi^\prime(u_{j_k}^{l_{j_k}})=0,\quad \Phi(u_{j_k}^{l_{j_k}})=c_{T+1-k},
\qquad C_{S^1,\; 2T-2k}(\Psi_a, \;S^1\cdot u_{j_k}^{l_{j_k}})\neq 0,
\lb{4.21}
\eea
where we denote by $u_{j_k}^{l_{j_k}}$ the corresponding
critical points of $\Phi$ (or $\Psi_a$). Note that by Proposition 2.9, the numbers $c_{T+1-k}$ 
for $1\le k\le 4$ are pairwise distinct critical values of $\Phi$, thus if $j_k=j_{k^\prime}$
for some $1\le k, k^\prime\le 4$,
we have $l_{j_k}\neq l_{j_k^\prime}$.
Hence we have $(j_k,\,l_{j_k})=(j_k,\,2m_{j_k})$ for $1\le k\le 3$
and  $j_1, j_2, j_3$ are pairwise distinct, and then $\{j_1, j_2, j_3\}=\{1, 2, 3\}$.
In fact,  by Proposition 2.3  and (\ref{4.17})-(\ref{4.20}), we have
$ C_{S^1,\; 2T-2k}(\Psi_a, \;S^1\cdot u_j^m)=0$ for $1\le k, \,j\le3$
and any integer $m\neq 2m_j$.  Thus $l_{j_k}=2m_{j_k}$,
and then $\Phi(u_{j_k}^{2m_{j_k}})=c_{T+1-k}$ for $1\le k\le 3$ by (\ref{4.21}),
hence $j_1, j_2, j_3$ are pairwise distinct.

{\bf Definition 4.3.} {\it  For any tuple $(T, m_1, m_2,  m_3)$ found by Theorem  3.9
and $j_1, j_2, j_3$ satisfying (\ref{4.15})-(\ref{4.21}), 
we define its common index jump interval to be 
$[2T-6,\,2T-2]$. For any even integer $2T-2s\in [2T-6,\,2T-2]$, let 
$\xi_T(s)\in\{1, 2, 3\}$ be the unique integer satisfying 
$c_{T+1-s}=\Phi(u_{\xi_T(s)}^{2m_{\xi_T(s)}})$, i.e., we have $\xi_T(s)=j_s$. }

Let $v\in\R^h$ be the vector given by  (\ref{3.37}) associated to the symplecic
 paths $\{\gamma_1,\gamma_2,\gamma_3\}$ and $A(v)$ be the set given by (\ref{3.40})
 associated to $v$.   By Theorem 1.3 of \cite{LoZ1}, there are at least
 $\varrho_4(\Sg)-1\ge2$ geometrically distinct closed characteristics
   on $\Sg$ processing irrational mean indices.
hence $v\in\R^h\setminus\Q^h$. Thus $\dim V\ge 1$ by Theorem 3.11,
 where $V$ is given by Theorem 3.10 associated to $v$.
For any  $a=(a_1,\ldots, a_h)\in A(V)$, let $\chi(a)\equiv(\psi(a_1),\ldots,\psi(a_h))$.
 By Theorems 3.10 and 3.11, we have $-a\in A(v)$ and $\chi(a)\neq \chi(-a)$.
 For any tuples $(T,\,\chi(a))$ and $(T^{\prime},\,\chi(-a))$
 satisfying (\ref{3.39}), let $m_k=\left(\left[\frac{T}{M\hat i(y_k)}\right]+\chi(a)_k\right)M$ and  $m^{\prime}_k=\left(\left[\frac{T^{\prime}}{M\hat i(y_k)}\right]+\chi(-a)_k\right)M$
 be given by (\ref{3.36}) for $1\le k\le 3$; 

Now we prove a lemma which will be essential in our discussion below.

{\bf Lemma 4.4.} {\it Let $1\le \alpha,\,\beta\le 3$ and $\alpha\neq \beta$ be fixed.
Then there exists an $a\in A(v)$ and $T\in\N$ satisfying (\ref{3.39}) such that
$m_{\alpha}\hat i(y_\alpha)>m_\beta\hat i(y_\beta)$
and $m^\prime_{\alpha}\hat i(y_\alpha)<m^\prime_\beta\hat i(y_\beta)$,
where $m_k, m_k^\prime$ are given as above.}

{\bf Proof.}  The proof is motivated by  Theorem 5.4   of \cite{LoZ1}.

 We prove the lemma by contradiction.  Then we may assume for all $a\in A(v)$, we always have
 \bea m_{\alpha}\hat i(y_\alpha)>m_\beta\hat i(y_\beta),
 \qquad m^\prime_{\alpha}\hat i(y_\alpha)>m^\prime_\beta\hat i(y_\beta).\lb{4.22}\eea

Now we fix an  $a=(a_1,\ldots, a_h)\in A(V)$. Let $\delta_1>0$ be small enough,  $\Lambda=\max_{1\le j\le 3}\hat i(y_j)$ and
 \bea t_0=\frac{\delta_1}{6(|a|+1)(M\Lambda+1)}.\lb{4.23}\eea
 Note that in the proof of Theorem 3.9 (cf. Step 2 of Theorem 4.1 of \cite{LoZ1}),
 we can further require $T$ such that the vector $\{T v\}-\chi(a)$
are located in a sufficiently small neighborhood inside the open ball in $V$
 with radius  $\delta_1/(6M\Lambda + 1)$ and centered at $at_0$ (cf. P. 360 of \cite{LoZ1}), i.e.,
  \bea \{T v\}-\chi(a)\in V,\quad |\{T v\}-\chi(a)-at_0|<\frac{\delta_1}{6M\Lambda+1}.\lb{4.24}\eea
 Then $|\{T v\}-\chi(a)|<|at_0|+\frac{\delta_1}{6M\Lambda+1}\le \frac{\delta_1}{3}$,
hence we still have (\ref{3.39}) for $\delta_1$ small enough.

 {\it Claim. } We have $a_\alpha\hat i(y_\alpha)-a_\beta\hat i(y_\beta)=0$.

We prove it by contradiction. In fact, we can further require
 $T\in\N $ so that the following  holds:
 \bea \left|\left\{\frac{T}{M\hat i(y_k)}\right\}-\chi(a)_k-a_kt_0\right|
 <\frac{t_0}{3\Lambda}\min_{a_i\hat i(y_i)-a_j\hat i(y_j)\neq 0}|a_i\hat i(y_i)-a_j\hat i(y_j)|
 \lb{4.25}\eea
 for $1\le k\le 3$. Then we have
 \bea &&m_\alpha\hat i(y_\alpha)-m_\beta\hat i (y_\beta)\nn\\
 =&&\left(\left[\frac{T}{M\hat i(y_\alpha)}\right]+\chi(a)_\alpha\right)M\hat i(y_\alpha)
 -\left(\left[\frac{T}{M\hat i(y_\beta)}\right]+\chi(a)_\beta\right)M\hat i(y_\beta)\nn\\
 =&&\left(\chi(a)_\alpha+\frac{T}{M\hat i(y_\alpha)}-\left\{\frac{T}{M\hat i(y_\alpha)}\right\}\right)M\hat i(y_\alpha)
 -\left(\chi(a)_\beta+\frac{T}{M\hat i(y_\beta)}-\left\{\frac{T}{M\hat i(y_\beta)}\right\}\right)M\hat i(y_\beta)\nn\\
=&&\left(\chi(a)_\alpha-\left\{\frac{T}{M\hat i(y_\alpha)}\right\}\right)M\hat i(y_\alpha)
 -\left(\chi(a)_\beta-\left\{\frac{T}{M\hat i(y_\beta)}\right\}\right)M\hat i(y_\beta)\nn\\
 =&&-Mt_0(a_\alpha\hat i(y_\alpha)-a_\beta\hat i(y_\beta))
 +\left(\chi(a)_\alpha-\left\{\frac{T}{M\hat i(y_\alpha)}\right\}+a_\alpha t_0\right)M\hat i(y_\alpha)\nn\\
&& -\left(\chi(a)_\beta-\left\{\frac{T}{M\hat i(y_\beta)}\right\}+a_\beta t_0\right)M\hat i(y_\beta). \lb{4.26}\eea
 By (\ref{4.22}) and (\ref{4.25}), we have
 \bea a_\alpha\hat i(y_\alpha)\le a_\beta \hat i(y_\beta). \lb{4.27}\eea
 On the other hand, we repeat  this argument for $(T^{\prime},\,\chi(-a))$
 and obtain
  \bea -a_\alpha\hat i(y_\alpha)\le -a_\beta\hat i(y_\beta). \lb{4.28}\eea
 Combing   (\ref{4.27}) and (\ref{4.28}) we obtain the claim.

  Let $V_{\alpha, \beta}=\{a\in V\,|\, a_\alpha\hat i(y_\alpha)=a_\beta\hat i(y_\beta)\}$  and
  \bea B(v)=A(v)\setminus V_{\alpha,\beta},\quad {\rm if}\quad V_{\alpha, \beta}\neq V.\lb{4.29}\eea
 Since $\dim V\ge 1$  and  $A(v)$  is obtained from $V$ by deleting finitely many proper linear subspaces of $V$ by Theorem 3.11, and so is $B(v)$. Hence $B(v)$ is nonempty.
  Now we choose an $a\in B(v)$. By the above Claim, we have
 $a_\alpha\hat i(y_\alpha)=a_\beta\hat i(y_\beta)$. By the definition of $a \in B(v)$ we have $V_{\alpha, \beta} = V$.

 By (\ref{4.24}), the vector $\{T v\}- \chi(a)$ belongs to $V$,
 and thus belongs to $V_{\alpha, \beta} $. Then by the definition of $V_{\alpha, \beta} $, this implies
 \bea (\{T v_\alpha\}-\chi(a)_\alpha )\hat i(y_\alpha)=(\{T v_\beta\}- \chi(a)_\beta )\hat i(y_\beta).
 \lb{4.30}\eea
 By (\ref{3.37}), this implies
 \bea \left(\left\{\frac{T}{M\hat i(y_\alpha)} \right\}-\chi(a)_\alpha\right )\hat i(y_\alpha)
= \left(\left\{\frac{T}{M\hat i(y_\beta)} \right\}-\chi(a)_\beta\right )\hat i(y_\beta).
 \lb{4.31}\eea
 By the third equality of (\ref{4.26}), we have
\bea m_\alpha\hat i(y_\alpha)=m_\beta \hat i(y_\beta).\lb{4.32}\eea
 This contradict to (\ref{4.22})   and proves the lemma.\hfill\hb

Now we can give the main result in this section. It states that the
closed characteristics on $\Sigma$  have certain commutative property in the 
common index jump intervals.

{\bf Proposition 4.5.} {\it Let $1\le \alpha,\,\beta\le 3$ and $\alpha\neq \beta$ be fixed.
Then there exists an $a\in A(v)$ and $T\in\N$ satisfying (\ref{3.39}) such that
$\Phi(u_\alpha^{2m_{\alpha}})>\Phi(u_\beta^{2m_{\beta}})$
and $\Phi(u_\alpha^{2m^\prime_{\alpha}})<\Phi(u_\beta^{2m^\prime_{\beta}})$,
where $m_k, m_k^\prime$ are given as above. In particular, 
we have the following diagram 
\bea
\begin{tabular}{lccccc}
&$c_{T+1-\xi_T^{-1}(\beta)}$&$c_{T+1-\xi_T^{-1}(\alpha)}$& $\qquad$&$c_{T^\prime+1-\xi_{T^\prime}^{-1}(\alpha)}$& $c_{T^\prime+1-\xi_{T^\prime}^{-1}(\beta)}$\\
&$u_\beta^{2m_\beta}$&$u_\alpha^{2m_\alpha}$& $\qquad$&$u_\alpha^{2m_\alpha^\prime}$&$u_\beta^{2m_\beta^\prime}$\\
\end{tabular}\nn\eea
Furthermore, we have  $2T-2\xi_T^{-1}(\beta)<2T-2\xi_T^{-1}(\alpha)$
and  $2T^\prime-2\xi_{T^\prime}^{-1}(\beta)>2T^\prime-2\xi_{T^\prime}^{-1}(\alpha)$,
i.e., the orders of the two closed characteristics in the common index jump intervals
interchanged. }

{\bf Proof.} This follows directly from Lemmas 4.2, 4.4 and Definition 4.3.\hfill\hb

\setcounter{equation}{0}
\section{Proof of the main theorem  }

In this section, we give the proof of Theorem 1.1  by using Morse theory,
the index iteration theory developed by Long and his coworkers,
the commutative property for closed characteristics in the common index jump intervals
 and Kronecker's uniform distribution theorem in number theory.

We continue to use the notations as in \S 4.
First note that we have $l_{j_4}=2m_{j_4}-1$ and 
\bea &&\Phi^\prime(u_{j_4}^{2m_{j_4}-1})=0,\quad \Phi(u_{j_4}^{2m_{j_4}-1})=c_{T-3},
\qquad C_{S^1,\; 2T-8}(\Psi_a, \;S^1\cdot u_{j_4}^{2m_{j_4}-1})\neq 0,
\lb{5.1}
\eea
for some $j_4\in \{1,\ldots,K\}\subset\{1, 2, 3\}$. In fact, we have
$l_{{j_4}}\in\left\{2m_{{j_4}}-1,\, 2m_{{j_4}}\right\}\equiv\Delta$ since we have
$C_{S^1,\; 2T-8}\left(\Psi_a, \;S^1\cdot u_{j_4}^m\right)=0$
for $m\notin\Delta$ by (\ref{4.17})-(\ref{4.20}) and Proposition 2.3.
On the other hand, since $\{j_1, j_2, j_3\}=\{1, 2, 3\}$, thus we have $j_4=j_k$ for some $1\le k\le 3$. This implies $\Phi\left(u_{j_4}^{2m_{j_4}}\right)=c_{T+1-k}$
 and
$\Phi\left(u_{j_4}^{l_{j_4}}\right)=c_{T+1-4}$ by (\ref{4.21}). Hence
$l_{{j_4}}=2m_{{j_4}}-1$ by Proposition 2.9, and then (\ref{5.1}) holds by (\ref{4.20}),  (\ref{4.21}) and Proposition 2.3.

Now we  fix a tuple $(T^\ast,\,m_1^\ast,\,m_2^\ast,\,m_3^\ast)$
and $(j_k^\ast,\,l^\ast_{j_k^\ast})$ for $1\le k\le 4$ satisfying (\ref{4.15})-(\ref{4.21}).
Since $j_4^\ast\in\{1,\ldots, K\}$, by a permutation of $\{1, \ldots, K\}$,
we may assume $j_4^\ast=1$ with out loss of generality.
Thus by (\ref{4.12}) and (\ref{4.13}), $(\tau_1, y_1)$ must belong to Case A or B in \S4.
Hence we separate the proof of Theorem 1.1 into several cases according 
to the classification  of $(\tau_1, y_1)$. 

{\bf Lemma 5.1.} {\it If  $(\tau_1, y_1)$ belongs to Case B in \S4, then we have $^\#\T(\Sg)\ge 4$.}

{\bf Proof.}   Suppose $(T,\,m_1,\,m_2,\,m_3)$ is any tuple found by Theorem 3.9 and $(j_k,\,l_{j_k})$
 satisfy (\ref{4.15})-(\ref{4.21}). As mentioned in
Case B, we have  $i(y_1,\,1)=5$. Thus by Theorems 3.6,  3.7
and (\ref{4.3}), we have $i(y_1^m)=m(i(y_1,\,1)+1)-1-4=6m-5$ and
$\nu(y_1^m)=4$ for $m\in\N$. By (\ref{4.18}), we have
$i(y_1^{2m_1-1})+\nu(y_1^{2m_1-1})-1=2T-8$. Hence we have
$i(y_1^{2m_1})+\nu(y_1^{2m_1})-1=2T-2$. Hence by Propositions 2.3
and  2,6, we have $K(y_1)=1$ and
 \bea &&\rank C_{S^1,\; 2T-2}(\Psi_a, \;S^1\cdot u_1^{2m_1})\nn\\
 =&&k_{\nu(u_1^{2m_1})-1}(u_1^{2m_1})
 =k_{\nu(u_1)-1}(u_1)=k_{\nu(u_1^{2m^\ast_1-1})-1}(u_1^{2m^\ast_1-1})\nn\\
 =&&\rank C_{S^1,\; 2T^\ast-8}(\Psi_a, \;S^1\cdot u_1^{2m^\ast_1-1})\neq 0,\lb{5.2}\eea
where the last inequality follows from (\ref{5.1}).
Hence
\bea\rank C_{S^1,\; 2T-2-l}(\Psi_a, \;S^1\cdot u_1^{2m_1})
=k_{\nu(u_1^{2m_1})-1-l}(u_1^{2m_1})=0
\lb{5.3}\eea
for $l\neq 0$ by (ii) of Proposition 2.7.
Hence by (\ref{4.21}), we have $c_T=\Phi(u_1^{2m_1})$,
and then $c_{T+1-\xi_T^{-1}(i)}=\Phi(u_i^{2m_i})$
for  $i=2, 3$ and $\xi_T^{-1}(i)\in\{2, 3\}$.
Thus  we have $\Phi(u_1^{2m_1})>\Phi(u_i^{2m_i})$
for $i=2, 3$.  In particular,  for any $a\in A(v)$ and $m_k, m_k^\prime$ as in Proposition 4.5, we always have 
$\Phi(u_1^{2m_1})>\Phi(u_i^{2m_i})$ and  $\Phi(u_1^{2m^\prime_1})>\Phi(u_i^{2m^\prime_i})$
for $i=2, 3$. This contradict to Proposition 4.5
and proves the lemma.\hfill\hb

{\bf Lemma 5.2.} {\it If $(\tau_1, y_1)$ belongs to Case A in \S4  and  the matrix $M_1$ can be connected within $\Omega^0(M_1)$ to
 $N_1(1,\,-1)^{\diamond 2}\diamond M^\prime$ with $M^\prime\in \Sp(2)$ and
 $\sigma(M^\prime)\cap\U=\emptyset$, i.e.,  $M^\prime$ is huperbolic, then we have $^\#\T(\Sg)\ge 4$.}

{\bf Proof.}  As in Lemma 5.1, suppose $(T,\,m_1,\,m_2,\,m_3)$ and $(j_k,\,l_{j_k})$
 satisfy (\ref{4.15})-(\ref{4.21}). As mentioned in
Case A, we have  $i(y_1,\,1)=4$, thus by Theorems 3.6,  3.7
and (\ref{4.3}), we have $i(y_{1}^m)=m(i(y_1,\,1)+1)-1-4=5m-5$ and
$\nu(y_1^m)=3$ for $m\in\N$. By (\ref{4.18}), we have
$i(y_1^{2m_1-1})+\nu(y_1^{2m_1-1})-1=2T-8$. Hence we have
$i(y_1^{2m_1})+\nu(y_1^{2m_1})-1=2T-3$. By Propositions 2.3  and
2.6, we have $K(y_1)=2$ and then
 \bea \rank C_{S^1,\; 2T-3-l}(\Psi_a, \;S^1\cdot u_1^{2m_1})
 =k_{\nu(u_1^{2m_1})-1-l}(u_1^{2m_1})=k_{\nu(u^2_1)-1-l}(u^2_1),\lb{5.4}\\
 k_{\nu(u_1)-1}(u_1)=k_{\nu(u_1^{2m^\ast_1-1})-1}(u_1^{2m^\ast_1-1})
 =\rank C_{S^1,\; 2T^\ast-8}(\Psi_a, \;S^1\cdot u_1^{2m^\ast_1-1})
\neq 0,\lb{5.5}\eea
where the last inequality in (\ref{5.5}) follows from (\ref{5.1}).

By Proposition 2.5, we have   $k_{l, \pm 1}(u_1)=k_{l, \pm 1}(u_1^2)$;
by Proposition 2.3 and Definition 2.4, we have $k_l(u_1^2)=k_{l, -1}(u_1^2)$ and $k_l(u_1)=k_{l, +1}(u_1)$
for $l\in\Z$. By (\ref{5.5}) and  Corollary 8.4 of \cite{MaW1}, $u_1$
is a local maximum of $\Psi_a$  in the local characteristic manifold $W(u_1)$
and then we have $k_{l, \pm 1}(u_1)=0$ for any $l\neq\nu(u_1)-1=2$
by Corollary 8.4 of \cite{MaW1} and Definition 2.4. Hence we have $k_l(u_1^2)=0$ for $l\neq 2$.
Then by (\ref{5.4}), we have
$C_{S^1,\; 2T-3-l}(\Psi_a, \;S^1\cdot u_1^{2m_1})=0$ for $l\neq 0$.
On the other hand, we have $C_{S^1,\;2T-2\xi_T^{-1}(1)}(\Psi_a,\;S^1\cdot u_1^{2m_1})\neq 0$
by (\ref{4.21}). This contradiction  proves the lemma. \hfill\hb

{\bf Lemma 5.3.} {\it If $(\tau_1, y_1)$ belongs to Case A in \S4 and  the matrix $M_1$ can be connected within $\Omega^0(M_1)$ to
 $N_1(1,\,-1)^{\diamond 2}\diamond R(\theta)$ with $\frac{\theta}{\pi}\notin\Q$,
 i.e.,  $R(\theta)$ is irrationally elliptic, 
 then we have $^\#\T(\Sg)\ge 4$.}

{\bf Proof. }  As in lemma 5.1, suppose $(T,\,m_1,\,m_2,\,m_3)$ and $(j_k,\,l_{j_k})$
 satisfy (\ref{4.15})-(\ref{4.21}). As mentioned in
Case A, we have  $i(y_1,\,1)=4$, thus by Theorems 3.6,  3.7 and (\ref{4.3}),
we have \bea i(y_1^m)&=&m(i(y_1,
1)+1-1)+2E\left(\frac{m\theta}{2\pi}\right)-2-4
=4m+2E\left(\frac{m\theta}{2\pi}\right)-6,\nn\\
\nu(y_1^m)&=&3,\quad\forall m\in\N.\lb{5.6}\eea
Hence we have $\hat i(y_1)=4+\frac{\theta}{\pi}\notin\Q$
and $\frac{\alpha_{1, 1}}{\hat i(y_1)}=\frac{\theta/\pi}{4+\theta/\pi}\notin\Q$,
where we use notations as in Theorem 3.9, i.e., $\alpha_{1, 1}=\frac{\theta}{\pi}$.
Denote by $\beta=4+\frac{\theta}{\pi}\notin\Q$. Then we have
$(\frac{1}{M\hat i(y_1)},\, \frac{\alpha_{1, 1}}{\hat i(y_1)})=(\frac{1}{M\beta}, \, 1-\frac{4}{\beta})$.
Thus if $\frac{T}{M\hat i(y_1)}=K+\mu$ for some $K\in\Z$ and $\mu\in (-1,\,1)$,
we have $\frac{T\alpha_{1, 1}}{\hat i(y_1)}=T-4MK-4M\mu$.
Hence by (\ref{3.37}) and  (\ref{3.39}), we have
\be
 \left\{\matrix{\chi_{1, 1}=1\quad{\rm if}\quad \chi_1=0, \cr
\chi_{1, 1}=0\quad{\rm if}\quad \chi_1=1. \cr}\right.  \lb{5.7}\ee
Thus either $(\chi_1,\,\chi_{1, 1})=(1,\, 0)$  or $(\chi_1,\,\chi_{1, 1})=(0,\, 1)$
holds. By (4.16)  and (4.17) of \cite{LoZ1},
we have
\bea
 &&\{m_1\alpha_{1, 1}\}=\left\{\left\{\frac{T\alpha_{1, 1}}{\hat i(y_1)}\right\}-\chi_{1, 1}
+\left(\chi_1-\left\{\frac{T}{M\hat i(y_1)}\right\}\right)M\alpha_{1, 1}\right\}
=\{A_{1, 1}(T)+B_{1, 1}(T)\}
\nn\\=&& \left\{\matrix{\left\{\left\{\frac{T\alpha_{1, 1}}{\hat i(y_1)}\right\}-\chi_{1, 1}
+\left(\chi_1-\left\{\frac{T}{M\hat i(y_1)}\right\}\right)M\alpha_{1, 1}\right\}
\quad{\rm if}\quad (\chi_1,\,\chi_{1, 1})=(1,\, 0), \cr
\left\{1+\left\{\frac{T\alpha_{1, 1}}{\hat i(y_1)}\right\}-\chi_{1, 1}
+\left(\chi_1-\left\{\frac{T}{M\hat i(y_1)}\right\}\right)M\alpha_{1, 1}\right\}
\quad{\rm if}\quad (\chi_1,\,\chi_{1, 1})=(0,\, 1), \cr}\right.
\nn\eea
where $A_{1, 1}(T)=\left\{\frac{T\alpha_{1, 1}}{\hat i(y_1)}\right\}-\chi_{1, 1}$
and $B_{1, 1}(T)=\left(\chi_1-\left\{\frac{T}{M\hat i(y_1)}\right\}\right)M\alpha_{1, 1}$.
In fact, we have  $A_{1, 1}(T)>0$,  $B_{1, 1}(T)>0$ for  $(\chi_1,\,\chi_{1, 1})=(1,\, 0)$,
and    $A_{1, 1}(T)<0$,  $B_{1, 1}(T)<0$ for  $(\chi_1,\,\chi_{1, 1})=(0,\, 1)$,
thus the last equality above holds.

Hence by (\ref{3.39}), we have
\bea\left\{\matrix{\{m_1\alpha_{1, 1}\}<(2M+1)\epsilon\quad&&{\rm if}\quad (\chi_1,\,\chi_{1, 1})=(1,\, 0), \cr
\{m_1\alpha_{1, 1}\}>1-(2M+1)\epsilon\quad&&{\rm if}\quad (\chi_1,\,\chi_{1, 1})=(0,\, 1), \cr}\right.
\lb{5.8}\eea
where we have used the fact that $\alpha_{1, 1}=\theta/\pi\in(0,\,2)$.
By choosing $\epsilon\in \left(0,\;\frac{1}{2M+1}\min\{\frac{\theta}{2\pi},\,1-\frac{\theta}{2\pi}\}\right)$,
 we have
\bea
 i(y_1^{2m_1+1})-i(y_1^{2m_1})=\left\{\matrix{4\quad&&{\rm if}\quad (\chi_1,\,\chi_{1, 1})=(1,\, 0), \cr
6\quad&&{\rm if}\quad (\chi_1,\,\chi_{1, 1})=(0,\, 1). \cr}\right.  \lb{5.9}\eea
In fact, by (\ref{5.6}), we have
\bea &&i(y_1^{2m_1+1})-i(y_1^{2m_1})\nn\\
=&&4(2m_1+1)+2E\left(\frac{(2m_1+1)\theta}{2\pi}\right)-6
-8m_1-2E\left(\frac{2m_1\theta}{2\pi}\right)+6\nn\\
=&&4+2\left(E\left(\frac{2m_1\theta}{2\pi}+\frac{\theta}{2\pi}\right)
-E\left(\frac{2m_1\theta}{2\pi}\right)\right)
\nn\\=&&4+2\left(E\left(\{m_1\alpha_{1, 1}\}+\frac{\theta}{2\pi}\right)
-E\left(\{m_1\alpha_{1, 1}\}\right)\right)
\nn\\=&&\left\{\matrix{4\quad{\rm if}\quad \{m_1\alpha_{1, 1}\}<(2M+1)\epsilon,\qquad \cr
6\quad{\rm if}\quad \{m_1\alpha_{1, 1}\}>1-(2M+1)\epsilon. \cr}\right.
\nn\eea
Hence (\ref{5.9}) holds by (\ref{5.8}).

Since $i(y_1^{2m_1+1})=2T$ by (\ref{4.8}), $i(y_1, 1)=4$ and Theorem 3.6, hence we have
\bea i(y_1^{2m_1})+\nu(y_1^{2m_1})-1=\left\{\matrix{2T-2\quad&&{\rm if}\quad (\chi_1,\,\chi_{1, 1})=(1,\, 0),\cr
2T-4\quad&&{\rm if}\quad (\chi_1,\,\chi_{1, 1})=(0,\, 1), \cr}\right.
\lb{5.10}\eea
by (\ref{5.6}) and (\ref{5.9}).

By Propositions 2.3  and  2,6, we have $K(y_1)=1$ and
 \bea &&\rank C_{S^1,\; i(y_1^{2m_1})+\nu(y_1^{2m_1})-1}(\Psi_a, \;S^1\cdot u_1^{2m_1})\nn\\
 =&&k_{\nu(u_1^{2m_1})-1}(u_1^{2m_1})
 =k_{\nu(u_1)-1}(u_1)=k_{\nu(u_1^{2m^\ast_1-1})-1}(u_1^{2m^\ast_1-1})\nn\\
 =&&\rank C_{S^1,\; 2T^\ast-8}(\Psi_a, \;S^1\cdot u_1^{2m^\ast_1-1})\neq 0,\lb{5.11}\eea
where the last inequality follows from (\ref{5.1}).
Thus we have
\bea\rank C_{S^1,\; l}(\Psi_a, \;S^1\cdot u_1^{2m_1})=0
\lb{5.12}\eea
for any integer $l\neq  i(y_1^{2m_1})+\nu(y_1^{2m_1})-1$ by Proposition 2.7.
Hence by (\ref{4.21}) and (\ref{5.10}), we have
\bea \Phi(u_1^{2m_1})=\left\{\matrix{c_T\quad&&{\rm if}\quad (\chi_1,\,\chi_{1, 1})=(1,\, 0),\cr
c_{T-1}\quad&&{\rm if}\quad (\chi_1,\,\chi_{1, 1})=(0,\, 1), \cr}\right.
\lb{5.13}\eea

For any $a\in A(v)$ fixed, denote by $(T, \,\chi(a))$, $m_k$ and $(T^\prime,\,\chi(-a))$,
$m_k^\prime$ as in Lemma 4.4.  Suppose $(T,\,m_1,\,m_2,\,m_3)$, $(j_k,\,l_{j_k})$
 and $(T^\prime,\,m_1^\prime,\,m_2^\prime,\,m_3^\prime)$, $(j_k^\prime,\,l_{j_k^\prime})$
 satisfy (\ref{4.15})-(\ref{4.21}).

We have the following two cases:

{\bf Case a.} {\it We have $(\chi(a)_1,\,\chi(a)_{1, 1})=(1,\, 0)$.}

By Theorems 3.10 and 3.11,  we have $(\chi(-a)_1,\,\chi(-a)_{1, 1})=(0,\, 1)$.
Then we have $c_T=\Phi(u_1^{2m_1})$ and $c_{T^\prime-1}=\Phi(u_1^{2m_1^\prime})$
by (\ref{5.13}). Thus by Lemma 4.2, we have
$$m_1\hat i(y_1)>\max\{m_2\hat i(y_2),\;m_3\hat i(y_3)\},\qquad
m_1^\prime\hat i(y_1)>m_{j_3^\prime}^\prime\hat i(y_{j_3^\prime}).$$
Since $\Phi(u_1^{2m_1^\prime})=c_{T^\prime-1}\neq c_{T^\prime-2}=\Phi(u_{j_3^\prime}^{2m^\prime_{j_3^\prime}})$, we have $j_3^\prime\neq 1$
by Proposition 2.9.
Hence by the same proof as the {\it Claim } in Lemma 4.4,
we have $a_1\hat i(y_1)=a_{j_3^\prime}\hat i(y_{j_3^\prime})$
with some $j_3^\prime\in\{2, 3\}$.

{\bf Case b.} {\it We have $(\chi(a)_1,\,\chi(a)_{1, 1})=(0,\, 1)$.}

By Theorems 3.10 and 3.11,  we have $(\chi(-a)_1,\,\chi(-a)_{1, 1})=(1,\, 0)$.
Then we have $c_{T-1}=\Phi(u_1^{2m_1})$ and $c_{T^\prime}=\Phi(u_1^{2m_1^\prime})$
by (\ref{5.13}). Thus by Lemma 4.2, we have
$$m_1\hat i(y_1)>m_{j_3}\hat i(y_{j_3}),\qquad
m_1^\prime\hat i(y_1)>\max\{m_2^\prime\hat i(y_2),\;m_3^\prime\hat i(y_3)\}.$$
Hence by the same proof as the {\it Claim } in Lemma 4.4,
we have $a_1\hat i(y_1)=a_{j_3}\hat i(y_{j_3})$
with some $j_3\in\{2, 3\}$.

Combining these two cases, at least one of the two equalities:
$a_1\hat i(y_1)=a_2\hat i(y_2)$,
 $a_1\hat i(y_1)=a_3\hat i(y_3)$ holds.

  Let $V_j=\{a\in V\,|\, a_1\hat i(y_1)=a_j\hat i(y_j)\}$  for $j=2, 3$ and
  \bea C(v)=A(v)\setminus \cup_{V_j\neq V,\,j=2, 3}V_j. \lb{5.14}\eea
 Since $\dim V\ge 1$  and  $A(v)$  is obtained from $V$ by deleting finitely
 many proper linear subspaces of $V$ by Theorem 3.11, and so is
 $C(v)$. Hence $C(v)$ is nonempty.
  Now we choose an $a\in C(v)$. By the above argument, we have
 $a\in V_2$ or $a\in V_3$. By the definition of $a \in C(v)$ we have
  $V_2 = V$ or $V_3=V$.

 Now by the same argument as in Lemma 4.4,  we have
$m_1\hat i(y_1)=m_2 \hat i(y_2)$ or $m_1\hat i(y_1)=m_3\hat i(y_3)$ holds.
 This contradict to Lemma 4.2, and then the lemma holds.\hfill\hb

In our study below, we need properties of sequences of vectors in $\R^n$ 
uniformly distributed mod one in number theory which can be found  in \cite{GrR}
or \S23.10 of \cite{HaW}.

{\bf Definition 5.4.} (cf. P.  5-6 of \cite{GrR}) {\it For given $v = (v_1, \ldots, v_n)\in\R^n$, define $v \;\mod\; 1$ to be the vector $\{v\} = (\{v_1\}, \ldots ,\{v_n\})$. 
The sequence of vectors $\{u_k\}_{k\in\N}$ with $u_k\in\R^n$ is {\bf uniformly distributed mod 
one} if for any $0\le b_j < c_j < 1$ for $ j = 1,2\ldots,n$, we have
$$\lim_{n\rightarrow\infty}\frac{1}{N} \;^\#\{k\le N|\{u_k\}\in\oplus_{j=1}^n[b_j,\;c_j)\}
=\Pi_{j=1}^n(c_j-b_j).$$}

{\bf Theorem 5.5.} (Kronecker's result, cf. P.  6 of \cite{GrR}) {\it 
If $1, v_1, \ldots , v_n$ are linearly independent over $\Q$, then the vectors 
$\{(kv_1, \ldots, kv_n)\}_{k\in\N}$ are uniformly distributed mod one on $[0,\; 1]^n$.}

{\bf Lemma 5.6.} {\it If  $(\tau_1, y_1)$ belongs to Case A in \S4 and  the matrix $M_1$ can be connected within $\Omega^0(M_1)$ to
 $N_1(1,\,-1)^{\diamond 2}\diamond R(\theta)$ with
 $\frac{\theta}{\pi}\in(0,\,2)\cap\Q$, i.e., $R(\theta)$ is rationally elliptic,  then we have $^\#\T(\Sg)\ge 4$.}

{\bf Proof. } As in Lemma 5.1, suppose $(T,\,m_1,\,m_2,\,m_3)$ and $(j_k,\,l_{j_k})$
 satisfy (\ref{4.15})-(\ref{4.21}). As mentioned in
Case A, we have  $i(y_1,\,1)=4$, thus by Theorems 3.6,  3.7 and
(\ref{4.3}), we have \bea &&i(y_1^m)=m(i(y_1,
1)+1-1)+2E\left(\frac{m\theta}{2\pi}\right)-2-4
=4m+2E\left(\frac{m\theta}{2\pi}\right)-6, \nn\\
&&\nu(y_1^m)=3+2-2\varphi\left(\frac{m\theta}{2\pi}\right),\qquad m\in\N.\lb{5.15}\eea
 By (\ref{4.18}), we have
$i(y_1^{2m_1-1})+\nu(y_1^{2m_1-1})-1=2T-8$. Hence we have
\bea i(y_1^{2m_1})=2T-6, \quad i(y_1^{2m_1})+\nu(y_1^{2m_1})-1=2T-2.\lb{5.16}\eea
In fact, by (\ref{3.36}), we have $\frac{m_1\theta}{\pi}\in\Z$,
this yields
\bea &&E\left(\frac{(2m_1-1)\theta}{2\pi}\right)
=E\left(\frac{2m_1\theta}{2\pi}-\frac{\theta}{2\pi}\right)
=E\left(\frac{2m_1\theta}{2\pi}\right),\nn\\
&&\varphi\left(\frac{(2m_1-1)\theta}{2\pi}\right)
=\varphi\left(\frac{2m_1\theta}{2\pi}-\frac{\theta}{2\pi}\right)=1,
\quad\varphi\left(\frac{2m_1\theta}{2\pi}\right)=0,\lb{5.17}\eea
since $\frac{\theta}{2\pi}\in(0,\,1)$. Clearly, (\ref{5.17})
implies (\ref{5.16}).

 Hence by Propositions 2.3
and  2.6, we have
 \bea &&\rank C_{S^1,\; 2T-2-l}(\Psi_a, \;S^1\cdot u_1^{2m_1})\nn\\
 =&&k_{\nu(u_1^{2m_1})-1-l}(u_1^{2m_1})
 =k_{\nu(u_1^{K(u_1)})-1-l}(u_1^{K(u_1)})
 =k_{\nu(u_1^{2m^\ast_1})-1-l}(u_1^{2m^\ast_1})\nn\\
 =&&\rank C_{S^1,\; 2T^\ast-2-l}(\Psi_a, \;S^1\cdot u_1^{2m_1^\ast}),\lb{5.18}\eea
for any $l\in\Z$, where in the second and third equality above, we have 
used the fact that $K(u_1)|2m_1$ and $K(u_1)|2m_1^\ast$, which follows from (\ref{3.36})
and Proposition 2.6.

By (\ref{4.21}), we have $ C_{S^1,\; 2T^\ast-2-l}(\Psi_a, \;S^1\cdot u_1^{2m_1^\ast})\neq 0$
for some $l\in\{0, 2, 4\}$. Thus we have the following three cases:

(i) If  $C_{S^1,\; 2T^\ast-2}(\Psi_a, \;S^1\cdot u_1^{2m_1^\ast})\neq 0$,
then we have $C_{S^1,\; 2T^\ast-2-l}(\Psi_a, \;S^1\cdot u_1^{2m_1^\ast})=0$
for $l\neq 0$ by Proposition 2.7. This implies
$C_{S^1,\; 2T-2-l}(\Psi_a, \;S^1\cdot u_1^{2m_1})=0$
for $l\neq 0$ by (\ref{5.18}). Hence by (\ref{4.21}), we have $c_T=\Phi(u_1^{2m_1})$,
and then $c_{T+1-\xi_T^{-1}(i)}=\Phi(u_i^{2m_i})$
for  $i=2, 3$ and $\xi_T^{-1}(i)\in\{2, 3\}$.
Thus  we have $\Phi(u_1^{2m_1})>\Phi(u_i^{2m_i})$
for $i=2, 3$.   This contradict to Proposition 4.5
and proves the lemma in this case.

(ii) If  $C_{S^1,\; 2T^\ast-6}(\Psi_a, \;S^1\cdot u_1^{2m_1^\ast})\neq 0$,
then we have $C_{S^1,\; 2T^\ast-6+l}(\Psi_a, \;S^1\cdot u_1^{2m_1^\ast})=0$
for $l\neq 0$ by Proposition 2.7. This implies
$C_{S^1,\; 2T-6+l}(\Psi_a, \;S^1\cdot u_1^{2m_1})=0$
for $l\neq 0$ by (\ref{5.18}).  Hence by (\ref{4.21}), we have $c_{T-2}=\Phi(u_1^{2m_1})$,
and then $c_{T+1-\xi_T^{-1}(i)}=\Phi(u_i^{2m_i})$
for  $i=2, 3$ and $\xi_T^{-1}(i)\in\{1, 2\}$.
Thus  we have $\Phi(u_1^{2m_1})<\Phi(u_i^{2m_i})$
for $i=2, 3$.   This contradict to Proposition 4.5
and proves the lemma in this case.

The following of Lemma 5.6 is devoted to study the following case:

(iii) If 
\bea C_{S^1,\; 2T^\ast-4}(\Psi_a, \;S^1\cdot u_1^{2m_1^\ast})\neq 0,\lb{5.19}\eea
then
\bea C_{S^1,\; 2T^\ast-2}(\Psi_a, \;S^1\cdot u_1^{2m_1^\ast})=0,
\quad C_{S^1,\; 2T^\ast-6}(\Psi_a, \;S^1\cdot u_1^{2m_1^\ast})=0,\lb{5.20}\eea
by Proposition 2.7.
This implies
\bea C_{S^1,\; 2T-2}(\Psi_a, \;S^1\cdot u_1^{2m_1})=0,
\quad C_{S^1,\; 2T-6}(\Psi_a, \;S^1\cdot u_1^{2m_1})=0,\lb{5.21}\eea
 by (\ref{5.18}). Hence we have $c_{T-1}=\Phi(u_1^{2m_1})$
by (\ref{4.21}), i.e., $\xi_T(2)=1$, and then we have $c_T=\Phi(u_{\xi_T(1)}^{2m_{\xi_T(1)}})$
and $c_{T-2}=\Phi(u_{\xi_T(3)}^{2m_{\xi_T(3)}})$ for $\xi_T(1), \xi_T(3)\in\{2, 3\}$
and $\xi_T(1)\neq \xi_T(3)$.

In order to prove Theorem 1.1 in this case, we
must study the properties of the critical modules
carefully, these properties are listed in the following five claims.

{\bf Claim 1.} {\it There exist two tuples $(T, m_1, m_2, m_3)$
and $(T^\prime, m_1^\prime, m_2^\prime, m_3^\prime)$
satisfying (\ref{4.15})-(\ref{4.21}) such that
$c_T=\Phi(u_2^{2m_2})$, $c_{T-2}=\Phi(u_3^{2m_3})$
and $c_{T^\prime}=\Phi(u_3^{2m_3^\prime})$,
$c_{T^\prime-2}=\Phi(u_2^{2m_2^\prime})$, i.e., 
we have the following diagram 
\bea
\begin{tabular}{lccccccc}
&$c_{T-2}$&$c_{T-1}$& $c_T$& $\qquad$&$c_{T^\prime-2}$&$c_{T^\prime-1}$& $c_{T^\prime}$\\
&$u_3^{2m_3}$&$u_1^{2m_1}$& $u_2^{2m_2}$& $\qquad$& $u_2^{2m_2^\prime}$&$u_1^{2m_1^\prime}$&$u_3^{2m_3^\prime}$\\
\end{tabular}\nn\eea
This implies that the orders of closed characteristics in the common index jump intervals 
have some commutative property. }

Suppose the contrary. Then we may assume $c_{T-2}=\Phi(u_2^{2m_2})$ for any tuple
$(T, m_1, m_2, m_3)$ satisfying (\ref{4.15})-(\ref{4.21}) without loss of generality.
Then by  (\ref{4.21}) and (\ref{5.21}), we have $c_{T-1}=\Phi(u_1^{2m_1})$ and $c_T=\Phi(u_3^{2m_3})$. This contradict to Proposition 4.5 and proves Claim 1.

{\bf Claim 2.} {\it The matrix $M_2, M_3$ can be connected within $\Omega^0(M_2),\,\Omega^0(M_3)$ to
 $R(\vartheta_1)\diamond R(\vartheta_2)\diamond M_2^\prime$
 and $R(\varphi_1)\diamond R(\varphi_2)\diamond M_3^\prime$
 with $\frac{\vartheta_i}{\pi},\, \frac{\varphi_i}{\pi}\notin\Q$
 for $i=1, 2$ and $M_2^\prime, M_3^\prime\in \Sp(2)$.
Moreover, $M_2^\prime, M_3^\prime\in\{I_2, N_1(1, -1), -I_2, N_1(-1, 1),
R(\vartheta)\}$. In fact,  in order to interchange  the orders of closed characteristics 
in the common index jump
intervals as in Claim 1, $M_2, M_3$ must have at least two irrational rotation components.  }

Firstly we prove $M_2$ has the required property.

We prove $M_2$ must have at least two irrational rotation components at first. 
In fact, by Claim 1 and (\ref{4.21}), we have
\bea C_{S^1,\; 2T-2}(\Psi_a, \;S^1\cdot u_2^{2m_2})\neq0,\quad
C_{S^1,\; 2T^\prime-6}(\Psi_a, \;S^1\cdot u_2^{2m_2^\prime})\neq0.\lb{5.22}\eea
for some tuples $(T, m_1, m_2, m_3)$ and
 $(T^\prime, m_1^\prime, m_2^\prime, m_3^\prime)$ satisfying (\ref{4.15})-(\ref{4.21}).
By (\ref{4.16}), (\ref{5.22}),  Propositions 2.3 and 2.7, we have
\bea i(y_2^{2m_2})+\nu(y_2^{2m_2})-1=2T-2,
\quad C_{S^1,\; 2T-2-l}(\Psi_a, \;S^1\cdot u_2^{2m_2})=0,\;\forall l\neq 0,\lb{5.23}\eea
i.e., $u_2^{2m_2}$ is a local maximum in the local characteristic manifold $W(u_2^{2m_2})$.
Note that  by (\ref{3.36}) and Proposition 2.6, we have
$K(u_2)|2m_2$ and $K(u_2)|2m_2^\prime$.
Hence by Propositions 2.3 and  2.6
 we have
 \bea &&\rank C_{S^1,\; 2T-2-l}(\Psi_a, \;S^1\cdot u_2^{2m_2})\nn\\
 =&&\rank C_{S^1,\; i(y_2^{2m_2})+\nu(y_2^{2m_2})-1-l}(\Psi_a, \;S^1\cdot u_2^{2m_2})\nn\\
 =&&k_{\nu(u_2^{2m_2})-1-l}(u_2^{2m_2})
 =k_{\nu(u_2^{K(u_2)})-1-l}(u_2^{K(u_2)})
 =k_{\nu(u_2^{2m^\prime_2})-1-l}(u_2^{2m^\prime_2})\nn\\
 =&&\rank C_{S^1,\; i(y_2^{2m^\prime_2})+\nu(y_2^{2m^\prime_2})-1-l}(\Psi_a, \;S^1\cdot u_2^{2m_2^\prime}),\lb{5.24}\eea
for any $l\in\Z$. Hence by (\ref{5.22}) and Proposition 2.7,
we have
\bea C_{S^1,\; i(y_2^{2m^\prime_2})+\nu(y_2^{2m^\prime_2})-1-l}(\Psi_a, \;S^1\cdot u_2^{2m_2^\prime})
=0, \qquad\forall l\neq 0.\lb{5.25}\eea
Hence by (\ref{5.22}), we have
\bea i(y_2^{2m^\prime_2})+\nu(y_2^{2m^\prime_2})-1=2T^\prime-6.\lb{5.26}\eea

By (\ref{4.16}), (\ref{5.22}) and Proposition 2.3,
we have $e(\gamma_2(\tau_2))=8$, i.e., $(\tau_2, y_2)$ is elliptic.
Assume $\gamma_2(\tau_2)$ can be connected within $\Omega^0(\gamma_2(\tau_2))$ to
\bea &&N_1(1,1)^{\diamond p_-} \diamond I_{2p_0}\diamond
N_1(1,-1)^{\diamond p_+}
\diamond N_1(-1,1)^{\diamond q_-} \diamond (-I_{2q_0})\diamond
N_1(-1,-1)^{\diamond q_+}\nn\\
&&\diamond R(\vartheta_1)\diamond\cdots\diamond R(\vartheta_r)
\diamond N_2(\omega_1, u_1)\diamond\cdots\diamond N_2(\omega_{r_*}, u_{r_*}) \nn\\
&&\diamond N_2(\lm_1, v_1)\diamond\cdots\diamond N_2(\lm_{r_0},
v_{r_0}), \lb{5.27}\eea where we use notations as in Theorem 3.7.
Then by (\ref{4.3}) and Theorem 3.7, we have \bea &&i(y_2,
2m_2+1)-(i(y_2, 2m_2)+\nu(y_2, 2m_2)-1) \nn\\=&&(2m_2+1)(i(y_2,\,
1)+p_-+p_0-r)+2\sum_{j=1}^r
E\left(\frac{(2m_2+1)\vartheta_j}{2\pi}\right)-r -p_--p_0\nn\\&&
\qquad-\frac{1+(-1)^{2m_2+1}}{2}(q_0+q_+)+2\left(
\sum_{j=1}^{r_*}\varphi\left(\frac{(2m_2+1)\alpha_j}{2\pi}\right)-r_*\right)\nn\\
&&-2m_2(i(y_2,\, 1)+p_-+p_0-r)-2\sum_{j=1}^r
E\left(\frac{2m_2\vartheta_j}{2\pi}\right)+r +p_-+p_0\nn\\&&
\qquad+\frac{1+(-1)^{2m_2}}{2}(q_0+q_+)-2\left(
\sum_{j=1}^{r_*}\varphi\left(\frac{2m_2\alpha_j}{2\pi}\right)-r_*\right)\nn\\
&&-\nu(y_2,\,
1)-\frac{1+(-1)^{2m_2}}{2}(q_-+2q_0+q_+)-2(r+r_*+r_0)\nn\\
&&\qquad+2\left(\sum_{j=1}^{r}\varphi\left(\frac{2m_2\vartheta_j}{2\pi}\right)+
\sum_{j=1}^{r_*}\varphi\left(\frac{2m_2\alpha_j}{2\pi}\right)
+\sum_{j=1}^{r_0}\varphi\left(\frac{2m_2\beta_j}{2\pi}\right)\right)+1\nn\\
=&&i(y_2,\, 1)-p_0-p_+-q_--q_0-r\nn\\
\qquad &&+2\sum_{j=1}^r \left(E\left(\frac{(2m_2+1)\vartheta_j}{2\pi}\right)
-E\left(\frac{2m_2\vartheta_j}{2\pi}\right)
+\varphi\left(\frac{2m_2\vartheta_j}{2\pi}\right)-1\right)\nn\\
\qquad &&+2\sum_{j=1}^{r_0}\left(\varphi\left(\frac{2m_2\beta_j}{2\pi}\right)-1\right)+1\nn\\
=&&i(y_2,\, 1)-p_0-p_+-q_--q_0-r\nn\\
\qquad &&+2\sum_{1\le j\le r,\,\vartheta_j/\pi\notin\Q} \left(E\left(\frac{(2m_2+1)\vartheta_j}{2\pi}\right)
-E\left(\frac{2m_2\vartheta_j}{2\pi}\right)\right)\nn\\
\qquad &&-2^\#\{\beta_j|1\le j\le r_0,\,\beta_j/\pi\in\Q\}+1\nn\\
=&&2\sum_{1\le j\le r,\,\vartheta_j/\pi\notin\Q}
\left(E\left(\frac{(2m_2+1)\vartheta_j}{2\pi}\right)
-E\left(\frac{2m_2\vartheta_j}{2\pi}\right)\right)+\Xi,
\lb{5.28}\eea where $\Xi$ is independent of $m_2$. Here in the
second equality, we have used the fact that $\nu(y_2,\,
1)=p_-+2p_0+p_+$ and
$\varphi\left(\frac{(2m_2+1)\alpha_j}{2\pi}\right)=1$. In fact, if
$\frac{\alpha_j}{\pi}\notin\Q$, then
$\frac{(2m_2+1)\alpha_j}{2\pi}\notin\Z$. If
$\frac{\alpha_j}{\pi}\in\Q$, then $\frac{m_2\alpha_j}{\pi}\in\Z$ by
(\ref{3.36}), and then $\frac{(2m_2+1)\alpha_j}{2\pi}\notin\Z$ since $\alpha_j\in (0,\,2\pi)$. Thus
we have $\varphi\left(\frac{(2m_2+1)\alpha_j}{2\pi}\right)=1$ by
(\ref{1.11}). In the third equality, we have used the fact that if
$\frac{\vartheta_j}{\pi}\in\Q$, then
$\frac{m_2\vartheta_j}{\pi}\in\Z$ by (\ref{3.36}) and then
$\varphi(\frac{m_2\vartheta_j}{\pi})=0$ together with \bea
E\left(\frac{(2m_2+1)\vartheta_j}{2\pi}\right)=
E\left(\frac{2m_2\vartheta_j}{2\pi}+\frac{\vartheta_j}{2\pi}\right)
=E\left(\frac{2m_2\vartheta_j}{2\pi}\right)+1, \lb{5.29}\eea and
$\frac{m_2\beta_j}{\pi}\in\Z$ by (\ref{3.36}) whenever
$\frac{\beta_j}{\pi}\in\Q$; $\frac{m_2\beta_j}{\pi}\notin\Z$
whenever $\frac{\beta_j}{\pi}\notin\Q$.

By (\ref{4.8}),  (\ref{5.23}), (\ref{5.28}) and Theorem
3.6, we have \bea &&2T+i(y_2,\,
1)-n-(2T-2)=i(y_2^{2m_2+1})-
(i(y_2^{2m_2})+\nu(y_2^{2m_2})-1)\nn\\
=&&2\sum_{1\le j\le r,\,\vartheta_j/\pi\notin\Q} \left(E\left(\frac{(2m_2+1)\vartheta_j}{2\pi}\right)
-E\left(\frac{2m_2\vartheta_j}{2\pi}\right)\right)+\Xi\ge \Xi.
\lb{5.30}\eea
Hence by (\ref{4.8}),  (\ref{5.28}), (\ref{5.30}) and Theorem 3.6, we have
\bea &&i(y_2^{2m^\prime_2})+\nu(y_2^{2m^\prime_2})-1
\nn\\=&&i(y_2^{2m^\prime_2+1})-2\sum_{1\le j\le r,\,\vartheta_j/\pi\notin\Q} \left(E\left(\frac{(2m^\prime_2+1)\vartheta_j}{2\pi}\right)
-E\left(\frac{2m^\prime_2\vartheta_j}{2\pi}\right)\right)-\Xi
\nn\\\ge&&2T^\prime+i(y_2,\, 1)-n-2^\#\{j|1\le j\le r,\,\vartheta_j/\pi\notin\Q\}\nn\\
&&-(2T^\prime+i(y_2,\, 1)-n-(2T^\prime-2))\nn\\
\ge&& 2T^\prime-2-2^\#\{j|1\le j\le r,\,\vartheta_j/\pi\notin\Q\},
\lb{5.31}\eea where $n=4$. Here in the second inequality, we have used the fact
that $i(y_2^{2m_2^\prime+1})=2T^\prime+i(y_2,\, 1)-n$; \bea
E\left(\frac{(2m_2^\prime+1)\vartheta_j}{2\pi}\right)=
E\left(\frac{2m_2^\prime\vartheta_j}{2\pi}+\frac{\vartheta_j}{2\pi}\right)
\le E\left(\frac{2m_2^\prime\vartheta_j}{2\pi}\right)+1; \nn\eea and
\bea 2T^\prime+i(y_2,\, 1)-n-(2T^\prime-2)=2T+i(y_2,\,
1)-n-(2T-2)\ge\Xi. \nn\eea Suppose $^\#\{j|1\le j\le
r,\,\vartheta_j/\pi\notin\Q\}\le 1$, then by (\ref{5.31}), we have \bea
i(y_2^{2m_2^\prime})+\nu(y_2^{2m_2^\prime})-1 \ge
2T^\prime-4.\lb{5.32}\eea This contradict to (\ref{5.26}) and proves
 $^\#\{j|1\le j\le r,\,\vartheta_j/\pi\notin\Q\}\ge 2$.
Hence the matrix $M_2$ can be connected within $\Omega^0(M_2)$
 to $R(\vartheta_1)\diamond R(\vartheta_2)\diamond M_2^\prime$
 with $\frac{\vartheta_i}{\pi}\notin\Q$
 for $i=1, 2$ and $M_2^\prime\in \Sp(2)$ is one of the basic normal
 form in (\ref{3.12})-(\ref{3.14}). Thus in order to prove Claim 2, it is sufficient to
 show that $M_2^\prime\notin\{D(\lambda), N_1(-1, -1), N_1(1, 1)\}$.

As mentioned above, $(\tau_2, y_2)$ is elliptic, hence $M_2^\prime$
is elliptic too, thus we have $M_2^\prime\neq D(\lambda)$. 

 Suppose $M_2^\prime=N_1(-1, -1)$, then by
(\ref{4.3}) and Theorems 3.5-3.7, we have \bea
i(y_2^{2m_2})&=&2m_2(i(y_2,\, 1)+1-2)+2\sum_{j=1}^2
E\left(\frac{2m_2\vartheta_j}{2\pi}\right)-2-1
-\frac{1+(-1)^{2m_2}}{2}-4 \nn\\&=&2m_2(i(y_2,\, 1)-1)+2\sum_{j=1}^2
E\left(\frac{2m_2\vartheta_j}{2\pi}\right)-8, \nn\\
\nu(y_2^{2m_2})&=&2,\lb{5.33} \eea with $i(y_2,\, 1)\in 2\N$. Hence
we have $i(y_2^{2m_2})+\nu(y_2^{2m_2})-1\in 2\N-1$. This contradict
to (\ref{5.23}) and shows that $M_2^\prime\neq N_1(-1, -1)$.

Suppose $M_2^\prime=N_1(1, 1)$, then by (\ref{4.3}) and Theorems 3.5-3.7,, we
have \bea i(y_2^{2m_2})&=&2m_2(i(y_2,\, 1)+1+1-2)+2\sum_{j=1}^2
E\left(\frac{2m_2\vartheta_j}{2\pi}\right)-2-1 -1-4
\nn\\&=&2m_2i(y_2,\, 1)+2\sum_{j=1}^2
E\left(\frac{2m_2\vartheta_j}{2\pi}\right)-8, \nn\\
\nu(y_2^{2m_2})&=&2,\lb{5.34} \eea with $i(y_2,\, 1)\in 2\N$. Hence
we have $i(y_2^{2m_2})+\nu(y_2^{2m_2})-1\in 2\N+1$. This contradict
to (\ref{5.23}) and shows that $M_2^\prime\neq N_1(1, 1)$. Hence $(\tau_2, y_2)$ has the required property.

Applying the above  argument to
\bea C_{S^1,\; 2T-6}(\Psi_a, \;S^1\cdot u_3^{2m_3})\neq0,\quad
C_{S^1,\; 2T^\prime-2}(\Psi_a, \;S^1\cdot u_3^{2m_3^\prime})\neq0.\lb{5.35}\eea
yields $(\tau_3, y_3)$ has the required property.
This proves Claim 2.

{\bf Claim 3.} {\it We have $C_{S^1,\,2k+1}(\Psi_a,\,S^1\cdot u_j^m)=0$
for $k\in\Z$, $m\in\N$ and $j=2, 3$. This implies the critical
modules of iterations of both $(\tau_2, y_2)$ and $(\tau_3, y_3)$ 
have no contribution to the number 
$$M_{2k+1}=\sum_{1\le j\le 3,\,m\in\N}
\rank C_{S^1,\,2k+1}(\Psi_a,\,S^1\cdot u_j^m).$$}

Note that by Theorem 3.7, the index iteration
formula of $I_2$ can be viewed
as that of a rotation matrix $R(\theta)$ with $\theta=2\pi$.
Similarly  $-I_2$  can be viewed
as a rotation matrix $R(\theta)$ with $\theta=\pi$.
Hence in the following, we will handle them together.

Firstly we prove $(\tau_2, y_2)$ has the required property, $(\tau_3, y_3)$
can be proved similarly. 

Due to Claim 2 above, we need to consider $M_2^\prime$ belonging
to one of the following cases:

If $M_2^\prime=R(\vartheta_3)$ with $\vartheta_3\in(0,\, 2\pi]$. By
(\ref{4.3}) and  Theorems 3.5-3.7, we have \bea
i(y_2^m)&=&m(i(y_2,\, 1)+1-3)+2\sum_{j=1}^3
E\left(\frac{m\vartheta_j}{2\pi}\right)-3-1-4 \nn\\&=&m(i(y_2,\,
1)-2)+2\sum_{j=1}^3 E\left(\frac{m\vartheta_j}{2\pi}\right)-8, \nn\\
\nu(y_2^m)&=&3-2\varphi\left(\frac{m\vartheta_3}{2\pi}\right),
\qquad m\in\N,\lb{5.36} \eea with $i(y_2,\, 1)\in 2\N$ and $i(y_2,\, 1)\ge 4$.

If $\vartheta_3/\pi\notin\Q$, then we have $\nu(y_2^m)=1$
and $i(y_2^m)$ is always even for $m\in\N$. Hence we have
$C_{S^1,\,p}(\Psi_a,\,S^1\cdot u_2^m)=0$
for $p\neq i(y_2^m)$ by Proposition 2.3. This implies
$C_{S^1,\,2k+1}(\Psi_a,\,S^1\cdot u_2^m)=0$
for $k\in\Z$ and $m\in\N$.

If $\vartheta_3/\pi\in\Q$, then we have $\nu(y_2^m)=3$
for $K(y_2)|m$ and $\nu(y_2^m)=1$ otherwise;
and $i(y_2^m)$ is always even for $m\in\N$.
Hence as above we have
$C_{S^1,\,2k+1}(\Psi_a,\,S^1\cdot u_2^m)=0$
for $k\in\Z$ and $m\notin K(y_2)\N$ by Proposition 2.3.

By Propositions 2.3
and  2.6, we have
 \bea &&\rank C_{S^1,\; i(y_2^{2m_2})+\nu(y_2^{2m_2})-1-l}(\Psi_a, \;S^1\cdot u_2^{2m_2})\nn\\
 =&&k_{\nu(u_2^{2m_2})-1-l}(u_2^{2m_2})
 =k_{\nu(u_2^{K(u_2)})-1-l}(u_2^{K(u_2)})
 =k_{\nu(u_2^m)-1-l}(u_2^m)\nn\\
 =&&\rank C_{S^1,\; i(y_2^m)+\nu(y_2^m)-1-l}(\Psi_a, \;S^1\cdot u_2^m),
 \quad \forall l\in\Z,\lb{5.37}\eea
for $K(y_2)|m$. Hence by (\ref{5.22}), (\ref{5.23}), (\ref{5.37}) and Proposition 2.7,
we have
\bea C_{S^1,\; i(y_2^m)+\nu(y_2^m)-1-l}(\Psi_a, \;S^1\cdot u_2^m)
=0, \qquad\forall l\neq 0,\lb{5.38}\eea
for $K(y_2)|m$. Hence we have
$C_{S^1,\,2k+1}(\Psi_a,\,S^1\cdot u_2^m)=0$
for $k\in\Z$ and $m\in K(y_2)\N$.

If $M_2^\prime=N_1(-1, 1)$, then by (\ref{4.3}) and  Theorems 3.5-3.7,, we have
\bea i(y_2^m)&=&m(i(y_2,\, 1)+1-2)+2\sum_{j=1}^2
E\left(\frac{m\vartheta_j}{2\pi}\right)-2-1-4 \nn\\&=&m(i(y_2,\,
1)-1)+2\sum_{j=1}^2 E\left(\frac{m\vartheta_j}{2\pi}\right)-7, \nn\\
\nu(y_2^m)&=&1+\frac{1+(-1)^m}{2}, \qquad m\in\N,\lb{5.39} \eea with
$i(y_2,\, 1)\in 2\N$ and $i(y_2, \,1)\ge 4$. Then we have $\nu(y_2^m)=1$ for $m\in 2\N-1$
and $\nu(y_2^m)=2$ for $m\in 2\N$; $i(y_2^m)$ is even for
$m\in2\N-1$ and $i(y_2^m)$ is odd for $m\in2\N$. Hence as above we have
$C_{S^1,\,2k+1}(\Psi_a,\,S^1\cdot u_2^m)=0$ for $k\in\Z$ and
$m\in2\N-1$ by Proposition 2.3.

By Propositions 2.3 and  2.6, we have
 \bea &&\rank C_{S^1,\; i(y_2^{2m_2})+\nu(y_2^{2m_2})-1-l}(\Psi_a, \;S^1\cdot u_2^{2m_2})\nn\\
 =&&k_{\nu(u_2^{2m_2})-1-l}(u_2^{2m_2})
 =k_{\nu(u_2^2)-1-l}(u_2^2)
 =k_{\nu(u_2^{2m})-1-l}(u_2^{2m})\nn\\
 =&&\rank C_{S^1,\; i(y_2^{2m})+\nu(y_2^{2m})-1-l}(\Psi_a, \;S^1\cdot u_2^{2m}),
 \qquad \forall l\in\Z,\; m\in\N.\lb{5.40}\eea
Hence by (\ref{5.22}), (\ref{5.23}), (\ref{5.40}) and Proposition 2.7,
we have
\bea C_{S^1,\; i(y_2^{2m})+\nu(y_2^{2m})-1-l}(\Psi_a, \;S^1\cdot u_2^{2m})
=0, \qquad\forall l\neq 0,\; m\in\N.\lb{5.41}\eea
 Hence we have
$C_{S^1,\,2k+1}(\Psi_a,\,S^1\cdot u_2^{2m})=0$
for $k\in\Z$ and $m\in\N$.

If $M_2^\prime=N_1(1, -1)$, then by (\ref{4.3}) and  Theorems 3.5-3.7,, we have
\bea i(y_2^m)&=&m(i(y_2,\, 1)+1-2)+2\sum_{j=1}^2
E\left(\frac{m\vartheta_j}{2\pi}\right)-2-1-4 \nn\\&=&m(i(y_2,\,
1)-1)+2\sum_{j=1}^2 E\left(\frac{m\vartheta_j}{2\pi}\right)-7, \nn\\
\nu(y_2^m)&=&2, \qquad m\in\N,\lb{5.42} \eea with $i(y_2,\, 1)\in
2\N-1$ and  $i(y_2, \,1)\ge 4$. Then $i(y_2^m)$ is always odd for $m\in\N$. By Propositions
2.3 and  2.6, we have
 \bea &&\rank C_{S^1,\; i(y_2^{2m_2})+\nu(y_2^{2m_2})-1-l}(\Psi_a, \;S^1\cdot u_2^{2m_2})\nn\\
 =&&k_{\nu(u_2^{2m_2})-1-l}(u_2^{2m_2})
 =k_{\nu(u_2)-1-l}(u_2)
 =k_{\nu(u_2^m)-1-l}(u_2^m)\nn\\
 =&&\rank C_{S^1,\; i(y_2^m)+\nu(y_2^m)-1-l}(\Psi_a, \;S^1\cdot u_2^m),
 \qquad \forall l\in\Z,\; m\in\N.\lb{5.43}\eea
Hence by (\ref{5.22}), (\ref{5.23}), (\ref{5.43}) and Proposition 2.7,
we have
\bea C_{S^1,\; i(y_2^{m})+\nu(y_2^{m})-1-l}(\Psi_a, \;S^1\cdot u_2^{m})
=0, \qquad\forall l\neq 0,\;m\in\N.\lb{5.44}\eea
 Hence we have
$C_{S^1,\,2k+1}(\Psi_a,\,S^1\cdot u_2^m)=0$
for $k\in\Z$ and $m\in\N$.
This proves Claim 3.

{\bf Claim 4. } {\it We have $\sum_{i\in\Z}(-1)^i\rank C_{S^1,\,i}(\Psi_a, \,S^1\cdot u_1^m)=1$
for any $m\in\N$. This implies that the critical modules of iterations of
$(\tau_1, y_1)$ behave like those of a non-degenerate
critical point in the sense that the alternative sum of their ranks
is $1$.
}

Write $\frac{\theta}{2\pi}=\frac{r}{s}$ with $r, s\in\N$ and $(r,\, s)=1$.
Then $K(y_1)=s$ by (\ref{5.15}), and then  $\nu(y_1^m)=5$ for $K(y_1)|m$ and $\nu(y_1^m)=3$ otherwise.
Since $\frac{\theta}{2\pi}\in(0,\;1)$, we have $K(y_1)=s\ge 2$.
Then we have the following two cases:

(i) If $m\notin K(y_1)\N$, then we have
\bea &&\rank C_{S^1,\; i(y_1^m)+\nu(y_1^m)-1}(\Psi_a, \;S^1\cdot u_1^m)\nn\\
 =&&k_{\nu(u_1^m)-1}(u_1^m)
 =k_{\nu(u_1)-1}(u_1)=k_{\nu(u_1^{2m^\ast_1-1})-1}(u_1^{2m^\ast_1-1})\nn\\
 =&&\rank C_{S^1,\; 2T^\ast-8}(\Psi_a, \;S^1\cdot u_1^{2m^\ast_1-1})=1,\lb{5.45}\eea
by Propositions 2.3, 2.5-2.7, (\ref{4.18})  and (\ref{5.1}).
Hence
\bea\rank C_{S^1,\;  i(y_1^m)+\nu(y_1^m)-1-l}(\Psi_a, \;S^1\cdot u_1^m)
=k_{\nu(u_1^m)-1-l}(u_1^m)=0,\quad\forall l\neq 0,
\lb{5.46}\eea
by Proposition 2.7. Since $i(y_1^m)+\nu(y_1^m)-1$
is even by (\ref{5.15}), Claim 4 holds in this case.

(ii) If  $m\in K(y_1)\N$. then by Propositions 2.3
and  2.6, we have
 \bea &&\rank C_{S^1,\; i(y_1^{2m_1})+\nu(y_1^{2m_1})-1-l}(\Psi_a, \;S^1\cdot u_1^{2m_1})\nn\\
 =&&k_{\nu(u_1^{2m_1})-1-l}(u_1^{2m_1})
 =k_{\nu(u_1^{K(u_1)})-1-l}(u_1^{K(u_1)})
 =k_{\nu(u_1^m)-1-l}(u_1^m)\nn\\
 =&&\rank C_{S^1,\; i(y_1^m)+\nu(y_1^m)-1-l}(\Psi_a, \;S^1\cdot u_1^m),
 \quad \forall l\in\Z.\nn\eea
Thus in order to prove Claim 4 in this case, it is sufficient to prove Claim 4 for $m=2m_1$.

By (\ref{4.16}), (\ref{5.35}) and Proposition 2.3, we have
\bea i(y_3^{2m_3^\prime})+\nu(y_3^{2m_3^\prime})-1=2T^\prime-2,\lb{5.47}\eea
i.e., $u_3^{2m_3^\prime}$ is a local maximum in the local  characteristic manifold
$W(u_3^{2m_3^\prime})$.

As in Claim 2, we have
$K(u_3)|2m_3$ and $K(u_3)|2m_3^\prime$.
Hence by Propositions 2.3, 2.6 and (\ref{5.47}), 
 we have
 \bea &&\rank C_{S^1,\; 2T^\prime-2-l}(\Psi_a, \;S^1\cdot u_3^{2m^\prime_3})\nn\\
 =&&\rank C_{S^1,\; i(y_3^{2m_3^\prime})+\nu(y_3^{2m_3^\prime})-1-l}(\Psi_a, \;S^1\cdot u_3^{2m_3^\prime})\nn\\
 =&&k_{\nu(u_3^{2m_3^\prime})-1-l}(u_3^{2m_3^\prime})
 =k_{\nu(u_3^{K(u_3)})-1-l}(u_3^{K(u_3)})
 =k_{\nu(u_3^{2m_3})-1-l}(u_3^{2m_3})\nn\\
 =&&\rank C_{S^1,\; i(y_3^{2m_3})+\nu(y_3^{2m_3})-1-l}(\Psi_a, \;S^1\cdot u_3^{2m_3}),\lb{5.48}\eea
for any $l\in\Z$. Thus by (\ref{5.35}), (\ref{5.48}) and Proposition 2.7,
we have
\bea C_{S^1,\; i(y_3^{2m_3})+\nu(y_3^{2m_3})-1-l}(\Psi_a, \;S^1\cdot u_3^{2m_3})
=0, \qquad\forall l\neq 0.\lb{5.49}\eea
Hence by (\ref{5.35}), we have
\bea i(y_3^{2m_3})+\nu(y_3^{2m_3})-1=2T-6.\lb{5.50}\eea
Thus by (\ref{5.35}), (\ref{5.50}) and Proposition 2.7, we have
\bea  \rank C_{S^1,\; 2T-6}(\Psi_a, \;S^1\cdot u_3^{2m_3})=1.
 \lb{5.51}\eea
Hence we have
\bea &&M_{2T-6}=\sum_{1\le j\le 3,\,p\in\N}\rank C_{S^1,\; 2T-6}(\Psi_a, \;S^1\cdot u_j^p)
\nn\\&&\qquad=\rank C_{S^1,\; 2T-6}(\Psi_a, \;S^1\cdot u_3^{2m_3})=1=b_{2T-6},
\lb{5.52}\\&&M_{2T-2}=\sum_{1\le j\le 3,\,p\in\N}\rank C_{S^1,\; 2T-2}(\Psi_a, \;S^1\cdot u_j^p)
\nn\\&&\qquad=\rank C_{S^1,\; 2T-2}(\Psi_a, \;S^1\cdot u_2^{2m_2})=1=b_{2T-2},\lb{5.53}
\eea
In fact, the second equality in (\ref{5.52})
follows from (\ref{4.17})-(\ref{4.20}), (\ref{5.21}),
(\ref{5.23}) and Proposition 2.3; the second equality in (\ref{5.53})
follows from (\ref{4.17})-(\ref{4.20}), (\ref{5.21}),
(\ref{5.50}) and Proposition 2.3.  The third equality in (\ref{5.53}) follows  from (\ref{5.22}), (\ref{5.23}) and Proposition 2.7.  The last equalities
in (\ref{5.52}) and (\ref{5.53}) follows from Theorem 2.8.

Hence by Theorem 2.8, we have
\bea M_{2T-2}-M_{2T-3}+\cdots-M_1+M_0
&\ge&b_{2T-2}-b_{2T-3}+\cdots-b_1+b_0,\lb{5.54}\\
M_{2T-3}-M_{2T-4}+\cdots+M_1-M_0
&\ge&b_{2T-3}-b_{2T-4}+\cdots+b_1-b_0,\lb{5.55}\\
M_{2T-6}-M_{2T-7}+\cdots-M_1+M_0
&\ge&b_{2T-6}-b_{2T-7}+\cdots-b_1+b_0,\lb{5.56}\\
M_{2T-7}-M_{2T-8}+\cdots+M_1-M_0
&\ge&b_{2T-7}-b_{2T-8}+\cdots+b_1-b_0,\lb{5.57}
\eea
Thus from (\ref{5.52}) and (\ref{5.53}), we have
\bea
M_{2T-3}-M_{2T-4}+\cdots+M_1-M_0
&=&b_{2T-3}-b_{2T-4}+\cdots+b_1-b_0,\lb{5.58}\\
M_{2T-6}-M_{2T-7}+\cdots-M_1+M_0
&=&b_{2T-6}-b_{2T-7}+\cdots-b_1+b_0,\lb{5.59}
\eea
Adding (\ref{5.58}) and (\ref{5.59}), then by Theorem 2.8, we have
\bea &&-1=b_{2T-3}-b_{2T-4}+b_{2T-5}=M_{2T-3}-M_{2T-4}+M_{2T-5}
\nn\\=&&\sum_{2T-3\le q\le 2T-5\atop 1\le j\le 3,\,m\in\N}(-1)^{q+1}\rank C_{S^1,\; q}(\Psi_a, \;S^1\cdot u_j^m)
\nn\\=&&\sum_{2T-3\le q\le 2T-5}(-1)^{q+1}\rank C_{S^1,\; q}(\Psi_a, \;S^1\cdot u_1^{2m_1})
\nn\\=&&\sum_{q\in\Z}(-1)^{q+1}\rank C_{S^1,\; q}(\Psi_a, \;S^1\cdot u_1^{2m_1}).
\lb{5.60}\eea
Here in the next to the last equality, we have used
 (\ref{4.17})-(\ref{4.20}), (\ref{5.23}), (\ref{5.50}) and
 Proposition 2.3. In the last equality, we have used (\ref{5.16}), (\ref{5.21})
 and Proposition 2.3.
This proves Claim 4.

{\bf Claim 5. } {\it It is impossible that
 $C_{S^1,\,2K}(\Psi_a, \,S^1\cdot u_1^m)\neq 0$
 and $C_{S^1,\,2K}(\Psi_a, \,S^1\cdot u_j^k)\neq 0$
hold simultaneously for some $K, m, k\in\N$ and some $j\in\{2,\, 3\}$. 
This implies that the critical modules of iterations of $(\tau_1, y_1)$
and $(\tau_j, y_j)$ for $j\in\{2, 3\}$ can not hit together.}

Suppose there exist some $K, m, k\in\N$ and $j\in\{2,\, 3\}$
such that  $C_{S^1,\,2K}(\Psi_a, \,S^1\cdot u_1^m)\neq 0$
 and $C_{S^1,\,2K}(\Psi_a, \,S^1\cdot u_j^k)\neq 0$ hold simultaneously.

By Proposition 2.3,
we have $i(y_1^m)\le 2K\le i(y_1^m)+\nu(y_1^m)-1$.
By (\ref{4.10}) and Theorem 3.6, we have
\bea i(y_1^{l+1})\ge i(y_1^l)+\nu(y_1^l)+1,\quad
i(y_1^{l-1})+\nu(y_1^{l-1})-1\le i(y_1^l)-2,\lb{5.61}\eea
for any integer $l\ge 2$.
Hence we have
\bea C_{S^1,\, i(y_1^m)-1}(\Psi_a,\, S^1\cdot u_1^l)=0,
\quad C_{S^1,\, i(y_1^m)+\nu(y_1^m)}(\Psi_a,\, S^1\cdot u_1^l)=0,\lb{5.62}\eea
for any $l\in\N$ by Proposition 2.3.
In fact, if $l=m$, then (\ref{5.62}) holds directly from Proposition 2.3.
If  $l>m$, then $i(y_1^l)>i(y_1^m)+\nu(y_1^m)$ by (\ref{5.61}), hence (\ref{5.62}) holds  from Proposition 2.3.
If $l<m$, then $i(y_1^l)+\nu(y_1^l)-1<i(y_1^m)-1$ by (\ref{5.61}), hence (\ref{5.62}) holds  from Proposition 2.3.
By (\ref{5.15}), we have $i(y_1^m),\, i(y_1^m)+\nu(y_1^m)-1\in 2\N$,
hence by Claim 3, (\ref{5.62}) and Theorem 2.8, we have
\bea &&M_{i(y_1^m)-1}=\sum_{1\le i\le 3,\,l\in\N}\rank C_{S^1,\; i(y_1^m)-1}(\Psi_a, \;S^1\cdot u_i^l)=0=b_{i(y_1^m)-1},\nn\\
&&M_{i(y_1^m)+\nu(y_1^m)}=\sum_{1\le i\le 3,\,l\in\N}\rank C_{S^1,\; i(y_1^m)+\nu(y_1^m)}(\Psi_a, \;S^1\cdot u_i^l)=0=b_{i(y_1^m)+\nu(y_1^m)}.
\lb{5.63}\eea
Thus by (\ref{5.63}), Theorem 2.8 and the proof of from (\ref{5.52})-(\ref{5.53})-(\ref{5.57})
to (\ref{5.58})-(\ref{5.60}), we have
\bea &&\frac{\nu(y^m_1)-1}{2}+1=\sum_{i(y_1^m)\le q \le i(y_1^m)+\nu(y_1^m)-1}(-1)^qb_q
=\sum_{i(y_1^m)\le q \le i(y_1^m)+\nu(y_1^m)-1}(-1)^qM_q
\nn\\=&&\sum_{i(y_1^m)\le q \le i(y_1^m)+\nu(y_1^m)-1\atop 1\le i\le 3,\,l\in\N}
(-1)^q\rank C_{S^1,\; q}(\Psi_a, \;S^1\cdot u_i^l)
\nn\\=&&\sum_{i(y_1^m)\le q \le i(y_1^m)+\nu(y_1^m)-1}
(-1)^q\rank C_{S^1,\; q}(\Psi_a, \;S^1\cdot u_1^m)
\nn\\&&\qquad+\sum_{i(y_1^m)\le q \le i(y_1^m)+\nu(y_1^m)-1\atop 2\le i\le 3,\,l\in\N}
(-1)^q\rank C_{S^1,\; q}(\Psi_a, \;S^1\cdot u_i^l),
\lb{5.64}\eea
where the last equality follows from (\ref{5.61}) and Proposition 2.3.
By (\ref{5.15}), we have $\nu(y_1^m)\le 5$,
thus we have
\bea C_{S^1,\; 2p}(\Psi_a, \;S^1\cdot u_1^m)=0,\quad \forall 2p\in[i(y_1^m),\;i(y_1^m)+\nu(y_1^m)-1]\setminus\{2K\}.
\lb{5.65}\eea
In fact, only one of the following possible cases holds: $2K=i(y_1^m)$,
$i(y_1^m)<2K<i(y_1^m)+\nu(y_1^m)-1$ or $2K=i(y_1^m)+\nu(y_1^m)-1$,
hence (\ref{5.65}) holds by Proposition 2.7 and the assumption that $C_{S^1,\,2K}(\Psi_a,\,S^1\cdot u_1^m)\neq 0$.
Thus by Theorem 2.8, we have
\bea &&1=b_{2p}\le M_{2p}=\sum_{1\le i\le 3,\,l\in\N}
\rank C_{S^1,\; 2p}(\Psi_a, \;S^1\cdot u_i^l)
\nn\\=&&\sum_{2\le i\le 3,\,l\in\N}
\rank C_{S^1,\; 2p}(\Psi_a, \;S^1\cdot u_i^l),
\lb{5.66}\eea
for $2p\in[i(y_1^m),\;i(y_1^m)+\nu(y_1^m)-1]\setminus\{2K\}$,
where in  the last equality we have used (\ref{5.61}), (\ref{5.65})
and Proposition 2.3.
Hence by the assumption that $C_{S^1,\,2K}(\Psi_a, \,S^1\cdot u_j^k)\neq 0$
and Claim 3, we have
\bea &&\sum_{i(y_1^m)\le q \le i(y_1^m)+\nu(y_1^m)-1\atop 2\le i\le 3,\,l\in\N}
(-1)^q\rank C_{S^1,\; q}(\Psi_a, \;S^1\cdot u_i^l)
\nn\\=&&\sum_{i(y_1^m)\le 2p \le i(y_1^m)+\nu(y_1^m)-1\atop 2\le i\le 3,\,l\in\N}
\rank C_{S^1,\; 2p}(\Psi_a, \;S^1\cdot u_i^l)\ge \frac{\nu(y_1^m)-1}{2}+1,
\lb{5.67}\eea
In fact, we get the last inequality by counting the number of even integers 
between $i(y_1^m)$ and $i(y_1^m)+\nu(y_1^m)-1$, since by (\ref{5.66}), we count 
the number of $2p$ for  $2p\neq 2K$, this number is $\frac{\nu(y_1^m)-1}{2}$, and there is at least $1$ for $2p=2K$ by the assumption that $C_{S^1,\,2K}(\Psi_a,\,S^1\cdot u_j^k)\neq 0$.

By Proposition 2.3 and Claim 4, we have
\bea &&\sum_{i(y_1^m)\le q \le i(y_1^m)+\nu(y_1^m)-1}
(-1)^q\rank C_{S^1,\; q}(\Psi_a, \;S^1\cdot u_1^m)
\nn\\=&&\sum_{q\in\Z}
(-1)^q\rank C_{S^1,\; q}(\Psi_a, \;S^1\cdot u_1^m)=1.
\lb{5.68}\eea
Combining (\ref{5.64}), (\ref{5.67}) and (\ref{5.68}),
we get a contradiction. This proves Claim 5.

Now we show that there is some  closed characteristic $(\tau_{i_0},\, y_{i_0})$ for 
$i_0\in\{2, 3\}$ such that $i(y_{i_0}, 1)=4$.
In fact, by (\ref{5.15}), (\ref{5.46}), (\ref{5.61}), Proposition 2.3 and
Theorem 2.8, we have \bea &&1=b_0\le M_0=\sum_{1\le i\le 3,\,l\in\N}
\rank C_{S^1,\; 0}(\Psi_a, \;S^1\cdot u_i^l) \nn\\=&&\sum_{2\le i\le
3,\,l\in\N} \rank C_{S^1,\; 0}(\Psi_a, \;S^1\cdot u_i^l),
\lb{5.69}\eea 
where the last equality follows from $ C_{S^1,\; 0}(\Psi_a, \;S^1\cdot u_1^m)=0$
for $m\in\N$. Thus there exist $i_0, l_0$ such that $C_{S^1,\;
0}(\Psi_a, \;S^1\cdot u_{i_0}^{l_0})\neq 0$. Note that
$i(y_{i_0})\ge 0$ since $i(y_{i_0})$ is the Morse index, thus we have
$l_0=1$ by (\ref{4.10}) and Proposition 2.3, and then
$i(y_{i_0})=0$. We may assume $i_0=2$ without loss of generality. Thus
by Theorem 3.6, we have \bea i(y_2,\,\,1)=4.\lb{5.70}\eea 

Up to now, by Claim 5, the problem is transformed to find appropriate $K, m, k\in\N$
such that 
 \bea C_{S^1,\,2K}(\Psi_a, \,S^1\cdot u_1^m)\neq 0,\quad
C_{S^1,\,2K}(\Psi_a, \,S^1\cdot u_2^k)\neq 0.\lb{5.71}\eea Using the precise index iteration formula
 (cf. Theorem 3.7), this is transformed further to a problem in number theory,
 i.e., whether an appropriate integer valued equation has integer solutions.
Thus in the following we separate our proof into several cases according to the possible
cases of $M_2$ appearing in Claim 2.

{\bf Case 1.} {\it The matrix $M_2$ can be connected within $\Omega^0(M_2)$
to $R(\vartheta_1)\diamond R(\vartheta_2)\diamond R(\vartheta_3)$
 with $\frac{\vartheta_1}{\pi},\, \frac{\vartheta_2}{\pi}\notin\Q$
and $\frac{\vartheta_3}{\pi}\in\Q\cap (0,\,2]$, i.e., we handle the case that $M_2^\prime\in
\{I_2, -I_2, R(\vartheta)\}$ with $\frac{\vartheta}{\pi}\in\Q$ together as in Claim 2. }

Write $\frac{\vartheta_3}{2\pi}=\frac{r_1}{s_1}$ with $r_1, s_1\in\N$ and
$(r_1,\, s_1)=1$.  Now we want to find some special $K, m, k\in\N$ such that (\ref{5.71}) holds., thus we suppose $k=ps_1$ and $m\in\{qs-1,\, qs,\,qs+1\}$ for some $p, q\in\N$,
where $s$ is given in Claim 4, i.e., $\frac{r}{s}=\frac{\theta}{2\pi}$.
By (\ref{5.38}), we have
\bea 2K=i(y_2^{ps_1})+\nu(y_2^{ps_1})-1,\lb{5.72}\eea
 since $K(y_2)=s_1$.
By (\ref{5.46}), we have
\bea 2K=i(y_1^{qs\pm1})+\nu(y_1^{qs\pm1})-1,
\quad {\rm if }\quad m=qs\pm1, \lb{5.73}\eea
note that here we have used the fact that $s\ge 2$, i.e., $y_1^{qs\pm1}$ is non-degenerate. If $m=qs$, by (\ref{5.16}), $s|2m_1^\ast$, Propositions 2.3
and  2.6, we have
 \bea &&\rank C_{S^1,\; i(y_1^{qs})+\nu(y_1^{qs})-1-l}(\Psi_a, \;S^1\cdot u_1^{qs})\nn\\
 =&&k_{\nu(u_1^{qs})-1-l}(u_1^{qs})
 =k_{\nu(u_1^s)-1-l}(u_1^s)
 =k_{\nu(u_1^{2m^\ast_1})-1-l}(u_1^{2m^\ast_1})\nn\\
 =&&\rank C_{S^1,\; i(y_1^{2m^\ast_1})+\nu(y_1^{2m^\ast_1})-1-l}(\Psi_a, \;S^1\cdot u_1^{2m^\ast_1})
\nn\\=&&\rank C_{S^1,\; 2T^\ast-2-l}(\Psi_a, \;S^1\cdot u_1^{2m^\ast_1}),
 \lb{5.74}\eea
for any $l\in\Z$. Hence by (\ref{5.19}) and (\ref{5.20}). we have
\bea 2K=i(y_1^{qs})+\nu(y_1^{qs})-1-2,\quad {\rm if }\quad m=qs.\lb{5.75}\eea
By (\ref{5.15}), we have
\bea &&i(y_1^{qs})+\nu(y_1^{qs})-1=4qs+2qr-6+5-1=2q(2s+r)-2,\nn\\
&&i(y_1^{qs-1})+\nu(y_1^{qs-1})-1=4(qs-1)+2qr-6+3-1=2q(2s+r)-8,\nn\\
&&i(y_1^{qs+1})+\nu(y_1^{qs+1})-1=4(qs+1)+2(qr+1)-6+3-1=2q(2s+r)+2.
\qquad\lb{5.76}\eea
By (\ref{5.36}) and (\ref{5.70}), we have
\bea &&i(y_2^{ps_1})+\nu(y_2^{ps_1})-1=2ps_1+2\sum_{j=1}^2
E\left(\frac{ps_1\vartheta_j}{2\pi}\right)+2pr_1-8+3-1\nn\\
=&&2\sum_{j=1}^2E\left(\frac{ps_1\vartheta_j}{2\pi}\right)
+2p(s_1+r_1)-6.\lb{5.77}\eea
By (\ref{5.72})-(\ref{5.77}), we have
\bea 2\sum_{j=1}^2E\left(\frac{ps_1\vartheta_j}{2\pi}\right)
+2p(s_1+r_1)-6=2q(2s+r)+2l,\lb{5.78}\eea
for some $l\in\{-4,\,-2,\,1\}$, where $l=-4$ if $m=qs-1$, $l=-2$ if $m=qs$ and $l=1$ if $m=qs+1$.
Now we suppose further that $p=(2s+r)p^\prime$ and $q=(s_1+r_1)q^\prime$
for some $p^\prime, q^\prime\in\N$, then we have
\bea 2\sum_{j=1}^2E\left(\frac{p^\prime(2s+r)s_1\vartheta_j}{2\pi}\right)
+2p^\prime(2s+r)(s_1+r_1)-6=2q^\prime(s_1+r_1)(2s+r)+2l,\lb{5.79}\eea
for some $ l\in\{-4,\,-2,\,1\}$.
Write $\frac{(2s+r)s_1\vartheta_j}{2\pi}=\alpha_j\notin\Q$
and $N=(s_1+r_1)(2s+r)\ge 3$,
then (\ref{5.79}) is equivalent to find $p^\prime\in\N$ such that
one of the following holds
\bea \sum_{j=1}^2E(p^\prime\alpha_j)\equiv l \;\mod \;N,\qquad
l\in\{-1,\,1,\,4\}.\lb{5.80}\eea
In fact, if we obtain $p^\prime$ and $l$ from (\ref{5.80}), then we can
substitute them into (\ref{5.79}) to get $q^\prime$, and then find a solution of
(\ref{5.78}), and consequently find a solution of (\ref{5.71}).

In order to solve (\ref{5.80}), we have to consider the following sub-cases:

{\bf Sub-case 1.1.} {\it We have $\{1,\,\alpha_1,\,\alpha_2\}$
are linearly independent over $\Q$. }

Clearly, $\{1,\,\frac{\alpha_1}{N},\,\frac{\alpha_2}{N}\}$
are linearly independent over $\Q$ also, then by Theorem 5.5,
 the vectors $\{n\frac{\alpha_1}{N},\,n\frac{\alpha_2}{N}\}_{n\ge 1}$
are uniformly distributed mod one. Hence we can choose $n\in\N$
such that $n\frac{\alpha_1}{N}\equiv \epsilon\;\mod \;1$ and
$n\frac{\alpha_2}{N}=-\epsilon^\prime\;\mod\;1$ for some
$\epsilon,\,\epsilon^\prime\in (0,\,\frac{1}{N})$.
Thus we have
\bea &&\sum_{j=1}^2E(n\alpha_j)=\sum_{j=1}^2E\left(Nn\frac{\alpha_j}{N}\right)
\nn\\\equiv&& E(N\epsilon)+E(-N\epsilon^\prime)\;\mod\;N
\nn\\\equiv&& 1\;\mod\;N,\lb{5.81}
\eea
where the last equality follows from $E(N\epsilon)=1$ and
$E(-N\epsilon^\prime)=0$.
Hence (\ref{5.80}) holds for $p^\prime=n$ and $l=1$. This proves Theorem 1.1 in this case.

{\bf Sub-case 1.2.} {\it We have $\{1,\,\alpha_1,\,\alpha_2\}$
are linearly dependent over $\Q$. }

Since $\alpha_1, \alpha_2\notin\Q$,
we can write $\alpha_2=\frac{r_2}{s_2}\alpha_1+\frac{r_3}{s_3}$
for some $r_2\in\Z\setminus\{0\}$, $r_3\in\Z$, $s_2, s_3\in\N$ with $(r_2,\, s_2)=1$ and
$(r_3,\, s_3)=1$. Note that $\{1,\,s_3\alpha_1\}$
are linearly independent over $\Q$ since $\alpha_1\notin\Q$, then by  Theorem 5.5,
 the vectors $\{ns_3\alpha_1\}_{n\ge 1}$
are uniformly distributed mod one. Hence we can choose $n\in\N$
such that $ns_3\alpha_1=\epsilon \;\mod\; 1$ for some
$\epsilon\in (-1,\,1)$ as we required below.
Then we have
\bea &&\sum_{j=1}^2E(ns_2s_3N\alpha_j)
=E(ns_2s_3N\alpha_1)+E\left(\frac{ns_2s_3Nr_2\alpha_1}{s_2}+ns_2r_3N\right)
\nn\\\equiv&&E\left(s_2N\left(ns_3\alpha_1\right)\right)
+E\left(s_2N\left(\frac{r_2}{s_2}\left(ns_3\alpha_1\right)\right)\right)\;\mod\;N
\nn\\\equiv&&E\left(s_2N\epsilon\right)+E\left(r_2N\epsilon\right)
\;\mod\;N.\lb{5.82}
\eea

We have the following cases:

(i) If $\frac{r_2}{s_2}<0$, then we require
 $\epsilon\in \left(0,\, \min\{\frac{1}{s_2N},\,-\frac{1}{r_2N}\}\right)$,
then (\ref{5.82}) becomes
\bea E\left(s_2N\epsilon\right)+E\left(r_2N\epsilon\right)
\equiv 1\;\mod\;N,\lb{5.83}
\eea
where (\ref{5.83}) follows from $E(s_2N\epsilon)=1$ and
$E(r_2N\epsilon)=0$.
Hence (\ref{5.80}) holds for $p^\prime=ns_2s_3N$ and $l=1$.
This proves Theorem 1.1 in this case.

(ii) If $\frac{r_2}{s_2}>0$, by a permutation of $\alpha_1, \alpha_2$
if necessary, we may assume $\frac{r_2}{s_2}\ge 1$. Then we have:

(ii-a) If $\frac{r_2}{s_2}=1$, then we require
 $\epsilon\in \left(\frac{1}{s_2N},\,\frac{2}{s_2N}\right)$,
then (\ref{5.82}) becomes
\bea 2E\left(s_2N\epsilon\right)
\equiv 4\;\mod\;N,\lb{5.84}
\eea
where (\ref{5.84}) follows from $E(s_2N\epsilon)=2$.
Hence (\ref{5.80}) holds for $p^\prime=ns_2s_3N$ and $l=4$.
This proves Theorem 1.1 in this case.

(ii-b) If $\frac{r_2}{s_2}>1$,
 then we require
 $\epsilon\in \left(\max\{\frac{-1}{s_2N},\,\frac{-2}{r_2N}\},\;\frac{-1}{r_2N}\right)$,
then (\ref{5.82}) becomes
\bea E\left(s_2N\epsilon\right)+E\left(r_2N\epsilon\right)
\equiv -1\;\mod\;N,\lb{5.85}
\eea
where (\ref{5.85}) follows from $E(s_2N\epsilon)=0$ and
$E(r_2N\epsilon)=-1$.
Hence (\ref{5.80}) holds for $p^\prime=ns_2s_3N$ and $l=-1$.
This proves Theorem 1.1 in this case.

{\bf Case 2.} {\it The matrix $M_2$ can be connected within $\Omega^0(M_2)$
to $R(\vartheta_1)\diamond R(\vartheta_2)\diamond R(\vartheta_3)$
 with $\frac{\vartheta_i}{\pi}\notin\Q$ for $1\le i\le 3$. }

As in Case 1, we want to find some $K, m, k\in\N$ such that
(\ref{5.71}) holds. By (\ref{5.36}), (\ref{5.70}) \ and Proposition 2.3, we have
\bea 2K=i(y_2^k)=2k+2\sum_{j=1}^3
E\left(\frac{k\vartheta_j}{2\pi}\right)-8.\lb{5.86}\eea
Thus as in Case 1, we have
\bea 2\sum_{j=1}^3E\left(\frac{k\vartheta_j}{2\pi}\right)
+2k-8=2q(2s+r)+2l,\lb{5.87}\eea
for some $l\in\{-4,\,-2,\,1\}$.
Suppose $k=(2s+r)p^\prime$, then we have
\bea 2\sum_{j=1}^3E\left(\frac{p^\prime(2s+r)\vartheta_j}{2\pi}\right)
+2p^\prime(2s+r)-8=2q(2s+r)+2l,\lb{5.88}\eea
for some $ l\in\{-4,\,-2,\,1\}$.
Write $\frac{(2s+r)\vartheta_j}{2\pi}=\alpha_j\notin\Q$
and $N=2s+r\ge 5$,
then (\ref{5.88}) is equivalent to find $p^\prime\in\N$ such that
one of the following holds
\bea \sum_{j=1}^3E(p^\prime\alpha_j)\equiv l \;\mod \;N,\qquad
l\in\{0,\,2,\,5\}.\lb{5.89}\eea

We have the following sub-cases:

{\bf Sub-case 2.1.} {\it We have $\{1,\,\alpha_1,\,\alpha_2,\,\alpha_3\}$
are linearly independent over $\Q$. }

Clearly, $\{1,\,\frac{\alpha_1}{N},\,\frac{\alpha_2}{N},\,\frac{\alpha_3}{N}\}$
are linearly independent over $\Q$, then by Theorem 5.5,
 the vectors
$\{n\frac{\alpha_1}{N},\,n\frac{\alpha_2}{N},\,n\frac{\alpha_3}{N}\}_{n\ge 1}$
are uniformly distributed mod one. Hence we can choose $n\in\N$
such that $n\frac{\alpha_j}{N}\equiv \epsilon_j\;\mod \;1$
for some $\epsilon_j\in (-\frac{1}{N},\;0)$ and $1\le j\le 3$.
Thus we have
\bea &&\sum_{j=1}^3E(n\alpha_j)=\sum_{j=1}^3E\left(Nn\frac{\alpha_j}{N}\right)
\nn\\\equiv&& \sum_{j=1}^3E(N\epsilon_j)\;\mod\;N
\nn\\\equiv&& 0\;\mod\;N,\lb{5.90}
\eea
where in the last equality we have used $E(N\epsilon_j)=0$ for $1\le j\le 3$.
Hence (\ref{5.89}) holds for $p^\prime=n$ and $l=0$.
This proves Theorem 1.1 in this case.

{\bf Sub-case 2.2.} {\it We have $\{1,\,\alpha_1,\,\alpha_2,\,\alpha_3\}$
are linearly dependent over $\Q$
together with $\dim_{\Q}({\rm span_{\Q}}\{1,\,\alpha_1,\,\alpha_2,\,\alpha_3\})=3$. }

Since $\dim_{\Q}({\rm span_{\Q}}\{1,\,\alpha_1,\,\alpha_2,\,\alpha_3\})=3$,
we may assume $\{1,\,\alpha_1,\,\alpha_2\}$ are linear independent over $\Q$
without loss of generality.
Thus we can write $\alpha_3=\frac{r_1}{s_1}\alpha_1+\frac{r_2}{s_2}\alpha_2+\frac{r_3}{s_3}$
for some $r_i\in\Z$, $s_i\in\N$ with $(r_i,\, s_i)=1$ for $1\le i\le 3$.
Note that $\{1,\,s_3\alpha_1,\,s_3\alpha_2\}$
are linearly independent over $\Q$, then by Theorem 5.5,
 the vectors $\{ns_3\alpha_1,\,ns_3\alpha_2\}_{n\ge 1}$
are uniformly distributed mod one. Hence we can choose $n\in\N$
such that $ns_3\alpha_1=\epsilon_1 \;\mod\; 1$ and
$ns_3\alpha_2=\epsilon_2 \;\mod\; 1$ for some
$\epsilon_1, \epsilon_2\in (-1,\,1)$ as we required below.

Then we have
\bea &&\sum_{j=1}^3E(ns_1s_2s_3N\alpha_j)
\nn\\=&&\sum_{j=1}^2E(ns_1s_2s_3N\alpha_j)
+E\left(\frac{ns_1s_2s_3Nr_1\alpha_1}{s_1}
+\frac{ns_1s_2s_3Nr_2\alpha_2}{s_2}+ns_1s_2r_3N\right)
\nn\\\equiv&&\sum_{j=1}^2E(s_1s_2N(ns_3\alpha_j))
+E\left(s_1s_2N\left(\frac{r_1}{s_1}\left(ns_3\alpha_1\right)
+\frac{r_2}{s_2}\left(ns_3\alpha_2\right)\right)\right)\;\mod\;N
\nn\\\equiv&&\sum_{j=1}^2E(s_1s_2N\epsilon_j)
+E(r_1s_2N\epsilon_1+s_1r_2N\epsilon_2)\;\mod\;N.
\lb{5.91}
\eea

We have the following cases:

(i) If $\frac{r_1}{s_1}\ge0$ and $\frac{r_2}{s_2}\ge0$, we have
$\frac{r_1}{s_1}+\frac{r_2}{s_2}> 0$ since $\alpha_3\notin\Q$,  and then $r_1s_2+r_2s_1>0$.
We require
 $\epsilon_1, \epsilon_2\in \left(\max\{\frac{-1}{(r_1s_2+r_2s_1)N},\,\frac{-1}{s_1s_2N}\},\;0\right)$
then (\ref{5.91}) becomes
\bea \sum_{j=1}^2E(s_1s_2N\epsilon_j)
+E(r_1s_2N\epsilon_1+s_1r_2N\epsilon_2)\equiv 0\;\mod\;N,\lb{5.92}\eea
where we have used $E(s_1s_2N\epsilon_j)=0$ for $1\le j\le 2$
and $E(r_1s_2N\epsilon_1+s_1r_2N\epsilon_2)=0$.
Hence (\ref{5.89}) holds for $p^\prime=ns_1s_2s_3N$ and $l=0$.
This proves Theorem 1.1 in this case.

(ii) If $\frac{r_1}{s_1}\le0$ and $\frac{r_2}{s_2}\le0$, we have
$\frac{r_1}{s_1}+\frac{r_2}{s_2}< 0$ since $\alpha_3\notin\Q$, and then $r_1s_2+r_2s_1<0$.
We require
 $\epsilon_1, \epsilon_2\in \left(0,\;\min\{\frac{-1}{(r_1s_2+r_2s_1)N},\,\frac{1}{s_1s_2N}\}\right)$
then (\ref{5.91}) becomes
\bea \sum_{j=1}^2E(s_1s_2N\epsilon_j)
+E(r_1s_2N\epsilon_1+s_1r_2N\epsilon_2)\equiv 2\;\mod\;N,\lb{5.93}\eea
where we have used $E(s_1s_2N\epsilon_j)=1$ for $1\le j\le 2$
and $E(r_1s_2N\epsilon_1+s_1r_2N\epsilon_2)=0$.
Hence (\ref{5.89}) holds for $p^\prime=ns_1s_2s_3N$ and $l=2$.
This proves Theorem 1.1 in this case.

(iii) If $\frac{r_1}{s_1}>0$ and $\frac{r_2}{s_2}<0$,
We require
 $\epsilon_1, \epsilon_2\in \left(0,\;\frac{1}{s_1s_2N}\right)$
 satisfies $\frac{-1}{s_1s_2N}<\frac{r_1}{s_1}\epsilon_1+\frac{r_2}{s_2}\epsilon_2<0$.
In fact, we first choose $\epsilon_1\in\left(0,\,\min\{\frac{-r_2}{2r_1s_2^2N},\,\frac{1}{s_1s_2N}\}\right)$
 sufficiently close to $0$, then we can choose
 $\epsilon_2\in\left(\frac{-s_2r_1}{s_1r_2}\epsilon_1,\;\min\{\frac{1}{s_1s_2N},\,\frac{-1}{s_1r_2N}+\frac{-s_2r_1}{s_1r_2}\epsilon_1\}\right)$,
 then the above inequality holds.
Hence (\ref{5.91}) becomes
\bea \sum_{j=1}^2E(s_1s_2N\epsilon_j)
+E(r_1s_2N\epsilon_1+s_1r_2N\epsilon_2)\equiv 2\;\mod\;N.\lb{5.94}\eea
where we have used $E(s_1s_2N\epsilon_j)=1$ for $1\le j\le 2$
and $E(r_1s_2N\epsilon_1+s_1r_2N\epsilon_2)=0$.
Hence (\ref{5.89}) holds for $p^\prime=ns_1s_2s_3N$ and $l=2$.
This proves Theorem 1.1 in this case.

(iv) Similarly, if $\frac{r_1}{s_1}<0$ and $\frac{r_2}{s_2}>0$, Theorem 1.1 holds.
This proves Theorem 1.1 in Sub-case 2.2.

{\bf Sub-case 2.3.} {\it We have $\{1,\,\alpha_1,\,\alpha_2,\,\alpha_3\}$
are linearly dependent over $\Q$
together with $\dim_{\Q}({\rm span_{\Q}}\{1,\,\alpha_1,\,\alpha_2,\,\alpha_3\})=2$. }

In this case, we can write
$\alpha_2=\frac{r_1}{s_1}\alpha_1+\frac{r_2}{s_2}$ and
$\alpha_3=\frac{r_3}{s_3}\alpha_1+\frac{r_4}{s_4}$
for some $r_1, r_3\in\Z\setminus\{0\}$, $r_2, r_4\in\Z$, $s_i\in\N$ with $(r_i,\, s_i)=1$
for $1\le i\le 4$. Note that $\{1,\,s_2s_4\alpha_1\}$
are linearly independent over $\Q$ since $\alpha_1\notin\Q$, then by Theorem 5.5,
  the vectors $\{ns_2s_4\alpha_1\}_{n\ge 1}$
are uniformly distributed mod one. Hence we can choose $n\in\N$
such that $ns_2s_4\alpha_1=\epsilon \;\mod\; 1$  for some
$\epsilon\in (-1,\,1)$ as we required below.
Then we have
\bea &&\sum_{j=1}^3E(ns_1s_2s_3s_4N\alpha_j)
\nn\\=&&E(ns_1s_2s_3s_4N\alpha_1)
+E\left(\frac{ns_1s_2s_3s_4Nr_1\alpha_1}{s_1}
+ns_1r_2s_3s_4N\right)\nn\\
&&+E\left(\frac{ns_1s_2s_3s_4Nr_3\alpha_1}{s_3}
+ns_1s_2s_3r_4N\right)\nn\\
\nn\\\equiv&&E(s_1s_3N(ns_2s_4\alpha_1))
+E\left(s_1s_3N\left(\frac{r_1}{s_1}\left(ns_2s_4\alpha_1\right)
\right)\right)\nn\\
&&+E\left(s_1s_3N\left(\frac{r_3}{s_3}\left(ns_2s_4\alpha_1\right)
\right)\right)\;\mod\;N\nn\\
\equiv&&E(s_1s_3N\epsilon)+E(r_1s_3N\epsilon)+E(s_1r_3N\epsilon)
\;\mod\;N.
\lb{5.95}
\eea

We have the following cases:

(i) If $\frac{r_1}{s_1}>0$ and $\frac{r_3}{s_3}>0$,
we require
 $\epsilon\in \left(\max\{\frac{-1}{r_1s_3N},\,\frac{-1}{r_3s_1N},\,\frac{-1}{s_1s_3N}\},\;0\right)$
then (\ref{5.95}) becomes
\bea E(s_1s_3N\epsilon)+E(r_1s_3N\epsilon)+E(s_1r_3N\epsilon)\equiv 0\;\mod\;N,\lb{5.96}\eea
where we have used $E(s_1s_3N\epsilon)=0$, $E(r_1s_3N\epsilon)=0$
and $E(s_1r_3N\epsilon)=0$.
Hence (\ref{5.89}) holds for $p^\prime=ns_1s_2s_3s_4N$ and $l=0$.
This proves Theorem 1.1 in this case.

(ii) If $\frac{r_1}{s_1}<0$ and $\frac{r_3}{s_3}<0$,
we require
 $\epsilon\in \left(\max\{\frac{1}{r_1s_3N},\,\frac{1}{r_3s_1N},\,\frac{-1}{s_1s_3N}\},\;0\right)$
then (\ref{5.95}) becomes
\bea E(s_1s_3N\epsilon)+E(r_1s_3N\epsilon)+E(s_1r_3N\epsilon)\equiv 2\;\mod\;N,\lb{5.97}\eea
where we have used $E(s_1s_3N\epsilon)=0$, $E(r_1s_3N\epsilon)=1$
and $E(s_1r_3N\epsilon)=1$.
Hence (\ref{5.89}) holds for $p^\prime=ns_1s_2s_3s_4N$ and $l=2$.
This proves Theorem 1.1 in this case.

(iii) If $\frac{r_1}{s_1}>0$ and $\frac{r_3}{s_3}<0$,
we require
 $\epsilon\in \left(0,\;\min\{\frac{1}{r_1s_3N},\,\frac{-1}{r_3s_1N},\,\frac{1}{s_1s_3N}\}\right)$
then (\ref{5.95}) becomes
\bea E(s_1s_3N\epsilon)+E(r_1s_3N\epsilon)+E(s_1r_3N\epsilon)\equiv 2\;\mod\;N,\lb{5.98}\eea
where we have used $E(s_1s_3N\epsilon)=1$, $E(r_1s_3N\epsilon)=1$
and $E(s_1r_3N\epsilon)=0$.
Hence (\ref{5.89}) holds for $p^\prime=ns_1s_2s_3s_4N$ and $l=2$.
This proves Theorem 1.1 in this case.

(iv) Similarly, if $\frac{r_1}{s_1}<0$ and $\frac{r_3}{s_3}>0$, Theorem 1.1 holds.
This proves Theorem 1.1 in Sub-case 2.3.

{\bf Case 3.} {\it The matrix $M_2$ can be connected within $\Omega^0(M_2)$
to $R(\vartheta_1)\diamond R(\vartheta_2)\diamond N_1(-1,\,1)$
 with $\frac{\vartheta_1}{\pi},\, \frac{\vartheta_2}{\pi}\notin\Q$. }

By (\ref{5.39}) and (\ref{5.70}), we have
\bea &&i(y_2^{2p})+\nu(y_2^{2p})-1=6p+2\sum_{j=1}^2
E\left(\frac{2p\vartheta_j}{2\pi}\right)-7+2-1\nn\\
=&&2\sum_{j=1}^2E\left(\frac{2p\vartheta_j}{2\pi}\right)
+6p-6.\lb{5.99}\eea
Hence as in Case 1, we have
\bea 2\sum_{j=1}^2E\left(\frac{2p\vartheta_j}{2\pi}\right)
+6p-6=2q(2s+r)+2l,\lb{5.100}\eea
for some $l\in\{-4,\,-2,\,1\}$.
Suppose $p=(2s+r)p^\prime$ and $q=3q^\prime$, then we have
\bea 2\sum_{j=1}^2E\left(\frac{p^\prime(2s+r)2\vartheta_j}{2\pi}\right)
+6p^\prime(2s+r)-6=6q^\prime(2s+r)+2l,\lb{5.101}\eea
for some $ l\in\{-4,\,-2,\,1\}$.
Write $\frac{(2s+r)2\vartheta_j}{2\pi}=\alpha_j\notin\Q$
and $N=3(2s+r)$,
then (\ref{5.101}) is equivalent to find $p^\prime\in\N$ such that
one of the following holds
\bea \sum_{j=1}^2E(p^\prime\alpha_j)\equiv l \;\mod \;N,\qquad
l\in\{-1,\,1,\,4\}.\lb{5.102}\eea
Then by the same proof as in Case 1, Theorem 1.1 holds in this case.

{\bf Case 4.} {\it The matrix $M_2$ can be connected within $\Omega^0(M_2)$
to $R(\vartheta_1)\diamond R(\vartheta_2)\diamond N_1(1,\,-1)$
 with $\frac{\vartheta_1}{\pi},\, \frac{\vartheta_2}{\pi}\notin\Q$. }

By (\ref{5.42}), we have $i(y_2,\, 1)\in 2\N-1$, this contradict to
(\ref{5.70}), so this case can not happen. Hence Theorem 1.1 holds
in this case.

The proof of Lemma 5.6 is complete.\hfill\hb

{\bf Lemma 5.7.} {\it  If  $(\tau_1, y_1)$ belongs to Case A in \S4 and  the matrix $M_1$ can be connected within $\Omega^0(M_1)$ to
 $N_1(1,\,-1)^{\diamond 2}\diamond N_1(-1,\,b)$ with
 $b=0, -1$, then we have $^\#\T(\Sg)\ge 4$.}

{\bf Proof.} Note that the case $N_1(-1,\,0)=-I_2$
has already been proved in Lemma 5.6 since $-I_2=R(\pi)$.
While the proof for $R(\pi)$ also apply to the case
$N_1(-1,\,-1)$ since the index iteration formulae for
$-I_2$ and $N_1(-1,\,-1)$ are the same by Theorem 3.7,
the only difference is their nullities for $2m$-th iteration,
but this will not affect our argument. In fact we replace $\nu(y_1^{2m})$ by
$\nu(y_1^{2m})+1$ in all the corresponding formulae, then the proof goes as before.
 This proves Lemma 5.7.
\hfill\hb

{\bf Lemma 5.8.} {\it  If  $(\tau_1, y_1)$ belongs to Case A in \S4 and the matrix $M_1$ can be connected within $\Omega^0(M_1)$ to
 $N_1(1,\,-1)^{\diamond 2}\diamond N_1(-1,\,1)$,
 then we have $^\#\T(\Sg)\ge 4$.}

{\bf Proof. } As in Lemma 5.1, suppose $(T,\,m_1,\,m_2,\,m_3)$ and $(j_k,\,l_{j_k})$
 satisfy (\ref{4.15})-(\ref{4.21}). As mentioned in
Case A, we have  $i(y_1,\,1)=4$, thus by Theorems 3.6,  3.7
and (\ref{4.3}), we have
\bea i(y_1^m)&=&m(i(y_1,\,1)+1)-1-4=5m-5, \nn\\
\nu(y_1^m)&=&3+\frac{1+(-1)^m}{2},\qquad m\in\N.\lb{5.103}\eea
 By (\ref{4.18}), we have
$i(y_1^{2m_1-1})+\nu(y_1^{2m_1-1})-1=2T-8$. Hence we have
$i(y_1^{2m_1})=2T-5$ and $i(y_1^{2m_1})+\nu(y_1^{2m_1})-1=2T-2$. Hence by Propositions 2.3
and  2.6 with $K(u_1)=2$, we have
 \bea &&\rank C_{S^1,\; i(y_1^{2m})+\nu(y_1^{2m})-1-l}(\Psi_a, \;S^1\cdot u_1^{2m})\nn\\
 =&&k_{\nu(u_1^{2m})-1-l}(u_1^{2m})
 =k_{\nu(u_1^2)-1-l}(u_1^2)
 =k_{\nu(u_1^{2m^\ast_1})-1-l}(u_1^{2m^\ast_1})\nn\\
 =&&\rank C_{S^1,\; 2T^\ast-2-l}(\Psi_a, \;S^1\cdot u_1^{2m_1^\ast}),
 \lb{5.104}\eea
 and \bea&&\rank C_{S^1,\; i(y_1^{2m-1})+\nu(y_1^{2m-1})-1-l}(\Psi_a, \;S^1\cdot u_1^{2m-1})\nn\\
 =&&k_{\nu(u_1^{2m-1})-1-l}(u_1^{2m-1})
 =k_{\nu(u_1)-1-l}(u_1)
 =k_{\nu(u_1^{2m^\ast_1-1})-1-l}(u_1^{2m^\ast_1-1})\nn\\
 =&&\rank C_{S^1,\; 2T^\ast-8-l}(\Psi_a, \;S^1\cdot u_1^{2m_1^\ast-1}),\lb{5.105}\eea
for any $m\in\N$ and $l\in\Z$.

By Proposition 2.3 and (\ref{4.21}), we have $ C_{S^1,\; 2T^\ast-2-l}(\Psi_a, \;S^1\cdot u_1^{2m_1^\ast})\neq 0$ for some $l\in\{0, 2\}$. Then we have the following two cases:

(i) If  $C_{S^1,\; 2T^\ast-2}(\Psi_a, \;S^1\cdot u_1^{2m_1^\ast})\neq 0$,
then we have $C_{S^1,\; 2T^\ast-2-l}(\Psi_a, \;S^1\cdot u_1^{2m_1^\ast})=0$
for $l\neq 0$ by Proposition 2.7 since $u_1^{2m_1^\ast}$ is a local maximum in 
the local characteristic manifold $W(u_1^{2m_1^\ast})$. This implies
$C_{S^1,\; 2T-2-l}(\Psi_a, \;S^1\cdot u_1^{2m_1})=0$
for $l\neq 0$ by (\ref{5.104}). Hence by (\ref{4.21}), we have $c_T=\Phi(u_1^{2m_1})$,
and then $c_{T+1-\xi_T^{-1}(i)}=\Phi(u_i^{2m_i})$
for  $i=2, 3$ and $\xi_T^{-1}(i)\in\{2, 3\}$.
Thus  we have $\Phi(u_1^{2m_1})>\Phi(u_i^{2m_i})$
for $i=2, 3$.   This contradict to Proposition 4.5
and proves the lemma in this case.

(ii) It remains to consider the case
\bea C_{S^1,\; 2T^\ast-4}(\Psi_a, \;S^1\cdot u_1^{2m_1^\ast})\neq 0,\lb{5.106}\eea
then $C_{S^1,\; 2T^\ast-2}(\Psi_a, \;S^1\cdot u_1^{2m_1^\ast})=0$
and $C_{S^1,\; 2T^\ast-6}(\Psi_a, \;S^1\cdot u_1^{2m_1^\ast})=0$
by Propositions 2.3 and  2.7.
This implies
$C_{S^1,\; 2T-2}(\Psi_a, \;S^1\cdot u_1^{2m_1})=0$
and $C_{S^1,\; 2T-6}(\Psi_a, \;S^1\cdot u_1^{2m_1})=0$
 by (\ref{5.104}). Hence we have $c_{T-1}=\Phi(u_1^{2m_1})$
by (\ref{4.21}), and then we have $c_T=\Phi(u_{\xi_T(1)}^{2m_{\xi_T(1)}})$
and $c_{T-2}=\Phi(u_{\xi_T(3)}^{2m_{\xi_T(3)}})$ for $\xi_T(1), \xi_T(3)\in\{2, 3\}$
and $\xi_T(1)\neq \xi_T(3)$.

Note that by(\ref{5.103}), we have
\bea i(y_1^{2m})=10m-5=10m-8+3=i(y_1^{2m-1})+\nu(y_1^{2m-1})-1+3,
\quad \forall m\in\N.\lb{5.107}\eea

By the same argument as in Lemma 5.6, Claims 1-4 and Claim 5 in \S 4 for $m\in 2\N-1$
hold in this case.

We remark that  Claim 5 in \S 4 for  $m\in2\N$ also holds in this case. In fact,
we can modify the proof of Claim 5 in \S 4 as the following:
By (\ref{5.61}) and (\ref{5.107}) we have
\bea C_{S^1,\, i(u_1^{2m})-2}(\Psi_a,\, S^1\cdot u_1^l)=0,
\quad C_{S^1,\, i(u_1^{2m})+\nu(u_1^{2m})}(\Psi_a,\, S^1\cdot u_1^l)=0,\lb{5.108}\eea
for any $l\in\N$ by Proposition 2.3.
By (\ref{5.103}), we have $i(y_1^{2m})-1,\, i(y_1^{2m})+\nu(y_1^{2m})-1\in 2\N$,
hence by Claim 3, (\ref{5.108}) and Theorem 2.8, we have
\bea &&M_{i(y_1^{2m})-2}=\sum_{1\le i\le 3,\,l\in\N}\rank C_{S^1,\; i(y_1^{2m})-2}(\Psi_a, \;S^1\cdot u_i^l)=0=b_{i(y_1^{2m})-2},\nn\\
&&M_{i(y_1^{2m})+\nu(y_1^{2m})}=\sum_{1\le i\le 3,\,l\in\N}\rank
C_{S^1,\; i(y_1^{2m})+\nu(y_1^{2m})}(\Psi_a, \;S^1\cdot
u_i^l)=0=b_{i(y_1^{2m})+\nu(y_1^{2m})}.\qquad \lb{5.109}\eea Thus as in
Claim 4 in \S 4, by Theorem 2.8, we have \bea
&&\frac{\nu(y^{2m}_1)}{2}+1=\sum_{i(y_1^{2m})-1\le q \le
i(y_1^{2m})+\nu(y_1^{2m})-1}(-1)^qb_q\nn\\
 =&&\sum_{i(y_1^{2m})-1\le q \le
i(y_1^{2m})+\nu(y_1^{2m})-1}(-1)^qM_q \nn\\=&&\sum_{i(y_1^{2m})-1\le
q \le i(y_1^{2m})+\nu(y_1^{2m})-1\atop 1\le i\le 3,\,l\in\N}
(-1)^q\rank C_{S^1,\; q}(\Psi_a, \;S^1\cdot u_i^l)
\nn\\=&&\sum_{i(y_1^{2m})-1\le q \le i(y_1^{2m})+\nu(y_1^{2m})-1}
(-1)^q\rank C_{S^1,\; q}(\Psi_a, \;S^1\cdot u_1^{2m})
\nn\\&&\qquad+\sum_{i(y_1^{2m})-1\le q \le
i(y_1^{2m})+\nu(y_1^{2m})-1\atop 2\le i\le 3,\,l\in\N} (-1)^q\rank
C_{S^1,\; q}(\Psi_a, \;S^1\cdot u_i^l), \lb{5.110}\eea where the
last equality follows from (\ref{5.61}), (\ref{5.107}) and
Proposition 2.3. By (\ref{5.103}), we have $\nu(y_1^{2m})=4$ and
\bea C_{S^1,\; 2p}(\Psi_a, \;S^1\cdot u_1^{2m})=0,\quad \forall
2p\in[i(y_1^{2m})-1,\;i(y_1^{2m})+\nu(y_1^{2m})-1]\setminus\{2K\}.
\lb{5.111}\eea In fact, either $2K=i(y_1^{2m})+\nu(y_1^{2m})-3$ or
$2K=i(y_1^{2m})+\nu(y_1^{2m})-1$ holds. Hence (\ref{5.111}) holds by
Propositions 2.3 and  2.7. Thus by Theorem 2.8, we have \bea &&1=b_{2p}\le
M_{2p}=\sum_{1\le i\le 3,\,l\in\N} \rank C_{S^1,\; 2p}(\Psi_a,
\;S^1\cdot u_i^l) \nn\\=&&\sum_{2\le i\le 3,\,l\in\N} \rank
C_{S^1,\; 2p}(\Psi_a, \;S^1\cdot u_i^l), \lb{5.112}\eea for
$2p\in[i(y_1^{2m})-1,\;i(y_1^{2m})+\nu(y_1^{2m})-1]\setminus\{2K\}$,
where in the last equality we have used (\ref{5.61}), (\ref{5.111})
and Proposition 2.3. Hence as in Claim 5 in \S 4,  by the assumption that
$C_{S^1,\,2K}(\Psi_a, \,S^1\cdot u_j^k)\neq 0$ and Claim 3, we have
\bea &&\sum_{i(y_1^{2m})-1\le q \le i(y_1^{2m})+\nu(y_1^{2m})-1\atop
2\le i\le 3,\,l\in\N} (-1)^q\rank C_{S^1,\; q}(\Psi_a, \;S^1\cdot
u_i^l) \nn\\=&&\sum_{i(y_1^{2m})-1\le 2p \le
i(y_1^{2m})+\nu(y_1^{2m})-1\atop 2\le i\le 3,\,l\in\N} \rank
C_{S^1,\; 2p}(\Psi_a, \;S^1\cdot u_i^l)\ge
\frac{\nu(y_1^{2m})}{2}+1. \lb{5.113}\eea By Proposition 2.3 and
Claim 4, we have \bea &&\sum_{i(y_1^{2m})-1\le q \le
i(y_1^{2m})+\nu(y_1^{2m})-1} (-1)^q\rank C_{S^1,\; q}(\Psi_a,
\;S^1\cdot u_1^{2m}) \nn\\=&&\sum_{q\in\Z} (-1)^q\rank C_{S^1,\;
q}(\Psi_a, \;S^1\cdot u_1^{2m})=1. \lb{5.114}\eea Combining
(\ref{5.110}), (\ref{5.113}) and (\ref{5.114}), we get a
contradiction. This proves Claim 5.

Thus as in Lemma 5.6,  we use Claim 5 to get a contradiction,
i.e., we want to find some $K, m, k\in\N$ such that
 \bea C_{S^1,\,2K}(\Psi_a, \,S^1\cdot u_1^m)\neq 0,\quad
C_{S^1,\,2K}(\Psi_a, \,S^1\cdot u_2^k)\neq 0.\lb{5.115}\eea
Suppose  $m\in\{2q-1,\, 2q,\,2q+1\}$ for $q\in\N$,
Then by Proposition 2.7, (\ref{5.1}) and (\ref{5.105}) or (\ref{5.104}) and (\ref{5.106}), we have
$2K=i(y_1^{2q})+\nu(y_1^{2q})-1-2$ or
$2K=i(y_1^{2q\pm1})+\nu(y_1^{2q\pm1})-1$.

By (\ref{5.103}), we have
\bea &&i(y_1^{2q})+\nu(y_1^{2q})-1=10q-5+4-1=10q-2,\nn\\
&&i(y_1^{2q-1})+\nu(y_1^{2q-1})-1=5(2q-1)-5+3-1=10q-8,\nn\\
&&i(y_1^{2q+1})+\nu(y_1^{2q+1})-1=5(2q+1)-5+3-1=10q+2.
\lb{5.115}\eea
Thus by the same argument as in Lemma 5.6, we can
transform (\ref{5.115}) to an appropriate  integer valued equation and
use Theorem 5.5 to get solutions.
 This proves Lemma 5.8.
\hfill\hb

{\bf Lemma 5.9.} {\it If  $(\tau_1, y_1)$ belongs to Case A in \S4 and the matrix $M_1$ can be connected within $\Omega^0(M_1)$ to
 $N_1(1,\,-1)^{\diamond 2}\diamond I_2$, then we have $^\#\T(\Sg)\ge 4$.}

{\bf Proof.}  As in Lemma 5.1, suppose $(T,\,m_1,\,m_2,\,m_3)$ and $(j_k,\,l_{j_k})$
satisfy (\ref{4.15})-(\ref{4.21}). As mentioned in
Case A above, we have  $i(y_1,\,1)=4$, thus by Theorems 3.6,  3.7
and (\ref{4.3}), we have $i(y_1^m)=m(i(y_1,\,1)+1+1)-1-1-4=6m-6$ and
$\nu(y_1^m)=5$ for $m\in\N$. By (\ref{4.18}), we have
$i(y_1^{2m_1-1})+\nu(y_1^{2m_1-1})-1=2T-8$. Hence we have
$i(y_1^{2m_1})+\nu(y_1^{2m_1})-1=2T-2$. Hence by Propositions 2.3
and  2,6, we have $K(y_1)=1$ and
 \bea &&\rank C_{S^1,\; 2T-2}(\Psi_a, \;S^1\cdot u_1^{2m_1})\nn\\
 =&&k_{\nu(u_1^{2m_1})-1}(u_1^{2m_1})
 =k_{\nu(u_1)-1}(u_1)=k_{\nu(u_1^{2m^\ast_1-1})-1}(u_1^{2m^\ast_1-1})\nn\\
 =&&\rank C_{S^1,\; 2T^\ast-8}(\Psi_a, \;S^1\cdot u_1^{2m^\ast_1-1})\neq 0,\lb{5.117}\eea
where the last equality follows from (\ref{5.1}).
Hence
\bea\rank C_{S^1,\; 2T-2-l}(\Psi_a, \;S^1\cdot u_1^{2m_1})
=k_{\nu(u_1^{2m_1})-1-l}(u_1^{2m_1})=0
\lb{5.118}\eea
for $l\neq 0$ by Proposition 2.7.
 Hence by (\ref{4.21}), we have $c_T=\Phi(u_1^{2m_1})$,
and then $c_{T+1-\xi_T^{-1}(i)}=\Phi(u_i^{2m_i})$
for  $i=2, 3$ and $\xi_T^{-1}(i)\in\{2, 3\}$.
Thus  we have $\Phi(u_1^{2m_1})>\Phi(u_i^{2m_i})$
for $i=2, 3$.   This contradict to Proposition 4.5
and proves the lemma.
\hfill\hb

{\bf Proof of Theorem 1.1.} Combining  Lemmas 5.1-5.3 and 5.6-5.9, we have $^\#T(\Sg)\ge 4$
for all the possible cases.
This proves Theorem 1.1.\hfill\hb

\noindent {\bf Acknowledgements.} I would like to sincerely thank my
Ph. D. thesis advisor, Professor Yiming Long, for introducing me to Hamiltonian
dynamics and for his valuable help and encouragement during my research. 
I would like to sincerely thank him for his valuable discussions  and suggestions during 
the write of this paper. 
I would like to say how enjoyable it is to work with him. I would like to sincerely thank Dr Hui Liu for
pointing out that I missed one possible case in the first version of this paper.

\bibliographystyle{abbrv}

\medskip

\end{document}